\documentclass[12pt]{amsart}
\newtheorem{theorem}{Theorem}
\newtheorem{lemma}{Lemma}
\newtheorem{proposition}{Proposition}
\newtheorem{corollary}{Corollary}
\begin{document}
\title{}
\author{}
\centerline{\Large Jacobi fields of the Tanaka-Webster} \vskip
0.2in \centerline{\Large connection on Sasakian manifolds} \vskip
0.4in \centerline{\large Elisabetta Barletta\footnote{Authors'
address: Universit\`a degli Studi della Basilicata, Dipartimento
di Matematica, Campus Macchia Romana, 85100 Potenza, Italy,
e-mail: {\tt barletta@unibas.it}, {\tt dragomir@unibas.it} }
\hspace{1cm} Sorin Dragomir}
\begin{abstract} We build a variational theory of geodesics of the
Ta\-na\-ka-Webster connection $\nabla$ on a strictly pseudoconvex
CR manifold $M$. Given a contact form $\theta$ on $M$ such that
$(M , \theta )$ has nonpositive pseudohermitian sectional
curvature $(k_\theta (\sigma ) \leq 0)$ we show that $(M, \theta
)$ has no horizontally conjugate points. Moreover, if $(M, \theta
)$ is a Sasakian manifold such that $k_\theta (\sigma ) \geq k_0 >
0$ then we show that the distance between any two consecutive
conjugate points on a lengthy geodesic of $\nabla$ is at most $\pi
/(2 \sqrt{k_0})$. We obtain the first and second variation
formulae for the Riemannian length of a curve in $M$ and show that
in general geodesics of $\nabla$ admitting horizontally conjugate
points do not realize the Riemannian distance.
\end{abstract}
\maketitle
\section{Introduction}
Sasakian manifolds possess a rich geometric structure (cf.
\cite{kn:Bla}, p. 73-80) and are perhaps the closest odd
dimensional analog of K\"ahlerian manifolds. In particular the
concept of holomorphic sectional curvature admits a Sasakian
counterpart, the so called $\varphi$-sectional curvature $H(X)$
(cf. \cite{kn:Bla}, p. 94) and it is a natural problem (as well as
in K\"ahlerian geometry, cf. e.g. \cite{kn:KoNo}, p. 171, and p.
368-373) to investigate how restrictions on $H(X)$ influence upon
the topology of the manifold. An array of findings in this
direction are described in \cite{kn:Bla}, p. 77-80. For instance,
by a result of M. Harada, \cite{kn:Har1}, for any compact regular
Sasakian manifold $M$ satisfying the inequality $h > k^2$ the
fundamental group $\pi_1 (M)$ is cyclic. Here $h = \inf \{ H(X) :
X \in T_x (M), \; \| X \| = 1, \; x \in M \}$ and it is also
assumed that the least upper bound of the sectional curvature of
$M$ is $1/k^2$. Moreover, if additionally $M$ has minimal diameter
$\pi$ then $M$ is isometric to the standard sphere $S^{2n+1}$, cf.
\cite{kn:Har2}, p. 200.
\par
In the present paper we embrace a different point of view, that of
pseudohermitian geometry (cf. \cite{kn:Web}). To describe it we
need to introduce a few basic objects (cf. \cite{kn:Bla}, p.
19-28). Let $M$ be a $(2n+1)$-dimensional $C^\infty$ manifold and
$(\varphi , \xi , \eta , g)$ a {\em contact metric structure} i.e.
$\varphi$ is an endomorphism of the tangent bundle, $\xi$ is a
tangent vector field, $\eta$ is a differential $1$-form, and $g$
is a Riemannian metric on $M$ such that
\[ \varphi^2 = - I + \eta \otimes \xi , \;\;\; \varphi (\xi ) = 0,
\;\;\; \eta (\xi ) = 1, \]
\[ g(\varphi X , \varphi Y) = g(X,Y) - \eta (X) \eta (Y), \;\;\;
X,Y \in T(M), \] and $\Omega = d \eta$ (the {\em contact
condition}) where $\Omega (X,Y) = g(X , \varphi Y)$. Any contact
Riemannian manifold $(M, (\varphi , \xi , \eta , g))$ admits a
natural almost CR structure
\[ T_{1,0}(M) = \{ X - i J X : X \in {\rm Ker}(\eta )\} \]
($i = \sqrt{-1}$) i.e. it satisfies (\ref{e:int1}) below. By a
result of S. Ianu\c{s}, \cite{kn:Ian}, if $(\varphi , \xi , \eta
)$ is {\em normal} (i.e. $[\varphi , \varphi ] + 2 (d \eta )
\otimes \xi = 0$) then $T_{1,0}(M)$ is integrable, i.e. it obeys
to (\ref{e:int2}) in Section 2. Cf. \cite{kn:Bla}, p. 57-61, for
the geometric interpretation of normality, as related to the
classical embeddability theorem for real analytic CR structures
(cf. \cite{kn:AnHi}). Integrability of $T_{1,0}(M)$ is required in
the construction of the Tanaka-Webster connection of $(M , \eta
)$, cf. \cite{kn:Tan}, \cite{kn:Web} and definitions in Section 2
(although many results in pseudohermitian geometry are known to
carry over to arbitrary contact Riemannian manifolds, cf.
\cite{kn:Tann} and more recently \cite{kn:BaDr2}, \cite{kn:BlDr}).
A manifold carrying a contact metric structure $(\varphi , \xi ,
\eta , g)$ whose underlying contact structure $(\varphi , \xi ,
\eta )$ is normal is a {\em Sasakian manifold} (and $g$ is a {\em
Sasakian metric}). The main tool in the Riemannian approach to the
study of Sasakian geometry is the availability of a variational
theory of geodesics of the Levi-Civita connection of $(M , g)$
(cf. e.g. \cite{kn:Har2}, 194-197). In this paper we start the
elaboration of a similar theory regarding the geodesics of the
Tanaka-Webster connection $\nabla$ of $(M , \eta )$ and give a few
applications (cf. Theorems \ref{t:conj1}-\ref{t:conj2} and
\ref{t:14} below). Our motivation is twofold. First, we aim to
study the topology of Sasakian manifolds under restrictions on the
curvature of $\nabla$ and conjecture that Carnot-Carath\'eodory
complete Sasakian manifolds whose pseudohermitian Ricci tensor
$\rho$ satisfies $\rho (X,X) \geq (2n-1) k_0 \| X \|^2$ for some
$k_0 > 0$ and any $X \in {\rm Ker}(\eta )$ must be compact.
Second, the relationship between the sub-Riemannian geodesics of
the sub-Riemannain manifold $(M , {\rm Ker}(\eta ), g)$ and the
geodesics of $\nabla$ (emphasized by our Corollary
\ref{c:relation}) together with R.S. Strichartz's arguments (cf.
\cite{kn:Str}, p. 245 and 261-262) clearly indicates that a
variational theory of geodesics of $\nabla$ is the key requirement
in bringing results such as those in \cite{kn:Str2} or
\cite{kn:Oba} into the realm of subelliptic theory. In
\cite{kn:Ba} one obtains a pseudohermitian version of the Bochner
formula (cf. e.g. \cite{kn:BGM}, p. 131) implying a lower bound on
the first nonzero eigenvalue $\lambda_1$ of the sublaplacian
$\Delta_b$ of a compact Sasakian manifold
\begin{equation}
- \lambda_1 \geq 2nk/(2n-1) \label{e:Lic}
\end{equation}
(a CR analog to the Lichnerowicz theorem, \cite{kn:Lic}). It is
likely that a theory of geodesics of $\nabla$ may be employed to
show that equality in (\ref{e:Lic}) implies that $M$ is CR
isomorphic to a sphere $S^{2n+1}$ (the CR analog to Obata's
result, \cite{kn:Oba}).

\vskip 0.1in {\small {\bf Acknowledgements}. The Authors are
grateful to the anonymous Referee who pointed out a few errors in
the original version of the manuscript. The Authors acknowledge
support from INdAM (Italy) within the interdisciplinary project
{\em Nonlinear subelliptic equations of variational origin in
contact geometry}.}

\section{Sub-Riemannian geometry on CR manifolds}
Let $M$ be an orientable $(2n+1)$-dimensional $C^\infty$ manifold.
A {\em CR structure} on $M$ is a complex distribution $T_{1,0}(M)
\subset T(M) \otimes {\mathbb C}$, of complex rank $n$, such that
\begin{equation}
T_{1,0}(M) \cap T_{0,1}(M) = (0)
\label{e:int1}
\end{equation}
and
\begin{equation}
\label{e:int2}
Z, W \in T_{1,0}(M) \Longrightarrow [Z, W] \in T_{1,0}(M)
\end{equation}
(the {\em formal integrability property}). Here $T_{0,1}(M) =
\overline{T_{1,0}(M)}$ (overbars denote complex conjugates). The
integer $n$ is the {\em CR dimension}. The pair $(M, T_{1,0}(M))$
is a {\em CR manifold} (of {\em hypersurface type}). Let $H(M) =
{\rm Re}\{ T_{1,0}(M) \oplus T_{0,1}(M) \}$ be the {\em Levi
distribution}. It carries the complex structure $J : H(M) \to
H(M)$ given by $J (Z + \overline{Z}) = i (Z - \overline{Z})$ ($i =
\sqrt{-1}$). Let $H(M)^\bot \subset T^* (M)$ the conormal bundle,
i.e. $H(M)^\bot_x = \{ \omega \in T_x^* (M) : {\rm Ker}(\omega )
\supseteq H(M)_x \}$, $x \in M$. A {\em pseudohermitian structure}
on $M$ is a globally defined nowhere zero cross-section $\theta$
in $H(M)^\bot$. Pseudohermitian structures exist as the
orientability assumption implies that $H(M)^\bot \approx M \times
{\mathbb R}$ (a diffeomorphism) i.e. $H(M)^\bot$ is a trivial line
bundle. For a review of the main notions of CR and pseudohermitian
geometry one may see \cite{kn:Dra}.
\par
Let $(M, T_{1,0}(M))$ be a CR manifold, of CR dimension $n$. Let
$\theta$ be a pseudohermitian structure on $M$. The {\em Levi
form} is
\[ L_\theta (Z , \overline{W}) = - i (d \theta )(Z ,
\overline{W}), \;\;\; Z,W \in T_{1,0}(M). \] $M$ is {\em
nondegnerate} if $L_\theta$ is nondegenerate for some $\theta$.
Two pseudohermitian structures $\theta$ and $\hat{\theta}$ are
related by
\begin{equation}
\hat{\theta} = f \; \theta \label{e:intro1}
\end{equation}
for some $C^\infty$ function $f : M \to {\mathbb R} \setminus \{ 0
\}$. Since $L_{\hat{\theta}} = f L_\theta$ nondegeneracy of $M$ is
a {\em CR invariant} notion, i.e. it is invariant under a
transformation (\ref{e:intro1}) of the pseudohermitian structure.
The whole setting bears an obvious analogy to conformal geometry
(a fact already exploited by many authors, cf. e.g.
\cite{kn:DrTo}, \cite{kn:Tan}-\cite{kn:Web}). If $M$ is
nondegenerate then any pseudohermitian structure $\theta$ on $M$
is actually a {\em contact form}, i.e. $\theta \wedge (d \theta
)^n$ is a volume form on $M$. By a fundamental result of N. Tanaka
and S. Webster (cf. {\em op. cit.}) on any nondegenerate CR
manifold on which a contact form $\theta$ has been fixed there is
a canonical linear connection $\nabla$ (the {\em Tanaka-Webster
connection} of $(M , \theta )$) compatible to the Levi
distribution and its complex structure, as well as to the Levi
form. Precisely, let $T$ be the globally defined nowhere zero
tangent vector field on $M$, transverse to $H(M)$, uniquely
determined by $\theta (T) = 1$ and $T \, \rfloor \, d \theta = 0$
(the {\em characteristic direction} of $d \theta$). Let
\[ G_\theta (X,Y) = (d \theta )(X , J Y), \;\;\; X,Y \in H(M), \]
(the {\em real} Levi form) and consider the semi-Riemannian metric
$g_\theta$ on $M$ given by
\[ g_\theta (X,Y) = G_\theta (X,Y), \;\;\; g_\theta (X,T) = 0,
\;\;\; g_\theta (T,T) = 1, \] for any $X,Y \in H(M)$ (the {\em
Webster metric} of $(M , \theta )$). Let us extend $J$ to an
endomorphism of the tangent bundle by setting $J T = 0$. Then
there is a unique linear connection $\nabla$ on $M$ such that i)
$H(M)$ is parallel with respect to $\nabla$, ii) $\nabla g_\theta
= 0$, $\nabla J = 0$, and iii) the torsion $T_\nabla$ of $\nabla$
is {\em pure}, i.e.
\begin{equation}
\label{e:intro2} T_\nabla (Z,W) = T_\nabla (\overline{Z},
\overline{W}) = 0, \;\;\; T_\nabla (Z, \overline{W}) = 2 i
L_\theta (Z, \overline{W}) T,
\end{equation}
for any $Z,W \in T_{1,0}(M)$, and
\begin{equation}\label{e:intro3}
\tau \circ J + J \circ \tau = 0,
\end{equation}
where $\tau (X) = T_\nabla (T , X)$ for any $X \in T(M)$ (the {\em
pseudohermitian torsion} of $\nabla$). The Tanaka-Webster
connection is a pseudohermitian analog to both the Levi-Civita
connection in Riemannian geometry and the Chern connection in
Hermitian geometry.
\par
A CR manifold $M$ is {\em strictly pseudoconvex} if $L_\theta$ is
positive definite for some $\theta$. If this is the case then the
Webster metric $g_\theta$ is a Riemannian metric on $M$ and if we
set $\varphi = J$, $\xi = - T$, $\eta = - \theta$ and $g =
g_\theta$ then $(\varphi , \xi , \eta , g)$ is a contact metric
structure on $M$. Also $(\varphi , \xi , \eta , g)$ is normal if
and only if $\tau = 0$. If this is the case $g_\theta$ is a
Sasakian metric and $(M , \theta )$ is a Sasakian manifold.

We proceed by recalling a few concepts from {\em sub-Riemannian
geometry} (cf. e.g. R.S. Strichartz, \cite{kn:Str}) on a strictly
pseudoconvex CR manifold. Let $(M , T_{1,0}(M))$ be a strictly
pseudoconvex CR manifold, of CR dimension $n$. Let $\theta$ be a
contact form on $M$ such that the Levi form $G_\theta$ is positive
definite. The Levi distribution $H(M)$ is {\em bracket generating}
i.e. the vector fields which are sections of $H(M)$ together with
all brackets span $T_x (M)$ at each point $x \in M$, merely as a
consequence of the nondegeneracy of the given CR structure.
Indeed, let $\nabla$ be the Tanaka-Webster connection of $(M ,
\theta )$ and let $\{ T_\alpha : 1 \leq \alpha \leq n \}$ be a
local frame of $T_{1,0}(M)$, defined on the open set $U \subseteq
M$. By the purity property (\ref{e:intro2})
\begin{equation}
\Gamma_{\alpha\overline{\beta}}^{\overline{\gamma}}
T_{\overline{\gamma}} - \Gamma_{\overline{\beta}\alpha}^\gamma
T_\gamma - [T_\alpha , T_{\overline{\beta}}] = 2 i
g_{\alpha\overline{\beta}} T, \label{e:29}
\end{equation}
where $\Gamma^A_{BC}$ are the coefficients of $\nabla$ with
respect to $\{ T_\alpha \}$
\[ \nabla_{T_B} T_C = \Gamma^A_{BC} T_A \]
and $g_{\alpha\overline{\beta}} = L_\theta (T_\alpha ,
T_{\overline{\beta}})$. Our conventions as to the range of indices
are $A,B,C, \cdots \in \{ 0, 1, \cdots , n , \overline{1}, \cdots
, \overline{n} \}$ and $\alpha , \beta , \gamma , \cdots  \in \{
1, \cdots , n \}$ (where $T_0 = T$). Note that $\{ T_\alpha ,
T_{\overline{\alpha}} , T \}$ is a local frame of $T(M) \otimes
{\mathbb C}$ on $U$. If $T_\alpha = X_\alpha - i J X_\alpha$ are
the real and imaginary parts of $T_\alpha$ then (\ref{e:29}) shows
that $\{ X_\alpha , J X_\alpha \}$ together with their brackets
span the whole of $T_x (M)$, for any $x \in U$. Actually more has
been proved. Given $x \in M$ and $v \in H(M)_x \setminus \{ 0 \}$
there is an open neighborhood $U \subseteq M$ of $x$ and a local
frame $\{ T_\alpha \}$ of $T_{1,0}(M)$ on $U$ such that $T_1 (x) =
v - i J_x v$, hence $v$ is a $2$-{\em step bracket generator} so
that $H(M)$ satisfies the {\em strong bracket generating
hypothesis} (cf. the terminology in \cite{kn:Str}, p. 224).
\par
Let $x \in M$ and $g(x) : T^*_x (M) \to H(M)_x$ determined by
\[ G_{\theta , x} (v , g(x) \xi ) = \xi (v), \;\;\; v \in H(M)_x ,
\;\; \xi \in T_x^* (M). \] Note that the kernel of $g$ is
precisely the conormal bundle $H(M)^\bot$. In other words
$G_\theta$ is a {\em sub-Riemannian metric} on $H(M)$ and $g$ its
alternative description (cf. also (2.1) in \cite{kn:Str}, p. 225).
If $\hat{\theta} = e^u \theta$ is another contact form such that
$G_{\hat{\theta}}$ is positive definite ($u \in C^\infty (M)$)
then $\hat{g} = e^{-u} g$. Clearly if the Levi form $L_\theta$ is
only nondegenerate then $(M , H(M), G_\theta )$ is a {\em
sub-Lorentzian manifold}, cf. the terminology in \cite{kn:Str}, p.
224.
\par
Let $\gamma : I \to M$ be a piecewise $C^1$ curve (where $I
\subseteq {\mathbb R}$ is an interval). Then $\gamma$ is a {\em
lengthy curve} if $\dot{\gamma}(t) \in H(M)_{\gamma (t)}$ for
every $t \in I$ such that $\dot{\gamma}(t)$ is defined. For
instance, any geodesic of $\nabla$ (i.e. any $C^1$ curve $\gamma
(t)$ such that $\nabla_{\displaystyle{\dot{\gamma}}} \dot{\gamma}
= 0$) of initial data $(x, v)$, $v \in H(M)_x$, is lengthy (as a
consequence of $\nabla g_\theta = 0$ and $\nabla T = 0$). A
piecewise $C^1$ curve $\xi : I \to T^* (M)$ is a {\em cotangent
lift} of $\gamma$ if $\xi (t) \in T_{\gamma (t)}^* (M)$ and
$g(\gamma (t)) \xi (t) = \dot{\gamma}(t)$ for every $t$ (where
defined). Clearly cotangent lifts of a given lengthy curve
$\gamma$ exist (cf. also Proposition \ref{p:xi0} below). Also,
cotangent lifts of $\gamma$ are uniquely determined modulo
sections of the conormal bundle $H(M)^\bot$ along $\gamma$. That
is, if $\eta : I \to T^* (M)$ is another cotangent lift of
$\gamma$ then $\eta (t) - \xi (t) \in H(M)^\bot_{\gamma (t)}$ for
every $t$. The {\em length} of a lengthy curve $\gamma : I \to M$
is given by
\[ L(\gamma ) = \int_I \{ \xi (t) \left[ g(\gamma (t))
\xi (t) \right]\}^{1/2} \; d t . \] The definition doesn't depend
upon the choice of cotangent lift $\xi$ of $\gamma$. The {\em
Carnot-Carath\'eodory distance} $\rho (x,y)$ among $x, y \in M$ is
the infimum of the lengths of all lengthy curves joining $x$ and
$y$. That $\rho$ is indeed a distance function on $M$ follows from
a theorem of W.L. Chow, \cite{kn:Cho}, according to which any two
points $x , y \in M$ may be joined by a lengthy curve (provided
that $M$ is connected).
\par
Let $g_\theta$ be the Webster metric of $(M , \theta )$. Then
$g_\theta$ is a {\em contraction} of the sub-Riemannian metric
$G_\theta$ ($G_\theta$ is an {\em expansion} of $g_\theta$), cf.
\cite{kn:Str}, p. 230. Let $d$ be the distance function
corresponding to the Webster metric. The length $L(\gamma )$ of a
lengthy curve $\gamma$ is precisely its length with respect to
$g_\theta$ hence
\begin{equation} d(x,y) \leq \rho (x,y), \;\;\; x,y \in M.
\label{e:30}
\end{equation}
While $\rho$ and $d$ are known to be inequivalent distance
functions, they do determine the same topology. For further
details on Carnot-Ca\-ra\-th\'e\-o\-do\-ry metrics see J.
Mitchell, \cite{kn:Mic}.
\par
Let $(U, x^1 , \cdots , x^{2n+1})$ be a system of local
coordinates on $M$ and let us set $G_{ij} = g_\theta (\partial_i ,
\partial_j )$ (where $\partial_i$ is short for $\partial /\partial
x^i$) and $[G^{ij}] = [G_{ij}]^{-1}$. Using
\[ G_\theta (X , g \; d x^i ) = (d x^i )(X), \;\;\; X \in H(M), \]
for $X = \partial_k - \theta_k T$ (where $\theta_i = \theta
(\partial_i )$) leads to
\begin{equation}\label{e:31} g^{ij}
(G_{jk} - \theta_j \theta_k ) = \delta^i_k - \theta_k T^i
\end{equation}
where $g \; d x^i = g^{ij} \partial_j$ and $T = T^i \partial_i$.
On the other hand $g^{ij} \theta_j = \theta (g \; d x^i ) = 0$ so
that (\ref{e:31}) yields
\begin{equation}
g^{ij} = G^{ij} - T^i T^j . \label{e:32}
\end{equation}
As an application we introduce a {\em canonical} cotangent lift of
a given lengthy curve on $M$.
\begin{proposition} Let $\gamma : I \to M$ be a lengthy curve and
let $\xi : I \to T^* (M)$ be given by $\xi (t) T_{\gamma (t)} = 1$
and $\xi (t) X = g_\theta (\dot{\gamma} , X)$, for any $X \in
H(M)_{\gamma (t)}$. Then $\xi$ is a cotangent lift of $\gamma$.
\label{p:xi0}
\end{proposition}
\noindent {\em Proof}. Let $x^i (t)$ be the components of $\gamma$
with respect to the chosen local coordinate system. By the very
definition of $\xi$
\begin{equation}
\label{e:xi0} \xi_j = G_{ij} \; \frac{d x^i}{d t} + \theta_j \, .
\end{equation}
Hence
\[ g \; \xi = \xi_j g^{ij} \partial_i = g^{ij} (G_{jk} \, \frac{d
x^k}{d t} + \theta_j ) \partial_i = g^{ij} G_{jk} \, \frac{d
x^k}{d t} \, \partial_i = \]
\[ = (G^{ij} - T^i T^j )G_{jk} \, \frac{d
x^k}{d t} \, \partial_i = (\delta^i_k - T^i \theta_k )\frac{d
x^k}{d t} \, \partial_i = \]
\[ = \dot{\gamma}(t) - \theta (\dot{\gamma}(t)) T =
\dot{\gamma}(t). \] We recall (cf. \cite{kn:Str}, p. 233) that a
{\em sub-Riemannian geodesic} is a $C^2$ curve $\gamma (t)$ in $M$
satisfying the Hamilton-Jacobi equations associated to the
Hamiltonian function $H(x, \xi ) = \frac{1}{2} \; g^{ij}(x) \xi_i
\xi_j$ that is
\begin{equation}
\frac{d x^i}{d t} = g^{ij}(\gamma (t)) \xi_j (t), \label{e:33}
\end{equation}
\begin{equation}
\frac{d \xi_k}{d t} = - \frac{1}{2} \; \frac{\partial
g^{ij}}{\partial x^k}(\gamma (t)) \xi_i (t) \xi_j (t),
\label{e:34}
\end{equation}
for some cotangent lift $\xi (t) \in T^* (M)$ of $\gamma (t)$. Our
purpose is to show that
\begin{theorem}
Let $M$ be a strictly pseudoconvex CR manifold and
$\theta$ a contact form on $M$ such that $G_\theta$ is positive
definite. A $C^2$ curve $\gamma (t) \in M$, $|t| < \epsilon$, is a
sub-Riemannian geodesic of $(M , H(M), G_\theta )$ if and only if
$\gamma (t)$ is a solution to
\begin{equation}
\nabla_{\displaystyle{\dot{\gamma}}} \dot{\gamma} = - 2 b(t) J
\dot{\gamma}, \;\;\; b^\prime (t) = A (\dot{\gamma},
\dot{\gamma}), \;\;\; |t| < \epsilon , \label{e:35}
\end{equation}
with $\dot{\gamma}(0) \in H(M)_{\gamma (0)}$, for some $C^2$
function $b : (-\epsilon , \epsilon ) \to {\mathbb R}$. Here $A
(X,Y) = g_\theta (\tau X , Y)$ is the pseudohermitian torsion of
$(M , \theta )$.
\label{t:2}
\end{theorem}
According to the terminology in \cite{kn:Str}, p. 237, the
canonical cotangent lift $\xi (t)$ of a given lengthy curve
$\gamma (t)$ is the one determined by the orthogonality
requirement
\begin{equation}
V_j (\xi ) \Gamma^j (\xi , v) = 0, \label{e:36}
\end{equation}
for any $v \in H(M)^\bot_{\gamma (t)}$ and any $|t| < \epsilon$,
where
\[ V_k (\xi ) = \frac{d \xi_k}{d t} + \frac{1}{2} \frac{\partial
g^{ij}}{\partial x^k} \xi_i \xi_j \, , \]
\[ \Gamma^i (\xi , v) = \Gamma^{ijk} \xi_j v_k \, , \;\;\;
\Gamma^{ijk} = \frac{1}{2} ( g^{\ell j} \frac{\partial
g^{ik}}{\partial x^\ell} + g^{\ell k} \frac{\partial
g^{ij}}{\partial x^\ell} - g^{\ell i} \frac{\partial
g^{jk}}{\partial x^\ell} ). \] Let $\gamma (t)$ be a lengthy curve
and $\xi_0 (t)$ the cotangent lift of $\gamma (t)$ furnished by
Proposition \ref{p:xi0}. Then any other cotangent lift $\xi (t)$
is given by
\begin{equation}
\xi (t) = \xi_0 (t) + a(t) \, \theta_{\gamma (t)} \, , \;\;\; |t|
< \epsilon , \label{e:37}
\end{equation}
for some $a : (- \epsilon , \epsilon ) \to {\mathbb R}$. We shall
need the following result (a replica of Lemma 4.4. in
\cite{kn:Str}, p. 237)
\begin{lemma}
The unique cotangent lift $\xi (t)$ of $\gamma (t)$ satisfying the
orthogonality condition {\rm (\ref{e:36})} is given by {\em
(\ref{e:37})} where
\[ a(t) = - \frac{1}{2} |\dot{\gamma}(t)|^{-2} g_\theta
(\nabla_{\displaystyle{\dot{\gamma}}} \dot{\gamma} \, , \, J
\dot{\gamma}) - 1, \;\;\; |t| < \epsilon . \] \label{l:2}
\end{lemma}
\noindent {\em Proof}. By (\ref{e:xi0}) and (\ref{e:37})
\[ V_k (\xi ) = V_k (\xi_0 ) + a^\prime (t) \theta_k + a(t)
\frac{\partial \theta_k}{\partial x^\ell} \frac{d x^\ell}{d t} +
\] \[ +  \frac{1}{2} \frac{\partial g^{ij}}{\partial x^k} [ a(t)
(\xi^0_i \theta_j + \xi^0_j \theta_i ) + a(t)^2 \theta_i \theta_j
]  \] (where $\xi_0 = \xi^0_i \, d x^i$) and using
\[ \frac{\partial g^{ij}}{\partial x^k} \, \theta_i \theta_j = 0
\]
we obtain
\begin{equation}
V_i (\xi ) = V_i (\xi_0 ) + a^\prime (t) \theta_i + 2 a (t) (d
\theta )(\dot{\gamma} , \partial_i ). \label{e:39}
\end{equation}
Note that $\Gamma^i (\xi , v) = \Gamma^i (\xi_0 , v)$ and
$\Gamma^{ijk} \theta_j v_k = 0$, for any $v \in H(M)^\bot_{\gamma
(t)}$. Let us contract (\ref{e:39}) with $\Gamma^i (\xi , v)$ and
use (\ref{e:36}) and $\Gamma^i (\xi_0 , v) \theta_i = 0$. This
ought to determine $a(t)$. Indeed
\begin{equation}
\label{e:40} V_i (\xi_0 ) \Gamma^i (\xi_0 , v) + 2a(t) (d \theta
)(\dot{\gamma} , \Gamma (\xi_0 , v)) = 0,
\end{equation}
where $\Gamma (\xi , v) = \Gamma^i (\xi , v) \partial_i$. On the
other hand, a calculation based on (\ref{e:32})-(\ref{e:xi0})
shows that
\[ V_k (\xi_0 ) = G_{k\ell} (\frac{d^2 x^\ell}{d t^2} + \left|
\begin{array}{c} \ell \\ ij \end{array} \right| \frac{d x^i}{d t}
\frac{dx^j}{d t} ) + 2 (d \theta )(\dot{\gamma} , \partial_k ), \]
\[ \left|
\begin{array}{c} \ell \\ ij \end{array} \right| = G^{\ell k}
|ij,k| , \;\;\; |ij,k| = \frac{1}{2} (\frac{\partial
G_{ik}}{\partial x^j} + \frac{\partial G_{jk}}{\partial x^i} -
\frac{\partial G_{ij}}{\partial x^k} ), \] hence
\begin{equation}
V_i (\xi_0 ) = G_{ij} (D_{\displaystyle{\dot{\gamma}}}
\dot{\gamma} )^j + 2 (d \theta )(\dot{\gamma} , \partial_i ) ,
\label{e:41}
\end{equation}
where $D$ is the Levi-Civita connection of $(M , g_\theta )$. Then
(\ref{e:40})-(\ref{e:41}) yield
\[ g_\theta (D_{\displaystyle{\dot{\gamma}}} \dot{\gamma} \, , \,
\Gamma (\xi_0 , v)) + 2 (a(t) +1) (d \theta )(\dot{\gamma} \, , \,
\Gamma (\xi_0 , v)) = 0, \]  for any $v \in H(M)^\bot_{\gamma
(t)}$. Yet $H(M)^\bot$ is the span of $\theta$ hence
\[ g_\theta (\Gamma (\xi_0 , \theta ) \, , \,
D_{\displaystyle{\dot{\gamma}}} \dot{\gamma} + 2(a(t)+1) J
\dot{\gamma}) = 0 \] and
\[ \Gamma^i (\xi_0 , \theta ) = -  G^{ij} (d \theta
)(\dot{\gamma} \, , \, \partial_j ), \] (because of $T \, \rfloor
\, d \theta = 0$) yields
\begin{equation}
\label{e:42} 2(a(t) + 1) |\dot{\gamma}(t)|^2 + g_\theta
(D_{\displaystyle{\dot{\gamma}}} \dot{\gamma} \, , \, J
\dot{\gamma}) = 0.
\end{equation}
Lemma \ref{l:2} is proved. At this point we may prove Theorem
\ref{t:2}. Let $\gamma (t) \in M$ be a sub-Riemannian geodesic of
$(M , H(M), G_\theta )$. Then there is a cotangent lift $\xi (t)
\in T^* (M)$ of $\gamma (t)$ (given by (\ref{e:37}) for some $a :
(-\epsilon , \epsilon ) \to {\mathbb R}$) such that $V (\xi ) = 0$
(where $V(\xi ) = V^i (\xi )
\partial_i$). In particular the orthogonality condition
(\ref{e:36}) is identically satisfied, hence $a(t)$ is determined
according to Lemma \ref{l:2}. Using (\ref{e:39}) and (\ref{e:41})
the sub-Riemannian geodesics equations are
\[ G_{ij} (D_{\displaystyle{\dot{\gamma}}} \dot{\gamma} )^j +
a^\prime (t) \theta_i + 2 (a(t)+1) (d \theta )(\dot{\gamma} \, ,
\, \partial_i ) = 0 \] or
\begin{equation}
D_{\displaystyle{\dot{\gamma}}} \dot{\gamma} + a^\prime (t) T +
2(a(t) + 1) J \dot{\gamma} = 0. \label{e:43}
\end{equation}
We recall (cf. e.g. \cite{kn:DrTo}) that $D = \nabla - (d \theta +
A) \otimes T$ on $H(M) \otimes H(M)$ hence (by the uniqueness of
the direct sum decomposition $T(M) = H(M) \oplus {\mathbb R} T$)
the equations (\ref{e:43}) become
\[ \nabla_{\displaystyle{\dot{\gamma}}} \dot{\gamma} + 2 (a(t)+1) J
\dot{\gamma} = 0, \;\;\; a^\prime (t) = A(\dot{\gamma} ,
\dot{\gamma}), \] (and we set $b = a + 1$). Theorem \ref{t:2} is
proved.
\begin{corollary} Let $M$ be a strictly pseudoconvex CR manifold
and $\theta$ a contact form on $M$ with vanishing pseudohermitian
torsion $(\tau = 0)$. Then any lengthy geodesic of the
Tanaka-Webster connection $\nabla$ of $(M , \theta )$ is a
sub-Riemannian geodesic of $(M , H(M), G_\theta )$. Viceversa, if
every lengthy geodesic $\gamma (t)$ of $\nabla$ is a
sub-Riemannian geodesic then $\tau = 0$.
\label{c:relation}
\end{corollary}
\noindent Indeed, if $\nabla_{\displaystyle{\dot{\gamma}}}
\dot{\gamma} = 0$ then the equations (\ref{e:35}) (with $b = 0$)
are identically satisfied.
\begin{proposition} Let $\gamma (t) \in M$ be a sub-Riemannian
geodesic and $s = \phi (t)$ a $C^2$ diffeomorphism. If $\gamma (t)
=  \overline{\gamma}(\phi (t))$ then $\overline{\gamma}(s)$ is a
sub-Riemannian geodesic if and only if $\phi$ is affine, i.e.
$\phi (t) = \alpha t + \beta$, for some $\alpha , \beta \in
{\mathbb R}$. In particular, every sub-Riemannian geodesic may be
reparametrized by arc length $\phi (t) = \int_0^t
|\dot{\gamma}(u)| \, d u$. \label{p:2}
\end{proposition}
\noindent {\em Proof}. Set $k = |\dot{\gamma}(0)|^2 > 0$. By
taking the inner product of the first equation in (\ref{e:35}) by
$\dot{\gamma}(t)$ it follows that $d |\dot{\gamma}(t)|^2 /d t =
0$, hence $|\dot{\gamma}(t)|^2 = k$, $|t| < \epsilon$. Throughout
the proof an overbar indicates the similar quantities associated
to $\overline{\gamma}(s)$. In particular $\overline{k} =
\phi^\prime (0)^{-2} k$. Locally
\begin{equation}
\frac{d^2 x^i}{d t^2} + \Gamma^i_{jk} \frac{d x^j}{dt} \frac{d
x^k}{dt} = - 2(a+1) J^i_j \, \frac{d x^j}{d t} . \label{e:44}
\end{equation}
On the other hand, using (\ref{e:42}) and
\[ \frac{d^2 x^i}{d t^2} + \Gamma^i_{jk} \frac{d x^j}{d t} \frac{d
x^k}{d t} = \phi^{\prime\prime}(t) \frac{d^2 \overline{x}^i}{d
s^2} + \phi^\prime (t)^2 (\frac{d^2 \overline{x}^i}{d s^2} +
\Gamma^i_{jk} \frac{d \overline{x}^j}{d s} \frac{d
\overline{x}^k}{d s} ) \] we obtain
\[ k(a+1) = \overline{k} (\overline{a} + 1) \phi^\prime (t)^3 . \]
Then (\ref{e:44}) may be written
\[ k \phi^{\prime\prime}(t) \frac{d \overline{\gamma}}{d s} + 2
(\overline{a} + 1) \phi^\prime (t)^2 [\overline{k} \phi^\prime
(t)^2 - k ] J \, \frac{d \overline{\gamma}}{d s} = 0 \] hence
$\phi^{\prime\prime}(t) = 0$. Proposition \ref{p:2} is proved.
\par
Let $S^1 \to C(M) \stackrel{\pi}{\longrightarrow} M$ be the
canonical circle bundle over $M$ (cf. e.g. \cite{kn:Dra}, p. 104).
Let $\Sigma$ be the tangent to the $S^1$-action. Next, let us
consider the $1$-form $\sigma$ on $C(M)$ given by
\[ \sigma = \frac{1}{n+2} \{ d r + \pi^* (i
\omega^\alpha_\alpha - \frac{i}{2} \, g^{\alpha\overline{\beta}} d
g_{\alpha\overline{\beta}} - \frac{R}{4(n+1)} \, \theta ) \} , \]
where $r$ is a local fibre coordinate on $C(M)$ (so that locally
$\Sigma = \partial /\partial r$) and $R =
g^{\alpha\overline{\beta}} R_{\alpha\overline{\beta}}$ is the
pseudohermitian scalar curvature of $(M , \theta )$. Then $\sigma$
is a connection $1$-form in $S^1 \to C(M) \to M$. Given a tangent
vector $v \in T_x (M)$ and a point $z \in \pi^{-1}(x)$ we denote
by $v^\uparrow$ its horizontal lift with respect to $\sigma$, i.e.
the unique tangent vector $v^\uparrow \in {\rm Ker}(\sigma_z )$
such that $(d_z \pi )v^\uparrow = v$. The {\em Fefferman metric}
of $(M , \theta )$ is the Lorentz metric on $C(M)$ given by
\[ F_\theta = \pi^* \tilde{G}_\theta + 2 (\pi^* \theta ) \odot
\sigma , \] where $\tilde{G}_\theta = G_\theta$ on $H(M) \otimes
H(M)$ and $\tilde{G}_\theta (X, T) = 0$, for any $X \in T(M)$.
Also $\odot$ is the symmetric tensor product. We close this
section by demonstrating the following geometric interpretation of
sub-Riemannian geodesics (of a strictly pseudoconvex CR manifold).
\begin{theorem}
Let $M$ be a strictly pseudoconvex CR manifold, $\theta$ a contact
form on $M$ such that $G_\theta$ is positive definite, and
$F_\theta$ the Fefferman metric of $(M , \theta )$. For any
geodesic $z : (-\epsilon , \epsilon ) \to C(M)$ of $F_\theta$ if
the projection $\gamma (t) = \pi (z(t))$ is lengthy then $\gamma :
(- \epsilon , \epsilon ) \to M$ is a sub-Riemannian geodesic of
$(M , H(M), G_\theta )$. Viceversa, let $\gamma (t) \in M$ be a
sub-Riemannian geodesic. Then any solution $z(t) \in C(M)$ to the
ODE
\begin{equation}
\dot{z}(t) = \dot{\gamma}(t)^\uparrow + ((n+2)/2) b(t)
\Sigma_{z(t)}, \label{e:45}
\end{equation}
where $b(t) = a(t) + 1$ is given by {\rm (\ref{e:42})}, is a
geodesic of $F_\theta$. \label{t:3}
\end{theorem}
\noindent Here $\dot{\gamma}(t)^\uparrow \in {\rm
Ker}(\sigma_{z(t)})$ and $(d_{z(t)} \pi ) \dot{\gamma}(t)^\uparrow
= \dot{\gamma}(t)$. To prove Theorem \ref{t:3} we shall need the
following
\begin{lemma} For any $X,Y \in H(M)$
\[ \nabla^{C(M)}_{X^\uparrow} Y^\uparrow = (\nabla_X Y)^\uparrow -
(d \theta )(X,Y) T^\uparrow - (A(X,Y) + (d \sigma )(X^\uparrow ,
Y^\uparrow )) \hat{\Sigma}, \]
\[ \nabla^{C(M)}_{X^\uparrow} T^\uparrow  = (\tau X + \phi
X)^\uparrow , \]
\[ \nabla^{C(M)}_{T^\uparrow} X^\uparrow = (\nabla_T X + \phi
X)^\uparrow + 2(d \sigma )(X^\uparrow , T^\uparrow ) \hat{\Sigma},
\]
\[ \nabla^{C(M)}_{X^\uparrow} \hat{\Sigma} = \nabla^{C(M)}_{\hat{\Sigma}} X^\uparrow = (J
X)^\uparrow , \]
\[ \nabla^{C(M)}_{T^\uparrow} T^\uparrow = V^\uparrow , \;\;
\nabla^{C(M)}_{\hat{\Sigma}} \hat{\Sigma} = 0, \]
\[ \nabla^{C(M)}_{\hat{\Sigma}} T^\uparrow = \nabla^{C(M)}_{T^\uparrow} \hat{\Sigma} = 0,
\]
where $\phi : H(M) \to H(M)$ is given by $G_\theta (\phi X , Y) =
(d \sigma )(X^\uparrow , Y^\uparrow )$, and $V \in H(M)$ is given
by $G_\theta (V , Y) = 2 (d \sigma )(T^\uparrow , Y^\uparrow )$.
Also $\hat{\Sigma} = ((n+2)/2) \Sigma$. \label{l:3}
\end{lemma}
\noindent This relates the Levi-Civita connection $\nabla^{C(M)}$
of $F_\theta$ to the Tanaka-Webster connection of $(M, \theta )$.
Cf. \cite{kn:Dra2} for a proof of Lemma \ref{l:3}.
\par
{\em Proof of Theorem} \ref{t:3}. Let $z(t) \in C(M)$ be a
geodesic of $\nabla^{C(M)}$ and $\gamma (t) = \pi (z(t))$. Assume
that $\dot{\gamma}(t) \in H(M)_{\gamma (t)}$. Note that
$\dot{z}(t) - \dot{\gamma}(t)^\uparrow \in {\rm Ker}(d_{z(t)} \pi
)$ hence $\dot{z}(t)$ is given by (\ref{e:45}), for some $b : (-
\epsilon , \epsilon ) \to {\mathbb R}$. Then (by Lemma \ref{l:3})
\[ 0 = \nabla^{C(M)}_{\displaystyle{\dot{z}}} \dot{z} =
\nabla^{C(M)}_{\displaystyle{\dot{\gamma}^\uparrow}}
\dot{\gamma}^\uparrow + b^\prime (t) \hat{\Sigma} + 2 b(t) (J
\dot{\gamma})^\uparrow = \]
\[ = (\nabla_{\displaystyle{\dot{\gamma}}} \dot{\gamma} )^\uparrow
+ [b^\prime (t) - A(\dot{\gamma} \, , \, \dot{\gamma})]
\hat{\Sigma} + 2 b(t) (J \dot{\gamma})^\uparrow \] hence (by
$T(C(M)) = {\rm Ker}(\sigma ) \oplus {\mathbb R} \Sigma$) $\gamma
(t)$ satisfies the equations (\ref{e:35}), i.e. $\gamma (t)$ is a
sub-Riemannian geodesic. The converse is obvious.

\section{Jacobi fields on CR manifolds}
Let $M$ be a strictly pseudoconvex CR manifold endowed with a
contact form $\theta$ such that $G_\theta$ is positive definite.
Let $\nabla$ be the Tanaka-Webster connection of $(M , \theta)$.
Let $\gamma (t) \in M$ be a geodesic of $\nabla$, parametrized by
arc length. A {\em Jacobi field} along $\gamma$ is vector field
$X$ on $M$ satisfying to the second order ODE \begin{equation}
\nabla^2_{\dot{\gamma}} X + \nabla_{\dot{\gamma}} T_\nabla (X ,
\dot{\gamma}) + R(X , \dot{\gamma}) \dot{\gamma} = 0. \label{e:J1}
\end{equation}
Let $J_\gamma$ be the real linear space of all Jacobi fields of
$(M , \nabla )$. Then $J_\gamma$ is $(4n+2)$-dimensional (cf.
Prop. 1.1 in \cite{kn:KoNo}, Vol. II, p. 63). We denote by
$\hat{\gamma}$ the vector field along $\gamma$ defined by
$\hat{\gamma}_{\gamma (t)} = t \dot{\gamma}(t)$ for every value of
the parameter $t$. Note that $\dot{\gamma}$, $\hat{\gamma} \in
J_\gamma$. We establish
\begin{theorem}  Every
Jacobi field $X$ along a lengthy geodesic $\gamma$ of $\nabla$ can
be uniquely decomposed in the following form
\begin{equation}
X = a \dot{\gamma} + b \hat{\gamma} + Y \label{e:J2}
\end{equation}
where $a , b \in {\mathbb R}$ and $Y$ is a Jacobi field along
$\gamma$ such that
\begin{equation}
g_\theta (Y, \dot{\gamma})_{\gamma (t)} = - \int^t_0 \theta
(X)_{\gamma (s)} A(\dot{\gamma} , \dot{\gamma})_{\gamma (s)} d s .
\label{e:J3}
\end{equation}
In particular, if {\rm i)} $X_{\gamma (t)} \in H(M)_{\gamma (t)}$
for every $t$, or {\rm ii)} $(M , \theta )$ is a Sasakian manifold
{\rm (}i.e. $\tau = 0${\rm )}, then $Y$ is perpendicular to
$\gamma$. \label{t:J1}
\end{theorem} \noindent We need the following
\begin{lemma}
For any Jacobi field $X \in J_\gamma$ \[ \frac{d}{d t} \{ g_\theta
(X , \dot{\gamma}) \} + \theta (X)_{\gamma (t)} A(\dot{\gamma} ,
\dot{\gamma})_{\gamma (t)} = {\rm const.} \] \label{l:J1}
\end{lemma}
\noindent {\em Proof}. Let us take the inner product of the Jacobi
equation (\ref{e:J1}) by $\dot{\gamma}$ and use the skew symmetry
of $g_\theta (R(X,Y)Z , W)$ in the arguments $(Z, W)$ (a
consequence of $\nabla g_\theta = 0$) so that to get
\[ \frac{d^2}{d t^2} \{ g_\theta (X, \dot{\gamma})\} + \frac{d}{d
t} \{ g_\theta (T_\nabla (X , \dot{\gamma}) , \dot{\gamma}) \} = 0
. \] On the other hand, let us set $X_H = X - \theta (X)T$ (so
that $X_H \in H(M)$). Then
\[ g_\theta (T_\nabla (X , \dot{\gamma}), \dot{\gamma}) = - 2 \Omega (X_H ,
\dot{\gamma}) g_\theta (T , \dot{\gamma}) + \theta (X) g_\theta
(\tau (\dot{\gamma}) , \dot{\gamma}) \] or (as $\gamma$ is
lengthy)
\[ g_\theta (T_\nabla (X , \dot{\gamma}), \dot{\gamma}) = \theta (X)
A(\dot{\gamma} , \dot{\gamma}). \] Lemma \ref{l:J1} is proved.
Throughout the section we adopt the notation $X^\prime =
\nabla_{\dot{\gamma}} X$ and $X^{\prime\prime} =
\nabla^2_{\dot{\gamma}} X$.
\par
{\em Proof of Theorem} \ref{t:J1}. We set by definition
\[ a = g_\theta (X , \dot{\gamma})_{\gamma (0)} \, , \;\;\;
b = g_\theta (X^\prime , \dot{\gamma})_{\gamma (0)} + \theta
(X)_{\gamma (0)} A(\dot{\gamma} , \dot{\gamma})_{\gamma (0)} \, ,
\] and $Y = X - a \dot{\gamma} - b \hat{\gamma}$. Clearly $Y \in
J_\gamma$. Then, by Lemma \ref{l:J1}
\[ \frac{d}{d t} \{ g_\theta (Y , \dot{\gamma}) \} + \theta (Y)
A(\dot{\gamma} , \dot{\gamma}) = \alpha , \] for some $\alpha \in
{\mathbb R}$. Next we integrate from $0$ to $t$
\[ g_\theta (Y , \dot{\gamma})_{\gamma (t)} - g_\theta (Y ,
\dot{\gamma})_{\gamma (0)} + \int_0^t \theta (Y)_{\gamma (s)}
A(\dot{\gamma} , \dot{\gamma})_{\gamma (s)} d s = \alpha t \] and
substitute $Y$ from (\ref{e:J2}) (and use $\dot{\gamma}$,
$\hat{\gamma} \in H(M)$). Differentiating the resulting relation
with respect to $t$ at $t = 0$ gives $\alpha = 0$. Hence
\[ g_\theta (Y , \dot{\gamma}) + \int_0^t \theta (X)_{\gamma (s)}
A(\dot{\gamma} , \dot{\gamma})_{\gamma (s)} d s = 0. \] The
existence statement in Theorem \ref{t:J1} is proved. We need the
following terminology. Given $X \in J_\gamma$ a Jacobi field $Y
\in J_\gamma$ satisfying (\ref{e:J3}) is said to be {\em slant at
$\gamma (t)$ relative to $X$}. Also $Y$ is {\em slant} if it slant
at any point of $\gamma$. To check the uniqueness statement let $X
= a^\prime \dot{\gamma} + b^\prime \hat{\gamma} + Z$ be another
decomposition of $X$, where $a^\prime$, $b^\prime \in {\mathbb R}$
and $Z \in J_\gamma$ is slant (relative to $X$). Then
\[ (a + b t) \dot{\gamma}(t) + Y_{\gamma (t)} = (a^\prime +
b^\prime t)\dot{\gamma}(t) + Z_{\gamma (t)} \] and taking the
inner product with $\dot{\gamma}(t)$ yields $a + b t = a^\prime +
b^\prime t$, i.e. $a = a^\prime$, $b = b^\prime$ and $Y_{\gamma
(t)} = Z_{\gamma (t)}$. Q.e.d.
\begin{corollary} Suppose a Jacobi field $X \in J_\gamma$ is slant
at $\gamma (r)$ and at $\gamma (s)$ relative to itself, for some
$r \neq s$. Then $X$ is slant. In particular, if {\rm i)}
$X_{\gamma (t)} \in H(M)_{\gamma (t)}$ for every $t$, or {\rm ii)}
$(M , \theta )$ is a Sasakian manifold, and $X$ is perpendicular
to $\gamma$ at two points, it is perpendicular to $\gamma$ at
every point of $\gamma$. \label{c:ppdue}
\end{corollary}
\noindent {\em Proof}. By Theorem \ref{t:J1} we may decompose $X =
a \dot{\gamma} + b \hat{\gamma} + Y$, where $Y \in J_\gamma$ is
slant (relative to $X$). Taking the inner product of $X_{\gamma
(r)} = (a + b r) \dot{\gamma}(r) + Y_{\gamma (r)}$ with
$\dot{\gamma}(r)$ gives $a + b r = 0$. Similarly $a + b s = 0$
hence (as $r \neq s$) $a = b = 0$, so that $X = Y$. Q.e.d.

\section{CR manifolds without conjugate points}
Two points $x$ and $y$ on a lengthy geodesic $\gamma (t)$ are {\em
horizontally conjugate} if there is a Jacobi field $X \in
J_\gamma$ such that $X_{\gamma (t)} \in H(M)_{\gamma (t)}$ for
every $t$ and $X_x = X_y = 0$. As $T_\nabla$ is pure, the Jacobi
equation (\ref{e:J1}) may also be written
\begin{equation}
\label{e:J4} X^{\prime\prime} - 2 \Omega (X^\prime , \dot{\gamma})
T + \theta (X^\prime ) \tau (\dot{\gamma}) + \theta (X)
(\nabla_{\dot{\gamma}} \tau ) \dot{\gamma} + R(X , \dot{\gamma})
\dot{\gamma} = 0.
\end{equation}
Given $X \in J_\gamma$ one has (by (\ref{e:J4}))
\[ \frac{d}{d t} \{ g_\theta (X^\prime , X) \} = g_\theta
(X^{\prime\prime} , X) + g_\theta (X^\prime , X^\prime ) = \]
\[ = |X^\prime |^2 + 2 \theta (X) \Omega (X^\prime , \dot{\gamma})
- \theta (X^\prime ) A(\dot{\gamma} , X) - \] \[ - \theta (X)
g_\theta (\nabla_{\dot{\gamma}} \tau \dot{\gamma} , X) - g_\theta
(R(X, \dot{\gamma}) \dot{\gamma} , X). \] On the other hand (again
by $\nabla g_\theta = 0$)
\[ \theta (X^\prime ) A(\dot{\gamma} , X) +
\theta (X) g_\theta (\nabla_{\dot{\gamma}} \tau \dot{\gamma} , X)
= \]
\[ = \theta (X^\prime ) A(\dot{\gamma} , X) + \theta (X)
\frac{d}{d t} \{ A(\dot{\gamma} , X)\} - \theta (X) A(\dot{\gamma}
, X^\prime ) = \]
\[ = \frac{d}{d t} \{ \theta (X) A(\dot{\gamma} , X) \} - \theta
(X) A(\dot{\gamma} , X^\prime ) \] hence
\begin{equation}
\frac{d}{d t} \{ g_\theta (X^\prime , X) + \theta (X)
A(\dot{\gamma} , X) \} =  \label{e:J5}
\end{equation}
\[ = |X^\prime |^2 - g_\theta (R(X,
\dot{\gamma}) \dot{\gamma} , X) + \theta (X) [A(\dot{\gamma} ,
X^\prime ) + 2 \Omega (X^\prime , \dot{\gamma})]. \] S. Webster
(cf. \cite{kn:Web}) has introduced a notion of pseudohermitian
sectional curvature by setting
\begin{equation} k_\theta (\sigma ) =  \frac{1}{4} \; G_\theta (X,X)^{-2}
\; g_{\theta , x} ( R_x (X , J_x X) J_x X , X ) , \label{e:res}
\end{equation}
for any holomorphic $2$-plane $\sigma$ (i.e. a $2$-plane $\sigma
\subset H(M)_x$ such that $J_x (\sigma ) = \sigma$), where $\{ X ,
J_x X \}$ is a basis of $\sigma$. The coefficient $1/4$ makes the
sphere $\iota : S^{2n+1} \subset {\mathbb C}^{n+1}$ (endowed with
the contact form $\theta_0 = \iota^* [\frac{i}{2}
(\overline{\partial} - \partial ) |z|^2 ]$) have constant
curvature $+1$. Clearly, this is a pseudohermitian analog to the
notion of holomorphic sectional curvature in Hermitian geometry.
On the other hand, for any $2$-plane $\sigma \subset T_x (M)$ one
may set
\[ k_\theta (\sigma ) = \frac{1}{4} \; g_{\theta , x} (R_x (X,Y)Y,
X) \] where $\{ X , Y \}$ is a $g_{\theta , x}$-orthonormal basis
of $\sigma$. Cf. \cite{kn:KoNo}, Vol. I, p. 200, the definition of
$k_\theta (\sigma )$ doesn't depend upon the choice of orthonormal
basis in $\sigma$ because the curvature $R(X,Y,Z,W) = g_\theta
(R(Z,W)Y, X)$ of the Tanaka-Webster connection is skew symmetric
in both pairs $(X,Y)$ and $(Z,W)$. We refer to $k_\theta$ as the
{\em pseudohermitian sectional curvature} of $(M , \theta )$. {\em
A posteriori} the restriction (\ref{e:res}) of $k_\theta$ to
holomorphic $2$-planes is referred to as the {\em holomorphic
pseudohermitian sectional curvature} of $(M , \theta )$. As an
application of (\ref{e:J5}) we may establish
\begin{theorem}
If $(M , \theta )$ has nonpositive pseudohermitian sectional
curvature then $(M , \theta )$ has no horizontally conjugate
points. \label{t:J2}
\end{theorem} \noindent We need
\begin{lemma} For every Jacobi field $X \in J_\gamma$
\[ \frac{d}{d t} \{ \theta (X) \} - 2 \Omega (X ,
\dot{\gamma})_{\gamma (t)} = c = {\rm const.} \] \label{l:J2}
\end{lemma} \noindent
To prove Lemma \ref{l:J2} one merely takes the inner product of
(\ref{e:J4}) by $T$. \par {\em Proof of Theorem} \ref{t:J2}. The
proof is by contradiction. If there is a lengthy geodesic $\gamma
(t) \in M$ (parametrized by arc length) and a Jacobi field $X \in
J_\gamma$ such that $X_{\gamma (a)} = X_{\gamma (b)} = 0$ for two
values $a$ and $b$ of the parameter then we may integrate in
(\ref{e:J5}) so that to obtain
\begin{equation}
\label{e:J6} \int_a^b \{ |X^\prime |^2 - g_\theta (R(X ,
\dot{\gamma}) \dot{\gamma} , X) + \theta (X) [A(\dot{\gamma} ,
X^\prime ) + 2 \Omega (X^\prime , \dot{\gamma})]\} dt = 0.
\end{equation}
On the other hand \[ \theta (X) \Omega (X^\prime , \dot{\gamma}) =
\theta (X) \frac{d}{d t} \{ \Omega (X , \dot{\gamma})\} = \] \[ =
\frac{d}{d t} \{ \theta (X) \Omega (X , \dot{\gamma}) \} - \Omega
(X , \dot{\gamma}) \theta (X^\prime ). \] Then (by Lemma
\ref{l:J2})
\[ 2 \int^b_a \theta (X) \Omega (X^\prime , \dot{\gamma}) d t = - 2 \int_a^b
\Omega (X , \dot{\gamma}) \frac{d}{d t} \{ \theta (X) \} d t = \]
\[ = c \int_a^b \frac{d}{d t} \{ \theta (X) \} d t - \int_a^b
\theta (X^\prime )^2 d t = - \int_a^b \theta (X^\prime )^2 d t \]
hence (\ref{e:J6}) becomes
\[ \int_a^b \{ |X^\prime |^2 - g_\theta (R(X , \dot{\gamma})
\dot{\gamma} , X) + \theta (X) A(\dot{\gamma}, X^\prime ) - \theta
(X^\prime )^2 \} d t = 0. \] Finally, if $X \in H(M)$ then
$X^\prime \in H(M)$ and then (under the assumptions of Theorem
\ref{t:J2}) $X^\prime = 0$, a contradiction.

\section{Jacobi fields on CR manifolds of constant
pseudohermitian sectional curvature} As well known (cf. Example
2.1 in \cite{kn:KoNo}, Vol. II, p. 71) one may determine a basis
of $J_\gamma$ for any elliptic space form (a Riemannian manifold
of positive constant sectional curvature). Similarly, we shall
prove
\begin{proposition} Let $M$ be a strictly pseudoconvex CR
manifold of CR dimension $n$, $\theta$ a contact form with
$G_\theta$ positive definite and constant pseudohermitian
sectional curvature. Let $\gamma (t) \in M$ be a lengthy geodesic
of the Tanaka-Webster connection $\nabla$ of $(M , \theta )$,
parametrized by arc length. For each $v \in T_{\gamma (0)} (M)$ we
let $E(v)$ be the space of all vector fields $X$ along $\gamma$
defined by $X_{\gamma (t)} = (a t + b)Y_{\gamma (t)}$, where $a,b
\in {\mathbb R}$, $\nabla_{\dot{\gamma}} Y = 0$, $Y_{\gamma (0)} =
v$. Assume that $(M, \theta )$ has parallel pseudohermitian
torsion, i.e. $\nabla \tau = 0$. Then $T \in J_\gamma$. Let $\{
v_1 , \cdots , v_{2n-2} \} \subset H(M)_{\gamma (0)}$ such that
$\{ \dot{\gamma}(0), J_{\gamma (0)} \dot{\gamma}(0) , v_1 , \cdots
, v_{2n-2} \}$ is a $G_{\theta , \gamma (0)}$-orthonormal basis of
$H(M)_{\gamma (0)}$. Then
\[ E(\dot{\gamma}(0)) \oplus E(v_1 ) \oplus \cdots \oplus
E(v_{2n-2}) \subseteq {\mathcal H}_\gamma := J_\gamma \cap
\Gamma^\infty (\gamma^{-1} H(M)) \] if and only if \[ A_{\gamma
(0)} (\dot{\gamma}(0), \dot{\gamma}(0)) = 0, \;\;\; A_{\gamma (0)}
(v_i , \dot{\gamma}(0)) = 0, \;\;\; 1 \leq i \leq 2n-2, \] where
$\gamma^{-1} H(M)$ is the pullback of $H(M)$ by $\gamma$. If
additionally $(M , \theta)$ has vanishing pseudohermitian torsion
{\rm (}i.e. $(M , \theta )$ is Tanaka-Webster flat{\rm )} then
$E(T_{\gamma (0)}) \subset J_\gamma$. \label{p:psh1}
\end{proposition}
\noindent The proof of Proposition \ref{p:psh1} requires the
explicit form of the curvature tensor of the Tanaka-Webster
connection of $(M , \theta )$ when $k_\theta =$ const. This is
provided by
\begin{theorem} Let $M$ be a strictly pseudoconvex CR manifold
and $\theta$ a contact form on $M$ such that $G_\theta$ is
positive definite and $k_\theta (\sigma ) = c$, for some $c \in
{\mathbb R}$ and any $2$-plane $\sigma \subset T_x (M)$, $x \in
M$. Then $c = 0$ and the curvature of the Tanaka-Webster
connection of $(M , \theta )$ is given by
\begin{equation} R(X,Y)Z  = \Omega (Z,Y) \tau (X) - \Omega (Z,X) \tau (Y)
+ \label{e:J7}
\end{equation}
\[ + A(Z,Y) J X - A(Z,X) J Y, \]
for any $X,Y,Z \in T(M)$. In particular, if $(M , \theta)$ has
constant pseudohermitian sectional curvature and CR dimension $n
\geq 2$ then the Tanaka-Webster connection of $(M , \theta )$ is
flat if and only if $(M , \theta )$ has vanishing pseudohermitian
torsion {\rm (}$\tau = 0${\rm )}. \label{t:J3}
\end{theorem}
\noindent The proof of Theorem \ref{t:J3} is given in Appendix A.
By Theorem \ref{t:J3} there are no ``pseudohermitian space forms''
except for those of zero pseudohermitian sectional curvature and
these aren't in general flat. Cf. \cite{kn:DrTo} the term {\em
pseudohermitian space form} is reserved for manifolds of constant
{\em holomorphic} pseudohermitian sectional curvature (and then
examples with arbitrary $c \in {\mathbb R}$ abound, cf.
\cite{kn:DrTo}, Chapter 1).
\par
{\em Proof of Proposition} \ref{p:psh1}. By (\ref{e:J7})
\[ R(X, \dot{\gamma})\dot{\gamma} = \Omega (X, \dot{\gamma}) \tau
(\dot{\gamma}) + A(\dot{\gamma} , \dot{\gamma}) J X - A(X ,
\dot{\gamma}) J \dot{\gamma} \] hence the Jacobi equation
(\ref{e:J4}) becomes
\begin{equation}
X^{\prime\prime} - 2 \Omega (X^\prime , \dot{\gamma}) T + \theta
(X^\prime ) \tau (\dot{\gamma}) + \label{e:psh1}
\end{equation}
\[ + \theta (X) (\nabla_{\dot{\gamma}} \tau )\dot{\gamma} + \Omega
(X, \dot{\gamma}) \tau (\dot{\gamma}) + A(\dot{\gamma} ,
\dot{\gamma}) J X - A(X , \dot{\gamma}) J \dot{\gamma} = 0. \] We
look for solutions to (\ref{e:psh1}) of the form $X_{\gamma (t)} =
f(t) T_{\gamma (t)}$. The relevant equation is
\[ f^{\prime\prime}(t) T + f^\prime (t) \tau (\dot{\gamma}) + f(t)
(\nabla_{\dot{\gamma}} \tau )\dot{\gamma} = 0 \] (by $\nabla T =
0$) or $f^{\prime\prime}(t) = 0$ and $f^\prime (t) \tau
(\dot{\gamma}) + f(t) (\nabla_{\dot{\gamma}} \tau )\dot{\gamma} =
0$. Therefore, if $\nabla \tau = 0$ then $T \in J_\gamma$ while if
$\tau = 0$ then $T, \hat{T} \in J_\gamma$, where $\hat{T}_{\gamma
(t)} = t T_{\gamma (t)}$. Next, we look for solutions to
(\ref{e:psh1}) of the form $X_{\gamma (t)} = f(t) Y_{\gamma (t)}$
where $Y$ is a vector field along $\gamma$ such that
$\nabla_{\dot{\gamma}} Y = 0$, $Y_{\gamma (0)} =: v \in
H(M)_{\gamma (0)}$, $|v| = 1$, and $g_{\theta , \gamma (0)}(v ,
J_{\gamma (0)} \dot{\gamma}(0)) = 0$. Substitution into
(\ref{e:psh1}) gives
\[ f^{\prime\prime}(t) Y + f(t) [A(\dot{\gamma}, \dot{\gamma}) J Y
- A(Y, \dot{\gamma}) J \dot{\gamma}] = 0 \] or (by taking the
inner product with $Y$) $f^{\prime\prime}(t) = 0$, i.e. $f(t) = a
t + b$, $a,b \in {\mathbb R}$. Therefore (with the notations in
Proposition \ref{p:psh1}) $E(v_i ) \subset J_\gamma \cap
\Gamma^\infty (\gamma^{-1} H(M))$ if and only if $A_{\gamma
(0)}(\dot{\gamma}(0), \dot{\gamma}(0)) = 0$ and $A_{\gamma
(0)}(v_i , \dot{\gamma}(0)) = 0$. Also, to start with,
$E(\dot{\gamma}(0))$ (the space spanned by $\dot{\gamma}$ and
$\hat{\gamma}$) consists of Jacobi fields lying in $H(M)$. As $\{
\dot{\gamma}(t) , J_{\gamma (t)} \dot{\gamma}(t) , Y_{1, \gamma
(t)} , \cdots , Y_{2n-2 , \gamma (t)} \}$ is an orthonormal basis
of $H(M)_{\gamma (t)}$ (where $Y_i$ is the unique solution to
$(\nabla_{\dot{\gamma}} Y)_{\gamma (t)} = 0$, $Y_{\gamma (0)} =
v_i$) it follows that the sum $E(\dot{\gamma}(0)) + E(v_1 ) +
\cdots + E(v_{2n-2})$ is direct. Q.e.d. \vskip 0.1in Let $(M,
(\varphi , \xi , \eta , g))$ be a contact Riemannian manifold. Let
$X \in T_x (M)$ be a unit tangent vector orthogonal to $\xi$ and
$\sigma \subset T_x (M)$ the $2$-plane spanned by $\{ X , \varphi
X \}$ (a $\varphi$-{\em holomorphic plane}). We recall (cf. e.g.
\cite{kn:Bla}, p. 94) that the $\varphi$-{\em sectional curvature}
is the restriction of the sectional curvature $k$ of $(M , g)$ to
the $\varphi$-holomorphic planes. Let us set $H(X) = k(\sigma )$.
A Sasakian manifold of constant $\varphi$-sectional curvature
$H(X) = c$, $c \in {\mathbb R}$, is a {\em Sasakian space form}.
Compact Sasakian space forms have been classified in
\cite{kn:Kam}. By a result in \cite{kn:Bla}, p. 97, the Riemannian
curvature $R^D$ of a Sasakian space form $M$ (of
$\varphi$-sectional curvature $c$) is given by
\begin{equation}
R^D (X,Y) Z = \frac{c+3}{4} \{ g(Y,Z) X - g(X,Z) Y \} +
\label{e:curv0}
\end{equation}
\[ + \frac{c-1}{4} \{ \eta (Z)[\eta (X) Y - \eta (Y) X ] + [g(X,Z)
\eta (Y) - g(Y,Z) \eta (X)] \xi + \]
\[ + \Omega (Z,Y) \varphi X - \Omega (Z,X) \varphi Y + 2 \Omega
(X,Y) \varphi Z \} \] for any $X,Y,Z \in T(M)$.
 Given a strictly pseudoconvex CR manifold $M$ and a
contact form $\theta$ we recall (cf. e.g. (1.59) in
\cite{kn:DrTo}) that
\begin{equation} D = \nabla + (\Omega - A) \otimes T + \tau
\otimes \theta + 2 \theta \odot J. \label{e:DeNabla}
\end{equation}
A calculation based on (\ref{e:DeNabla}) leads to
\[ R^D (X,Y) Z = R(X,Y) Z + (L X \wedge L Y ) Z - 2 \Omega (X,Y) J
Z - \]
\[ - g_\theta (S(X,Y) , Z) T + \theta (Z) S(X,Y) - \]
\[ - 2 g_\theta ((\theta \wedge {\mathcal O})(X,Y), Z) T + 2
\theta (Z) (\theta \wedge {\mathcal O})(X,Y) \] for any $X,Y,Z \in
T(M)$, relating the Riemannian curvature $R^D$ of $(M , g_\theta
)$ to the curvature $R$ of the Tanaka-Webster connection. Here
\[ L = \tau + J, \;\;\; {\mathcal O} = \tau^2 + 2 J \tau - I, \]
and $(X \wedge Y)Z = g_\theta (X,Z) Y - g_\theta (Y,Z) X$. Also
$S(X,Y) = (\nabla_X \tau )Y - (\nabla_Y \tau ) X$. Let us assume
that $(M , \theta )$ is a Sasakian manifold ($\tau = 0$) whose
Tanaka-Webster connection is flat ($R = 0$). Then $S = 0$, $L = J$
and ${\mathcal O} = - I$ hence
\[ R^D (X,Y) Z = (J X \wedge J Y) Z - 2 \Omega (X,Y) J Z + \]
\[ + 2 g_\theta ((\theta \wedge I)(X,Y) , Z) T - 2 \theta (Z)
(\theta \wedge I)(X,Y) \] and a comparison to (\ref{e:curv0})
shows that
\begin{proposition} Let $(M , \theta )$ be a Sasakian manifold.
Then its Tanaka-Webster connection is flat if and only if $(M ,
(J, - T, - \theta , g_\theta ))$ is a Sa\-sa\-ki\-an space form of
$\varphi$-sectional curvature $c = - 3$. \label{p:menotre}
\end{proposition}
By Lemma \ref{l:conj6} below the dimension of ${\mathcal
H}_\gamma$ is at most $4n$. On a Sasakian space form we may
determine $4n-1$ independent vectors in ${\mathcal H}_\gamma$.
Indeed, by combining Propositions \ref{p:psh1} and \ref{p:menotre}
we obtain
\begin{corollary} $\,$ \par\noindent
Let $(M , \theta )$ be a Sasakian space form of
$\varphi$-sectional curvature $c = - 3$ and $\gamma (t) \in M$ a
lengthy geodesic of the Tanaka-Webster connection $\nabla$,
parametrized by arc length. Let $\{ v_1 , \cdots , v_{2n-2} \}
\subset H(M)_{\gamma (0)}$ such that $\{ \dot{\gamma}(0),
J_{\gamma (0)} \dot{\gamma}(0), v_1 , \cdots , v_{2n-2} \}$ is a
$G_{\theta , \gamma (0)}$-orthonormal basis of $H(M)_{\gamma
(0)}$. Let $X_i$ be the vector field along $\gamma$ determined by
\[ \nabla_{\dot{\gamma} (t)} X_i = 0, \;\;\; X_i (\gamma (0)) = v_i
, \] for $1 \leq i \leq 2n-2$. Then ${\mathcal S} = \{
\dot{\gamma}, \hat{\gamma}, J \dot{\gamma}, X_i , \hat{X}_i : 1
\leq i \leq 2n-2 \}$ is a free system in ${\mathcal H}_\gamma$
while ${\mathcal S} \cup \{ T , \hat{T} \}$ is free in $J_\gamma$.
Here if $Y$ is a vector field along $\gamma (t)$ we set
$\hat{Y}_{\gamma (t)} = t Y_{\gamma (t)}$ for every $t$.
\end{corollary}

\section{Conjugate points on Sasakian manifolds}
Let $(M , \theta )$ be a Sasakian manifold and $\gamma : [a,b] \to
M$ a geodesic of the Tanaka-Webster connection $\nabla$,
parametrized by arc length. Given a piecewise differentiable
vector field $X$ along $\gamma$ we set
\[ I_a^b (X) = \int_a^b \{ g_\theta (\nabla_{\dot{\gamma}} X , \nabla_{\dot{\gamma}}
X) - g_\theta (R(X, \dot{\gamma})\dot{\gamma} , X) \}_{\gamma (t)}
d t \] where $R$ is the curvature of $\nabla$. We shall prove the
following
\begin{proposition} Let $(M , \theta)$ be a Sasakian manifold and
$\gamma (t) \in M$, $a \leq t \leq b$, a lengthy geodesic of
$\nabla$, parametrized by arc length, such that $\gamma (a)$ has
no conjugate point along $\gamma$. Let $Y \in {\mathcal H}_\gamma$
be a horizontal Jacobi field along $\gamma$ such that $Y_{\gamma
(a)} = 0$ and $Y$ is perpendicular to $\gamma$. Let $X$ be a
piecewise differentiable vector field along $\gamma$ such that
$X_{\gamma (a)} = 0$ and $X$ is perpendicular to $\gamma$. If
$X_{\gamma (b)} = Y_{\gamma (b)}$ then
\begin{equation}
I_a^b (X) \geq I_a^b (Y) \label{e:conj1}
\end{equation}
and the equality holds if and only if $X = Y$. \label{p:conj1}
\end{proposition}
{\em Proof}. Let $J_{\gamma , a}$ be the space of all Jacobi
fields $Z \in J_\gamma$ such that $Z_{\gamma (a)} = 0$. By Prop.
1.1 in \cite{kn:KoNo}, Vol. II, p. 63, $J_{\gamma , a}$ has
dimension $2n+1$. Moreover, let $J_{\gamma , a, \bot}$ be the
space of all $Z \in J_{\gamma , a}$ such that $g_\theta (Z ,
\dot{\gamma})_{\gamma (t)} = 0$, for every $t$. Then by Theorem
\ref{t:J1} it follows that $J_{\gamma , a , \bot}$ has dimension
$2n$. We shall need the following
\begin{lemma} For every Sasakian manifold $(M , \theta )$ the
characteristic direction $T$ of $(M , \theta )$ is a Jacobi field
along any geodesic $\gamma : [a,b] \to M$ of $\nabla$. Also, if
$T_a$ is the vector  field along $\gamma$ given by $T_{a , \gamma
(t)} = (t-a) T_{\gamma (t)}$, $a \leq t \leq b$, and $\gamma$ is
lengthy then $T_a \in J_{\gamma , a , \bot}$ and $T_{a , \gamma
(t)} \neq 0$, $a < t \leq b$. \label{l:conj1}
\end{lemma}
{\em Proof}. Let ${\mathcal J}_\gamma$ be the Jacobi operator.
Then
\[ {\mathcal J}_\gamma T = T^{\prime\prime} - 2 \Omega (T^\prime ,
\dot{\gamma}) T + R(T , \dot{\gamma}) \dot{\gamma} = R(T,
\dot{\gamma})\dot{\gamma} \] as $T^\prime = \nabla_{\dot{\gamma}}
T = 0$. On the other hand, on any nondegenerate CR manifold with
$S = 0$ (i.e. $S(X,Y) \equiv (\nabla_X \tau )Y - (\nabla_Y \tau )X
= 0$, for any $X,Y \in T(M)$) the curvature of the Tanaka-Webster
connection satisfies
\begin{equation}
R(T, X) X = 0, \;\;\; X \in T(M), \label{e:conjR}
\end{equation}
hence ${\mathcal J}_\gamma T = 0$. As $R(T,T) = 0$ and $R(X,Y) T =
0$ it suffices to check (\ref{e:conjR}) for $X \in H(M)$, i.e.
locally $X = Z^\alpha T_\alpha + Z^{\overline{\alpha}}
T_{\overline{\alpha}}$. Then
\[ R(T,X)X = \{ {{R_\beta}^\gamma}_{0\alpha} Z^\alpha Z^\beta +
{{R_\beta}^\gamma}_{0\overline{\alpha}} Z^{\overline{\alpha}}
Z^\beta \} T_\gamma + \]
\[ + {{R_{\overline{\beta}}}^{\overline{\gamma}}}_{0\alpha} Z^\alpha Z^{\overline{\beta}} +
{{R_{\overline{\beta}}}^{\overline{\gamma}}}_{0\overline{\alpha}}
Z^{\overline{\alpha}} Z^{\overline{\beta}} \}
T_{\overline{\gamma}} \] and (by (1.85)-(1.86) in \cite{kn:DrTo},
section 1.4)
\[ {{R_\beta}^\gamma}_{0\alpha} = g^{\gamma\overline{\lambda}}
g_{\alpha\overline{\mu}}
S^{\overline{\mu}}_{\beta\overline{\lambda}} \, , \;\;\;
{{R_\beta}^\gamma}_{0\overline{\alpha}} =
g_{\lambda\overline{\alpha}} g^{\gamma\overline{\mu}}
S^\lambda_{\beta\overline{\mu}} \, . \] To complete the proof of
Lemma \ref{l:conj1} let $u(t) = t - a$. Then (by $T \, \rfloor \,
\Omega = 0$ and (\ref{e:conjR}))
\[ {\mathcal J}_\gamma T_a  = u^{\prime\prime} T - 2 u^\prime
\Omega (T , \dot{\gamma}) T + u R(T , \dot{\gamma}) \dot{\gamma} =
0. \]
\begin{lemma} Let $(M , \theta )$ be a Sasakian manifold and
$\gamma (t) \in M$ a geodesic of $\nabla$. If $X \in J_\gamma$
then $X_H \equiv X - \theta (X) T$ satisfies the second order ODE
\begin{equation}
\nabla^2_{\dot{\gamma}} X_H + R(X_H , \dot{\gamma}) \dot{\gamma} =
0. \label{e:conj3}
\end{equation}
\label{l:conj2}
\end{lemma}
{\em Proof}.
\[ 0 = {\mathcal J}_\gamma X = \nabla^2_{\dot{\gamma}} X_H +
\theta (X^{\prime\prime}) T - 2 \Omega (\nabla_{\dot{\gamma}} X_H
, \dot{\gamma}) T + R(X_H , \dot{\gamma}) \dot{\gamma} \] hence
(by the uniqueness of the direct sum decomposition $T(M) = H(M)
\oplus {\mathbb R} T$) $X_H$ satisfies (\ref{e:conj3}). \par Let
us go back to the proof of Proposition \ref{p:conj1}. Let us
complete $T_a$ to a linear basis $\{ T_a , Y_2 , \cdots , Y_{2n}
\}$ of $J_{\gamma , a , \bot}$ and set $Y_1 = T_a$ for simplicity.
Then $Y = a^i Y_i$ for some $a^i \in {\mathbb R}$, $1 \leq i \leq
2n$. Let us observe that for each $a < t \leq b$ the tangent
vectors
\[ \{ T_{a , \gamma (t)} , Y^H_{2 , \gamma (t)} , \cdots , Y^H_{2n
, \gamma (t)} \} \subset [{\mathbb R} \dot{\gamma}(t)]^\bot
\subset T_{\gamma (t)} (M) \] are linearly independent, where
$Y_j^H := Y_j - \theta (Y_j ) T$, $2 \leq j \leq 2n$. Indeed
\[ 0 = \alpha T_{a , \gamma (t)} + \sum_{j=2}^{2n} \alpha^j
Y^H_{j, \gamma (t)} = \{ \alpha - \sum_{j=2}^{2n}
\frac{\alpha^j}{t-a} \theta (Y_j )_{\gamma (t)} \} T_{a, \gamma
(t)} + \sum_{j=2}^{2n} \alpha^j Y_{j, \gamma (t)} \] implies
$\alpha^j = 0$, and then $\alpha = 0$, because $\{ Y_{i , \gamma
(t)} : 1 \leq i \leq 2n \}$ are linearly independent, for any $a <
t \leq b$. At their turn, the vectors $Y_{i , \gamma (t)}$ are
independent because $\gamma (a)$ has no conjugate point along
$\gamma$. The proof is by contradiction. Assume that
\begin{equation}
\lambda^i Y_{i, \gamma (t_0 )} = 0, \label{e:conj5}
\end{equation}
for some $a < t_0 \leq b$ and some $\lambda = (\lambda^1 , \cdots
, \lambda^{2n}) \in {\mathbb R}^{2n} \setminus \{ 0 \}$. Let us
set $Z_0 = \lambda^i Y_i \in J_\gamma$. Then
\[ \lambda \neq 0 \Longrightarrow Z_0 \neq 0, \]
\[ Z_0 \in J_{\gamma , a} \Longrightarrow Z_{0, \gamma (a)} = 0,
\;\;\; {\rm (\ref{e:conj5})} \Longrightarrow Z_{0, \gamma (t_0 )}
= 0,
\]
hence $\gamma (a)$ and $\gamma (t_0 )$ are conjugate along
$\gamma$, a contradiction. Yet $[{\mathbb R}\dot{\gamma}(t)]^\bot$
has dimension $2n$ hence
\[ X_{\gamma (t)} = f(t) T_{a, \gamma (t)} + \sum_{j=2}^{2n} f^j
(t) Y^H_{j , \gamma (t)} \, , \] for some piecewise differentiable
functions $f(t)$, $f^j (t)$. We set $f^1 = f$, $Z_1 = T_a$ and
$Z_j = Y^H_j$, $2 \leq j \leq 2n$, for simplicity. Then
\begin{equation}
|X^\prime |^2 = | \frac{d f^i}{d t} Z_i |^2 + |f^i Z^\prime_i |^2
+ 2 g_\theta (\frac{d f^i}{d t} Z_i , f^j Z^\prime_j ).
\label{e:conj6}
\end{equation}
Also (by (\ref{e:conjR}) and Lemma \ref{l:conj2})
\[ - g_\theta (R(X , \dot{\gamma})\dot{\gamma} , X) = - f^i
g_\theta (R(Z_i , \dot{\gamma})\dot{\gamma} , X) = \]
\[ = - \sum_{j=2}^{2n} f^j g_\theta (R(Z_j ,
\dot{\gamma})\dot{\gamma} , X) = \sum_{j=2}^{2n} f^j g_\theta
(Z^{\prime\prime}_j , X)  \] or (as $T_a^{\prime\prime} = 0$)
\begin{equation}
- g_\theta (R(X, \dot{\gamma})\dot{\gamma} , X) = g_\theta (f^i
Z_i^{\prime\prime} , f^j Z_j ). \label{e:conj7}
\end{equation}
Finally, note that
\begin{equation}
g_\theta (\frac{d f^i}{d t} Z_i , f^j Z^\prime_j ) + |f^i
Z^\prime_i |^2 + g_\theta (f^i Z_i , f^j Z_j^{\prime\prime}) =
\label{e:conj8}
\end{equation}
\[ = \frac{d}{d t} g_\theta (f^i Z_i , f^j Z^\prime_j ) - g_\theta
(f^i Z_i , \frac{d f^j}{d t} Z^\prime_j ). \] Summing up (by
(\ref{e:conj6})-(\ref{e:conj8}))
\begin{equation}
|X^\prime |^2 - g_\theta (R(X, \dot{\gamma})\dot{\gamma} , X) =
\label{e:conj9}
\end{equation}
\[ = \frac{d}{d t} g_\theta (f^i Z_i , f^j Z^\prime_j ) - g_\theta
(f^i Z_i , \frac{d f^j}{d t} Z^\prime_j ) + g_\theta (\frac{d
f^i}{d t} Z_i , f^j Z^\prime_j ) + | \frac{d f^i}{d t} Z_i |^2 \,
. \]
\begin{lemma} Let $(M , \theta )$ be a Sasakian manifold and
$\gamma (t) \in M$ a geodesic of $\nabla$. If $X$ and $Y$ are
solutions to $\nabla^2_{\dot{\gamma}} Z + R(Z,
\dot{\gamma})\dot{\gamma} = 0$ then
\begin{equation}
\label{e:adj4} \frac{d}{d t} \{ g_\theta (X,Y^\prime ) - g_\theta
(X^\prime , Y) \} = 0.
\end{equation}
In particular, if $X_{\gamma (a)} = 0$ and $Y_{\gamma (a)} = 0$ at
some point $\gamma (a)$ of $\gamma$ then
\[ g_\theta (X, Y^\prime ) - g_\theta (X^\prime , Y) = 0. \]
\label{l:conj4}
\end{lemma}
{\em Proof}. As $\tau = 0$ the $4$-tensor $R(X,Y,Z,W)$ possesses
the symmetry property $R(X,Y,Z,W) = R(Z,W,X,Y)$ (cf. (\ref{e:A8})
in Appendix A) one may subtract the identities
\[ \frac{d}{d t} g_\theta (X, Y^\prime ) = g_\theta (X^\prime , Y^\prime )
- g_\theta (X , R(Y, \dot{\gamma}) \dot{\gamma}), \]
\[ \frac{d}{d t} g_\theta (X^\prime , Y ) = g_\theta (X^\prime ,
Y^\prime ) - g_\theta (R(X, \dot{\gamma})\dot{\gamma}, Y) \] so
that to obtain (\ref{e:adj4}). Q.e.d. \par By Lemma \ref{l:conj2}
the fields $Z_j$, $2 \leq j \leq 2n$ satisfy
$\nabla_{\dot{\gamma}}^2 Z_j + R(Z_j , \dot{\gamma})\dot{\gamma} =
0$. Then we may apply Lemma \ref{l:conj4} to conclude that
\[ g_\theta (\frac{df^i}{d t} Z_i , f^j Z^\prime_j ) - g_\theta
(f^i Z_i , \frac{d f^j}{d t} Z^\prime_j ) = \] \[ = f^i \frac{d
f^j}{d t} \{ g_\theta (Z_j , Z^\prime_i ) - g_\theta (Z^\prime_j ,
Z_i ) \} = 0 \] so that (\ref{e:conj9}) becomes
\[ |X^\prime |^2 - g_\theta (R(X, \dot{\gamma})\dot{\gamma} , X) =
\frac{d}{d t} g_\theta (f^i Z_i , f^j Z^\prime_j ) + |\frac{d
f^i}{d t} Z_i |^2 \] and integration gives
\begin{equation}
I_a^b (X) = g_\theta (f^i Z_i , f^j Z^\prime_j )_{\gamma (b)} +
\int_a^b |\frac{d f^i}{d t} Z_i |^2 d t . \label{e:conj10}
\end{equation}
We wish to apply (\ref{e:conj10}) to the vector field $X = Y$. If
this is the case the functions $f^j$ are $f^1 (t) = a^1 +
(1/(t-a)) \sum_{j=2}^{2n} a^j \theta (Y_j )_{\gamma (t)} = 0$
(because of $Y_{\gamma (t)} \in H(M)_{\gamma (t)}$) and $f^j =
a^j$ (so that $d f^j /d t = 0$) for $2 \leq j \leq 2n$. Then (by
(\ref{e:conj10}))
\begin{equation}
I_a^b (Y) = g_\theta (\sum_{i=2}^{2n} a^i Z_i , \sum_{j=2}^{2n}
a^j Z_j^\prime )_{\gamma (b)} \, . \label{e:conj11}
\end{equation}
As $X_{\gamma (b)} = Y_{\gamma (b)}$ it follows that $f^1 (b) = 0$
and $f^j (b) = a^j$, $2 \leq j \leq 2n$, so that by subtracting
(\ref{e:conj10}) and (\ref{e:conj11}) we get
\[ I_a^b (X) - I_a^b (Y) = \int_a^b |\frac{d f^i}{d t} Z_i |^2 d t
\geq 0 \] and (\ref{e:conj1}) is proved. The equality $I_a^b (X) =
I_a^b (Y)$ yields $d f^i /d t = 0$, i.e. $f^1 (t) = f^1 (b) = 0$
and $f^j (t) = f^j (b) = a^j$, $2 \leq j \leq 2n$, hence
\[ X_{\gamma (t)} = \sum_{j=2}^{2n} a^j Y^H_{j, \gamma (t)} =
\sum_{j=2}^{2n} a^j \{ Y_{j , \gamma (t)} - \theta (Y_j )_{\gamma
(t)} T_{\gamma (t)} \} = \]
\[ = \sum_{j=2}^{2n} a^j Y_{j, \gamma (t)} + (t-a) a^1 T_{\gamma
(t)} = a^i Y_{i, \gamma (t)} = Y_{\gamma (t)} \, . \] Q.e.d.
\vskip 0.1in Setting $Y = 0$ in Proposition \ref{p:conj1} leads to
\begin{corollary}
Let $(M , \theta )$ be a Sasakian manifold and $\gamma : [a,b] \to
M$ a lengthy geodesic of the Tanaka-Webster connection,
parametrized by arc length and such that $\gamma (a)$ has no
conjugate point along $\gamma$. If $X$ is a piecewise
differentiable vector field along $\gamma$ such that $X_{\gamma
(a)} = X_{\gamma (b)} = 0$ and $X$ is perpendicular to $\gamma$
then $I_a^b (X) \geq 0$ and equality holds if and only if $X = 0$.
\label{c:conj1}
\end{corollary}
Corollary \ref{c:conj1} admits the following application
\begin{theorem} Let $(M , \theta )$ be a Sasakian manifold and
$\nabla$ its Tanaka-Webster connection. Assume that the
pseudohermitian sectional curvature satisfies $k_\theta (\sigma )
\geq k_0 > 0$, for any $2$-plane $\sigma \subset T_x (M)$, $x \in
M$. Then for any lengthy geodesic $\gamma (t) \in M$ of $\nabla$,
parametrized by arc length, the distance between two consecutive
conjugate points of $\gamma$ is less equal than $\pi
/(2\sqrt{k_0})$. \label{t:conj1}
\end{theorem}
{\em Proof}. Let $\gamma : [a, c] \to M$ be a geodesic of
$\nabla$, parametrized by arc length, such that $\gamma (c)$ is
the first conjugate point of $\gamma (a)$ along $\gamma$. Let $b
\in (a, c)$ and let $Y$ be a unit vector field along $\gamma$ such
that $(\nabla_{\dot{\gamma}} Y)_{\gamma (t)} = 0$ and $Y$ is
perpendicular to $\gamma$. Let $f(t)$ be a nonzero smooth function
such that $f(a) = f(b) = 0$. Then we may apply Corollary
\ref{c:conj1} to the vector field $X = f Y$ so that
\[ 0 \leq I_a^b (X) = \int_a^b \{ f^\prime (t)^2 |Y|^2 - f(t)^2 g_\theta (R(Y,
\dot{\gamma})\dot{\gamma} , Y) \} d t = \] \[ =  \int_a^b \{
f^\prime (t)^2 - 4 f(t)^2 k_\theta (\sigma )\} d t \leq \int_a^b
\{ f^\prime (t)^2 - 4 k_0 f(t)^2 \} d t \] where $\sigma \subset
T_{\gamma (t)} (M)$ is the $2$-plane spanned by $\{ Y_{\gamma (t)}
, \dot{\gamma}(t) \}$. Finally, we may choose $f(t) = \sin [\pi
(t-a)/(b-a)]$ and use $\int_0^\pi \cos^2 x \; d x = \int_0^\pi
\sin^2 x \; d x = \pi /2$. We get $b-a \leq \pi /\sqrt{4k_0}$ and
let $b \to c$. Q.e.d. \vskip 0.1in We may establish the following
more general version of Theorem \ref{t:conj1}
\begin{theorem} \label{t:conj2} Let $(M, \theta )$ be a Sasakian
manifold of CR dimension $n$ such that the Ricci tensor $\rho$ of
the Tanaka-Webster connection $\nabla$ satisfies \[ \rho (X,X)
\geq (2n-1) k_0 g_\theta (X,X), \;\;\; X \in H(M), \] for some
constant $k_0 > 0$. Then for any geodesic $\gamma$ of $\nabla$,
parametrized by arc length, the distance between any two
consecutive conjugate points of $\gamma$ is less than $\pi
/\sqrt{k_0}$.
\end{theorem}
\vskip 0.1in {\bf Remark}. The assumption on $\rho$ in Theorem
\ref{t:conj2} involves but the pseudohermitian Ricci curvature.
Indeed (cf. (1.98) in \cite{kn:DrTo}, section 1.4)
\[ {\rm Ric}(T_\alpha , T_{\overline{\beta}}) =
g_{\alpha\overline{\beta}} - \frac{1}{2} \;
R_{\alpha\overline{\beta}} \, , \]
\[ R_{\alpha\beta} = i(n-1) A_{\alpha\beta} \, , \;\;\; R_{0\beta}
= S^{\overline{\alpha}}_{\overline{\alpha}\beta} \, , \;\;\;
R_{\alpha 0} = R_{00} = 0, \] hence (by $\tau = 0$) $\rho (X,X) =
2 R_{\alpha\overline{\beta}} Z^\alpha Z^{\overline{\beta}}$, for
any $X = Z^\alpha T_\alpha + Z^{\overline{\alpha}}
T_{\overline{\alpha}} \in H(M)$. Here ${\rm Ric}$ is the Ricci
tensor of the Riemannian manifold $(M , g_\theta )$ (whose
symmetry yields $R_{\alpha\overline{\beta}} =
R_{\overline{\beta}\alpha}$). Note that $S = 0$ alone implies $T
\, \rfloor \, \rho = 0$. Also, if $(M, g_\theta )$ is Ricci flat
then $(M , \theta )$ is pseudo-Einstein (of pseudohermitian scalar
curvature $R = 2$), in the sense of \cite{kn:Lee}. \vskip 0.1in
{\em Proof of Theorem} \ref{t:conj2}. Let $\gamma (t) \in M$ as in
the proof of Theorem \ref{t:conj1}. Let $\{ Y_1 , \cdots ,
Y_{2n-1} \}$ be parallel (i.e. $(\nabla_{\dot{\gamma}} Y_i
)_{\gamma (t)} = 0$) vector fields such that $Y_i \in H(M)$ and
$\{ \dot{\gamma}(t) , Y_{1 , \gamma (t)} , \cdots , Y_{2n-1,
\gamma (t)} \}$ is an orthonormal basis of $H(M)_{\gamma (t)}$ for
every $t$. Let $f(t)$ be a nonzero smooth function such that $f(a)
= f(b) = 0$ and let us set $X_i = f Y_i$. Then (by Corollary
\ref{c:conj1})
\[ 0 \leq \sum_{i=1}^{2n-1} I_a^b (X_i ) = \sum_{i=1}^{2n-1} \int_a^b \{ f^\prime
(t)^2 |Y_i |^2 - f(t)^2 g_\theta (R(Y_i ,
\dot{\gamma})\dot{\gamma} , Y_i )\} d t = \] \[ = \int_a^b \{
(2n-1) f^\prime (t)^2 - f(t)^2 \rho (\dot{\gamma} , \dot{\gamma})
\} d t \leq \] \[ \leq (2n-1) \int_a^b \{ f^\prime (t)^2 - k_0
f(t)^2 \} dt
\] and the proof may be completed as that of Theorem
\ref{t:conj1}. \vskip 0.1in {\bf Remark}. The assumption in
Theorem \ref{t:conj2} is weaker than that in Theorem
\ref{t:conj1}. Indeed, let $X \in H(M)$, $X \neq 0$, and $V  =
|V|^{-1} V$. Let $\{ X_j : 1 \leq j \leq 2n \}$ be a local
orthonormal frame of $H(M)$ and $\sigma_j \subset T_x (M)$ the
$2$-plane spanned by $\{ Y_{j,x} , X_x \}$, where $Y_j := X_j -
g_\theta (V, X_j ) V$. Then $k_\theta (\sigma_j ) = \frac{1}{4}
g_{\theta} (R(V_j , V)V , V_j )_x$ where $V_j = |Y_j |^{-1} Y_j$
and $k_\theta (\sigma_j ) \geq k_0 /4$ yields
\[ \rho (X,X)_x = 4 |X|_x^2 \sum_{j=1}^{2n} k_\theta (\sigma_j )
|Y_j|^2_x \geq (2n-1) k_0 |X|^2_x \, . \]

\vskip 0.1in As another application of Proposition \ref{p:conj1}
we establish
\begin{theorem} Let $(M , \theta )$ be a Sasakian manifold, of CR
dimension $n$. Let $\gamma : [a,b] \to M$ be a lengthy geodesic of
the Tanaka-Webster connection $\nabla$, parametrized by arch
length. Assume that {\rm i)} there is $c \in (a,b)$ such that the
points $\gamma (a)$ and $\gamma (c)$ are horizontally conjugate
along $\gamma$ and {\rm ii)} for any $\delta > 0$ such that
$[c-\delta , c + \delta ] \subset (a,b)$ one has $\dim_{\mathbb R}
{\mathcal H}_{\gamma_\delta} = 4n$, where $\gamma_\delta$ is the
restriction of $\gamma$ to $[c-\delta , c+ \delta ]$. Then there
is a piecewise differentiable horizontal vector field $X$ along
$\gamma$ such that {\rm 1)} $X$ is perpendicular to $\dot{\gamma}$
and $J \dot{\gamma}$, {\rm 2)} $X_{\gamma (a)} = X_{\gamma (b)} =
0$, and {\rm 3)} $I_a^b (X) < 0$. \label{t:conj3}
\end{theorem}
In general, we have
\begin{lemma} Let $(M, \theta )$ be a Sasakian manifold of CR
dimension $n$ and $\gamma (t) \in M$ a lengthy geodesic of
$\nabla$, parametrized by arch length. Then \[ 2n+1 \leq
\dim_{\mathbb R} {\mathcal H}_\gamma \leq 4n. \] \label{l:conj6}
\end{lemma}
Hence the hypothesis in Theorem \ref{t:conj3} is that ${\mathcal
H}_{\gamma_\delta}$ has maximal dimension. We shall prove Lemma
\ref{l:conj6} later on. As to the converse of Theorem
\ref{t:conj3}, Corollary \ref{c:conj1} guarantees only that the
existence of a piecewise differentiable vector field $X$ as above
implies that there is some point $\gamma (c)$ conjugate to $\gamma
(a)$ along $\gamma$.
\par
{\em Proof of Theorem} \ref{t:conj3}. Let $a < c < b$ such that
$\gamma (a)$ and $\gamma (c)$ are horizontally conjugate and let
$Y \in {\mathcal H}_\gamma$ such that $Y_{\gamma (a)} = Y_{\gamma
(c)} = 0$. By Corollary \ref{c:ppdue} (as $(M, \theta )$ is
Sasakian) $Y$ is perpendicular to $\gamma$. Let $(U, x^i )$ be a
normal (with respect to $\nabla$) coordinate neighborhood with
origin at $\gamma (c)$. By Theorem 8.7 in \cite{kn:KoNo}, Vol. I,
p. 149, there is $R > 0$ such that for any $0 < r < R$ the open
set
\[ U(\gamma (c); r) \equiv \{ y \in U : \sum_{i=1}^{2n+1} x^i (y)^2
< r^2 \} \] is convex\footnote{That is any two points of $U(\gamma
(c); r)$ may be joined by a geodesic of $\nabla$ lying in
$U(\gamma (c); r)$.} and each point of $U(\gamma (c); r)$ has a
normal coordinate neighborhood containing $U(\gamma (c); r)$. By
continuity there is $\delta > 0$ such that $\gamma (t) \in
U(\gamma (c); r)$ for any $c - \delta \leq t \leq c + \delta$. Let
$\gamma_\delta$ denote the restriction of $\gamma$ to the interval
$[c-\delta , c+ \delta ]$. We need the following
\begin{lemma} The points $\gamma (c \pm \delta )$ are not
conjugate along $\gamma_\delta$. \label{l:conj5}
\end{lemma}
The proof is by contradiction. If $\gamma (c + \delta )$ is
conjugate to $\gamma (c - \delta )$ along $\gamma_\delta$ then (by
Theorem 1.4 in \cite{kn:KoNo}, Vol. II, p. 67) there is $v \in
T_{\gamma (c-\delta )}(M)$ such that $\exp_{\gamma (c-\delta )} v
= \gamma (c + \delta )$ and the linear map
\[ d_v \exp_{\gamma (c -\delta )} : T_v (T_{\gamma (c-\delta
)}(M)) \to T_{\gamma (c + \delta )} (M) \] is singular, i.e. ${\rm
Ker}(d_v \exp_{\gamma (c - \delta )} ) \neq 0$. Yet $\gamma (c -
\delta ) \in U(\gamma (c); r)$ hence there is a normal (relative
to $\nabla$) coordinate neighborhood $V$ with origin at $\gamma (c
- \delta )$ such that $V \supseteq U(\gamma (c); r)$. In
particular $\exp_{\gamma (c - \delta )} : V \to M$ is a
diffeomorphism on its image, so that $d_v \exp_{\gamma (c -\delta
)}$ is a linear isomorphism, a contradiction. Lemma \ref{l:conj5}
is proved.
\par
Let us go back to the proof of Theorem \ref{t:conj3}. The linear
map
\[ \Phi : J_{\gamma_\delta} \to T_{\gamma (c - \delta )}(M) \oplus
T_{\gamma (c + \delta )}(M), \;\;\; Z \mapsto (Z_{\gamma (c-\delta
)} \, , \, Z_{\gamma (c+ \delta )}), \] is a monomorphism. Indeed
${\rm Ker}(\Phi ) = 0$, otherwise $\gamma (c \pm \delta )$ would
be conjugate (in contradiction with Lemma \ref{l:conj5}). Both
spaces are $(4n+2)$-dimensional so that $\Phi$ is an epimorphism,
as well. By hypothesis ${\mathcal H}_{\gamma_\delta}$ is
$4n$-dimensional hence $\Phi$ descends to an isomorphism
\[ {\mathcal H}_{\gamma_\delta} \approx H(M)_{\gamma (c-\delta )}
\oplus H(M)_{\gamma (c + \delta )} \, . \] Let then $Z \in
{\mathcal H}_{\gamma_\delta}$ be a horizontal Jacobi field such
that
\[ Z_{\gamma (c - \delta )} = Y_{\gamma (c - \delta )}, \;\;\;
Z_{\gamma (c + \delta )} = 0. \] We set
\[ X = \begin{cases} Y & {\rm on} \;\; \left. \gamma \right|_{[a,
c-\delta ]}, \cr Z & {\rm on} \;\; \gamma_\delta ,  \cr 0 & {\rm
on} \;\; \left. \gamma \right|_{[c+\delta , b]}. \cr \end{cases}
\]
By the very definition $X$ is horizontal, i.e. $X_{\gamma (t)} \in
H(M)_{\gamma (t)}$ for every $t$. Moreover (by ${\mathcal
J}_\gamma Y = 0$ and $\theta (Y) = 0$)
\[ I^c_a (Y ) = \int_a^c \{ |\nabla_{\dot{\gamma}} Y |^2 - g_\theta (R(Y ,
\dot{\gamma})\dot{\gamma} , Y) \} d t = \]
\[ = \int_a^c \{ |\nabla_{\dot{\gamma}} Y |^2 + g_\theta
(\nabla_{\dot{\gamma}}^2 Y, Y)\} d t = \]
\[ = g_\theta (\nabla_{\dot{\gamma}} Y , Y )_{\gamma (c)} -
g_\theta (\nabla_{\dot{\gamma}} Y , Y )_{\gamma (a)}  = 0 \] i.e.
$I_a^{c-\delta}(Y ) = - I_{c-\delta}^c (Y )$. Hence
\[ I_a^b (X) = I_a^{c-\delta}(Y ) + I_{c-\delta}^{c+\delta}(Z) = -
I_{c-\delta}^c (Y ) + I_{c-\delta}^{c+\delta}(Z) . \] Finally, let
us consider the vector field along $\gamma_\delta$
\[ W = \begin{cases} Y & {\rm on} \;\; \left. \gamma
\right|_{[c-\delta , c]} , \cr 0, & {\rm on} \;\; \left. \gamma
\right|_{[c , c + \delta ]} . \cr \end{cases} \] Note that
$W_{\gamma (c+\delta )} = 0$, $W_{\gamma (c-\delta )} = Z_{\gamma
(c - \delta )}$ and $W$ is perpendicular to $\gamma$. Thus we may
apply Proposition \ref{p:conj1} to $W$ and to $Z \in {\mathcal
H}_{\gamma_\delta}$ to conclude that $I_{c-\delta}^c (Y) =
I_{c-\delta}^{c+\delta}(W) \geq I_{c-\delta}^{c+\delta}(Z)$.
Consequently $I_a^b (X) < 0$. Let us show that $X$ is orthogonal
to $J \dot{\gamma}$. By Lemma \ref{l:J2} (as $Y \in J_\gamma$)
\[ \theta (Y^\prime )_{\gamma (t)} - 2 \Omega (Y ,
\dot{\gamma})_{\gamma (t)} = {\rm const}. = \theta (Y^\prime
)_{\gamma (a)} - 2 \Omega (Y, \dot{\gamma})_{\gamma (a)} \] hence
(as $Y_{\gamma (a)} = 0$ and $Y_{\gamma (t)} \in H(M)_{\gamma (t)}
\Longrightarrow Y^\prime_{\gamma (t)} \in H(M)_{\gamma (t)}$)
\[ 2 \Omega (Y, \dot{\gamma})_{\gamma (t)} = \theta (Y^\prime
)_{\gamma (t)} - \theta (Y^\prime )_{\gamma (a)} = 0 \] for any $a
\leq t \leq c-\delta$. Similarly (as $Z_{\gamma (c+\delta )} = 0$
and $Z$ is horizontal) $\Omega (Z, \dot{\gamma})_{\gamma (t)} = 0$
for any $c-\delta \leq t \leq c + \delta$. Therefore $\Omega (X,
\dot{\gamma})_{\gamma (t)} = 0$ for every $t$. Theorem
\ref{t:conj3} is proved.
\par
It remains that we prove Lemma \ref{l:conj6}. Let $\gamma (t) \in
M$, $|t| < \epsilon$, be a lengthy geodesic of $\nabla$. Let $X
\in {\mathcal H}_\gamma$ and $\{ Y_j : 1 \leq j \leq 4n+2 \}$ a
linear basis in $J_\gamma$. Then $X = c^j Y_j = c^j Y^H_j + c^j
\theta (Y_j ) T$ (where $Y^H_j \equiv Y_j - \theta (Y_j ) T$) for
some $c^j \in {\mathbb R}$. As $X_{\gamma (t)} \in H(M)_{\gamma
(t)}$ one has i) $c^j \theta (Y_j )_{\gamma (t)} = 0$ on one hand,
and ii) $c^j f_j^a (\gamma (t)) = f^a (\gamma (t))$, $1 \leq a
\leq 2n$, on the other, where $X = f^a X_a$, $Y_j^H = f_j^a X_a$
and $\{ X_a : 1 \leq a \leq 2n \}$ is a local frame of $H(M)$. One
may think of (i)-(ii) as a linear system in the unknowns $c^j$.
Let $r(t)$ be its rank. Then $\dim_{\mathbb R} {\mathcal H}_\gamma
= 4n+2 - r(t) \geq 2n+1$. To prove the remaining inequality in
Lemma \ref{l:conj6} it suffices to observe that ${\mathcal
H}_\gamma$ is contained in the space of all solutions to
$X^{\prime\prime} + R(X , \dot{\gamma}) \dot{\gamma} = 0$ obeying
$X_{\gamma (0)} \in H(M)_{\gamma (0)}$ and $X^\prime_{\gamma (0)}
\in H(M)_{\gamma (0)}$, which is $4n$-dimensional.

\section{The first variation of the length integral}
Let $M$ be a strictly pseudoconvex CR manifold and $y, z \in M$.
Let $\Gamma$ be the set of all piecewise differentiable curves
$\gamma : [a,b] \to M$ parametrized proportionally to arc length,
such that $\gamma (a) = y$ and $\gamma (b) = z$. As usual, for
each $\gamma \in \Gamma$ we let $T_\gamma (\Gamma )$ be the space
of all piecewise differentiable vector fields along $\gamma$ such
that $X_y = X_z = 0$. Given $X \in T_\gamma (\Gamma )$ let
$\gamma^s : [a,b] \to M$, $|s| < \epsilon$, be a family of curves
such that i) $\gamma^s \in \Gamma$, $|s| < \epsilon$, ii)
$\gamma^0 = \gamma$, iii) there is a partition $a = t_0 < t_1 <
\cdots < t_k = b$ such that the map $(t,s) \mapsto \gamma^s (t)$
is differentiable on each rectangle $[t_j , t_{j+1}] \times
(-\epsilon , \epsilon )$, $0 \leq j \leq k-1$, and iv) for each
fixed $t \in [a,b]$ the tangent vector to \[ \sigma_t : (-\epsilon
, \epsilon ) \to M, \;\;\; \sigma_t (s) = \gamma^s (t), \;\;\; |s|
< \epsilon , \] at the point $\gamma (t)$ is $X_{\gamma (t)}$. We
set as usual
\[ (d_\gamma L) X = \frac{d}{d s} \left\{ L(\gamma^s ) \right\}_{s=0} \, . \]
Here $L(\gamma^s )$ is the Riemannian length of $\gamma^s$ with
respect to the Webster metric $g_\theta$ (so that $\gamma^s$ need
not be lengthy to start with). One scope of this section is to
establish the following
\begin{theorem} Let $\gamma^s : [a,b] \to M$, $|s| < \epsilon$, be a $1$-parameter
family of curves such that $(t,s) \mapsto \gamma^s (t)$ is
differentiable on $[a,b] \times (-\epsilon , \epsilon )$ and each
$\gamma^s$ is parametrized proportionally to arc length. Let us
set $\gamma = \gamma^0$. Then
\begin{equation}
\frac{d}{d s} \{ L(\gamma^s ) \}_{s=0} = \frac{1}{r} \{ g_\theta
(X , \dot{\gamma})_{\gamma (b)} - g_\theta (X ,
\dot{\gamma})_{\gamma (a)} - \label{e:V1}
\end{equation}
\[ - \int_a^b [g_\theta (X, \nabla_{\dot{\gamma}} \dot{\gamma}) -
g_\theta (T_\nabla (X , \dot{\gamma}), \dot{\gamma})]_{\gamma (t)}
\; d t \} \] where $X_{\gamma (t)} = \dot{\sigma}_t (0)$, $a \leq
t \leq b$, and $r = |\dot{\gamma}(t)|$ is the common length of all
tangent vectors along $\gamma$. \label{t:V1}
\end{theorem}
This will be shortly seen to imply
\begin{theorem} Let $\gamma \in \Gamma$ and $X \in T_\gamma
(\Gamma )$. Let $a = c_0 < c_1 < \cdots < c_h < c_{h+1} = b$ be a
partition such that $\gamma$ is differentiable on each $[c_j ,
c_{j+1}]$, $0 \leq j \leq h$. Then
\begin{equation}
(d_\gamma L) X = \frac{1}{r} \{ \sum_{j=1}^h g_{\theta , \gamma
(c_j )} (X_{\gamma (c_j )} , \dot{\gamma}(c_j^{-}) -
\dot{\gamma}(c_j^{+})) - \label{e:V2}
\end{equation}
\[ - \int_a^b [g_\theta (X , \nabla_{\dot{\gamma}} \dot{\gamma}) -
g_\theta (T_\nabla (X , \dot{\gamma}), \dot{\gamma}) ]_{\gamma
(t)} \; d t \} \] where $\dot{\gamma}(c_j^{\pm}) = \lim_{t \to
c_j^{\pm}} \dot{\gamma}(t)$. \label{t:V2}
\end{theorem}
Consequently, we shall prove
\begin{corollary} A lengthy curve $\gamma \in \Gamma$ is a geodesic
of the Tanaka-Webster connection if and only if
\begin{equation}\label{e:V3} (d_\gamma L) X = \frac{1}{r} \int_a^b \theta
(X)_{\gamma (t)} A(\dot{\gamma} , \dot{\gamma})_{\gamma (t)} \; d
t \end{equation} for all $X \in T_\gamma (\Gamma )$. In
particular, if $(M , \theta )$ is a Sasakian manifold then lengthy
geodesics belonging to $\Gamma$ are the critical points of $L$ on
$\Gamma$. \label{c:V1}
\end{corollary}
The remainder of this section is devoted to the proofs of the
results above. We adopt the principal bundle approach in
\cite{kn:KoNo}, Vol. II, p. 80-83. The proof is a {\em verbatim}
transcription of the arguments there, except for the presence of
torsion terms.
\par
Let $\pi : O(M , g_\theta ) \to M$ be the ${\rm O}(2n+1)$-bundle
of $g_\theta$-orthonormal frames tangent to $M$. Let $Q =  [a,b]
\times (-\epsilon , \epsilon )$. Let $f : Q \to O(M, g_\theta )$
be a parametrized surface in $O(M, g_\theta )$ such that i) $\pi
(f(t,s)) = \gamma^s (t)$, $(t,s) \in Q$, and ii) $f^0 : [a,b] \to
O(M, g_\theta )$, $f^0 (t) = f(t,0)$, $a \leq t \leq b$, is a
horizontal curve. Precisely, the Tanaka-Webster connection
$\nabla$ of $(M , \theta )$ induces an infinitesimal connection in
the principal bundle ${\rm GL}(2n+1, {\mathbb R}) \to L(M) \to M$
(of all linear frames tangent to $M$) descending (because of
$\nabla g_\theta = 0$) to a connection $H$ in ${\rm O}(2n+1) \to
O(M , g_\theta ) \to M$. The requirement is that $(d f^0 /dt) (t)
\in H_{f^0 (t)}$, $a \leq t \leq b$.
\par
Let ${\mathbb S} , {\mathbb T} \in {\mathcal X}(Q)$ be given by
${\mathbb S} = \partial /\partial s$ and ${\mathbb T} = \partial
/\partial t$. Let \[ \mu \in \Gamma^\infty (T^* (O(M, g_\theta ))
\otimes {\mathbb R}^{2n+1}), \;\;\; \Theta = D \mu \, , \] \[
\omega \in \Gamma^\infty (T^* (O(M, g_\theta )) \otimes {\bf
o}(2n+1)), \;\;\; \Omega = D \omega \, , \] be respectively the
canonical $1$-form, the torsion $2$-form, the connection $1$-form,
and the curvature $2$-form of $H$ on $O(M, g_\theta )$. We denote
by
\[ \mu^* = f^* \mu , \;\;\; \Theta^* = f^* \Theta , \;\;\;
\omega^* = f^* \omega , \;\;\; \Omega^* = f^* \Omega , \] the
pullback of these forms to the rectangle $Q$. We claim that
\begin{equation}
[{\mathbb S} , {\mathbb T}] = 0, \label{e:A}
\end{equation}
\begin{equation}
\omega^* ({\mathbb T})_{(t,0)} = 0, \;\;\; a \leq t \leq b.
\label{e:B}
\end{equation}
Indeed (\ref{e:A}) is obvious. To check (\ref{e:B}) one needs to
be a bit pedantic and introduce the injections
\[ \alpha^s : [a,b] \to Q, \;\;\; \beta_t : (-\epsilon , \epsilon
) \to Q, \]
\[ \alpha^s (t) = \beta_t (s) = (t,s), \;\;\; a \leq t \leq b,
\;\;\; |s| < \epsilon , \] so that $f^0 = f \circ \alpha^0$. Then
\[ H_{f^0 (t)} \ni
\left. \frac{d f^0}{d t}(t) = (d_{(t,0)} f)(d_t \alpha^0 )
\frac{d}{d t} \right|_t = (d_{(t,0)} f) {\mathbb T}_{(t,0)} \, ,
\]
\[ \omega^* ({\mathbb T})_{(t,0)} = \omega_{f(t,0)} ((d_{(t,0)} f)
{\mathbb T}_{(t,0)}) = 0. \] Next, we claim that
\begin{equation}
{\mathbb S}(\mu^* ({\mathbb T})) = {\mathbb T}(\mu^* ({\mathbb
S})) + \omega^* ({\mathbb T}) \cdot \mu^* ({\mathbb S}) - \omega^*
({\mathbb S}) \cdot \mu^* ({\mathbb T}) + 2 \, \Theta^* ({\mathbb
S},{\mathbb T}), \label{e:C}
\end{equation}
\begin{equation}
{\mathbb S}(\omega^* ({\mathbb T})) = {\mathbb T}(\omega^*
({\mathbb S})) + \omega^* ({\mathbb T}) \omega^* ({\mathbb S}) -
\omega^* ({\mathbb S}) \omega^* ({\mathbb T}) + 2 \; \Omega^*
({\mathbb S}, {\mathbb T}). \label{e:D}
\end{equation}
The identities (\ref{e:C})-(\ref{e:D}) follow from Prop. 3.11 in
\cite{kn:KoNo}, Vol. I, p. 36, our identity (\ref{e:A}), and the
first and second structure equations for a linear connection (cf.
e.g. Theor. 2.4 in \cite{kn:KoNo}, Vol. I, p. 120). Let us
consider the $C^\infty$ function $F : Q \to [0, + \infty )$ given
by
\[ F(t,s) = \langle \mu^* ({\mathbb T})_{(t,s)} \, , \,
\mu^* ({\mathbb T})_{(t,s)} \rangle^{1/2} \, , \;\;\; (t,s) \in Q.
\] Here $\langle \xi , \eta \rangle$ is the Euclidean scalar
product of $\xi , \eta \in {\mathbb R}^{2n+1}$. Note that
\[ \mu^* ({\mathbb T})_{(t,s)} = \mu_{f(t,s)}((d_{(t,s)} f)
{\mathbb T}_{(t,s)}) = \] \[ = f(t,s)^{-1} (d_{f(t,s)} \pi
)(d_{(t,s)} f) {\mathbb T}_{(t,s)} = \left.  f(t,s)^{-1} d_t (\pi
\circ f \circ \alpha^s ) \frac{d}{d t} \right|_t  \] i.e.
\begin{equation} \mu^* ({\mathbb T})_{(t,s)} = f(t,s)^{-1} \dot{\gamma}^s
(t). \label{e:I}
\end{equation}
Yet $f(t,s) \in O(M , g_\theta )$, i.e. $f(t,s)$ is a linear
isometry of $({\mathbb R}^{2n+1} , \langle \; , \; \rangle )$ onto
$(T_{\gamma^s (t)} (M), g_{\theta , \gamma^s (t)})$, so that
\[ F(t,s) = g_{\theta , \gamma^s (t)} (\dot{\gamma}^s (t) ,
\dot{\gamma}^s (t))^{1/2}  \] and then
\[ L(\gamma^s ) = \int_a^b F(t, s) \; d t. \]
As $\gamma^s$ is parametrized proportionally to arc length
$F(t,s)$ doesn't depend on $t$. In particular
\begin{equation} \label{e:E} F(t,0) = r.
\end{equation}
We claim that
\begin{equation}
{\mathbb S}(F) = \frac{1}{r} \{ \langle {\mathbb T}(\mu^*
({\mathbb S})) , \mu^* ({\mathbb T}) \rangle + 2 \langle \Theta^*
({\mathbb S},{\mathbb T}), \mu^* ({\mathbb T}) \rangle \}
\label{e:F}
\end{equation}
at all points $(t,0) \in Q$. Indeed, by (\ref{e:C})
\[ 2\, F\, {\mathbb S}(F) = {\mathbb S}(F^2 ) =
{\mathbb S}(\langle \mu^* ({\mathbb T}) ,  \mu^* ({\mathbb
T})\rangle ) = 2 \, \langle {\mathbb S}(\mu^* ({\mathbb T})) ,
\mu^* ({\mathbb T}) \rangle = \]
\[ = 2 \langle {\mathbb T}(\mu^* ({\mathbb S})) , \mu^*
({\mathbb T})\rangle + 2 \langle \omega^* ({\mathbb T}) \cdot
\mu^* ({\mathbb S}) , \mu^* ({\mathbb T}) \rangle - \] \[ - 2
\langle \omega^* ({\mathbb S}) \cdot \mu^* ({\mathbb T}) , \mu^*
({\mathbb T}) \rangle +  4 \langle \Theta^* ({\mathbb S} ,
{\mathbb T}), \mu^* ({\mathbb T}) \rangle . \] On the other hand
$\omega$ is ${\bf o}(2n+1)$-valued (where ${\bf o}(2n+1)$ is the
Lie algebra of ${\rm O}(2n+1)$), i.e. $\omega^* ({\mathbb
S})_{(t,s)} : {\mathbb R}^{2n+1} \to {\mathbb R}^{2n+1}$ is skew
symmetric, hence the last-but-one term vanishes. Therefore
(\ref{e:F}) follows from (\ref{e:B}) and (\ref{e:E}). We may
compute now the first variation of the length integral
\[ \frac{d}{d t} \{ L(\gamma^s ) \}_{s=0} = \int_a^b {\mathbb
S}(F)_{(t,0)} \; d t = \;\;\; ({\rm by} \; {\rm (\ref{e:F})}) \]
\[ = \frac{1}{r} \int_a^b \{ \langle {\mathbb T}(\mu^*
({\mathbb S})) , \mu^* ({\mathbb T}) \rangle_{(t,0)} + 2 \langle
\Theta^* ({\mathbb S},{\mathbb T}), \mu^* ({\mathbb T})
\rangle_{(t,0)} \} \; d t. \] On the other hand
\[ \mu^* ({\mathbb S})_{(t,0)} = \mu_{f^0 (t)} ((d_{(t,0)}
f) {\mathbb S}_{(t,0)}) = \] \[ = \left. f(t,0)^{-1} d_0 (\pi
\circ f \circ \beta_t ) \frac{d}{d s} \right|_0 = f(t,0)^{-1}
\frac{d \sigma_t}{d s}(0) \] i.e.
\begin{equation} \label{e:H} \mu^* ({\mathbb S})_{(t,0)} = f^0
(t)^{-1} X_{\gamma (t)} \, .
\end{equation}
Note that given $u \in C^\infty (Q)$ one has ${\mathbb
T}(u)_{(t,0)} = (u \circ \alpha^0 )^\prime (t)$. Then
\[ {\mathbb T}(\mu^* ({\mathbb T}))_{(t,0)} = \lim_{h \to 0}
\frac{1}{h} \{ \mu^* ({\mathbb T})_{(t+h,0)} - \mu^* ({\mathbb
T})_{(t,0)} \} = \;\;\; ({\rm by} \; {\rm (\ref{e:I})})
\]
\[ = \lim_{h \to 0} \frac{1}{h} \{ f^0 (t+h)^{-1}
\dot{\gamma}(t+h) - f^0 (t)^{-1} \dot{\gamma}(t) \} . \] Yet, as
$f^0$ is an horizontal curve
\[ f^0 (t+h)^{-1} \dot{\gamma}(t+h) = f^0 (t)^{-1} \tau^{t+h}_t
\dot{\gamma}(t+h), \] where $\tau^{t+h}_t : T_{\gamma (t+h)} (M)
\to T_{\gamma (t)} (M)$ is the parallel displacement operator
along $\gamma$ from $\gamma (t+ h)$ to $\gamma (t)$. Hence
\[ {\mathbb T}(\mu^* ({\mathbb T}))_{(t,0)}
= f^0 (t)^{-1} \left( \lim_{h \to 0} \frac{1}{h} \{ \tau^{t+h}_t
\dot{\gamma}(t+h) - \dot{\gamma}(t)\} \right) \]  i.e.
\begin{equation} \label{e:J}
{\mathbb T}(\mu^* ({\mathbb T}))_{(t,0)} = f^0 (t)^{-1}
(\nabla_{\dot{\gamma}} \dot{\gamma})_{\gamma (t)} \, .
\end{equation}
To compute the torsion term we recall (cf. \cite{kn:KoNo}, Vol. I,
p. 132)
\[ T_{\nabla , x} (X,Y) = 2 v (\Theta_v (X^* , Y^* )), \]
for any $X,Y \in T_x (M)$, where $v$ is a linear frame at $x$ and
$X^* , Y^* \in T_v (L(M))$ project respectively on $X, Y$. Note
that $(d_{(t,0)} f) {\mathbb S}_{(t,0)}$ and $(d_{(t,0)}
f){\mathbb T}_{(t,0)}$ project on $X_{\gamma (t)}$ and
$\dot{\gamma}(t)$, respectively. Then
\begin{equation}
2 \Theta^* ({\mathbb S}, {\mathbb T})_{(t,0)} = f^0 (t)^{-1}
T_\nabla (X , \dot{\gamma})_{\gamma (t)} \, . \label{e:P}
\end{equation}
Finally (by (\ref{e:I}) and (\ref{e:H})-(\ref{e:P}))
\[ \frac{d}{d t} \{ L(\gamma^s ) \}_{s=0} = \frac{1}{r} \int_a^b
\{ {\mathbb T}(\langle \mu^* ({\mathbb S}), \mu^* ({\mathbb T})
\rangle ) - \] \[ - \langle \mu^* ({\mathbb S}) , {\mathbb
T}(\mu^* ({\mathbb T})) \rangle + 2 \langle \Theta^* ({\mathbb S},
{\mathbb T}), \mu^* ({\mathbb T}) \}_{(t,0)} \; d t = \]
\[ = \frac{1}{r} \{ \langle \mu^* ({\mathbb S}), \mu^*
({\mathbb T}) \rangle_{(b,0)} - \langle \mu^* ({\mathbb S}), \mu^*
({\mathbb T}) \rangle_{(a,0)} \} - \]
\[ - \frac{1}{r} \int_a^b \{ g_\theta (X , \nabla_{\dot{\gamma}}
\dot{\gamma}) - g_\theta (T_\nabla (X , \dot{\gamma}),
\dot{\gamma}) \}_{\gamma (t)} \; d t \] and (\ref{e:V1}) is
proved.
\par
{\em Proof of Theorem} \ref{t:V2}. Let $c_j = t_0^{(j)} <
t_1^{(j)} < \cdots < t^{(j)}_{k_j} = c_{j+1}$ be a partition of
$[c_j , c_{j+1}]$ such that $X$ is differentiable along the
restriction of $\gamma$ at each $[t^{(j)}_i , t^{(j)}_{i+1}]$, $0
\leq i \leq k_j - 1$. Moreover, let $\{ \gamma^s \}_{|s| <
\epsilon}$ be a family of curves $\gamma^s \in \Gamma$ such that
$\gamma^0 = \gamma$, the map $(t,s) \mapsto \gamma^s (t)$ is
differentiable on $[c_j , c_{j+1}] \times (-\epsilon , \epsilon )$
for every $0 \leq j \leq h$, and $X_{\gamma (t)} = (d \sigma_t /d
s)(0)$ for every $t$ (with $\sigma_t (s) = \gamma^s (t)$). Let
$\gamma^s_j$ (respectively $\gamma^s_{ji}$) be the restriction of
$\gamma^s$ (respectively of $\gamma^s_j$) to $[c_j , c_{j+1}]$
(respectively to $[t^{(j)}_i , t^{(j)}_{i+1}]$). We may apply
Theorem \ref{t:V1} (to the interval $[t^{(j)}_i , t^{(j)}_{i+1}]$
rather than $[a,b]$) so that to get
\[ \frac{d}{ds} \{ L(\gamma^s_{ji} )\}_{s=0} = \frac{1}{r} \{
g_\theta (X , \dot{\gamma})_{\gamma (t^{(j)}_{i+1})} - g_\theta (X
, \dot{\gamma})_{\gamma (t^{(j)}_i )}  -
\int_{t^{(j)}_i}^{t^{(j)}_{i+1}} F(X , \dot{\gamma}) d t \} \]
where $F(X , \dot{\gamma})$ is short for $g_\theta (X ,
\nabla_{\dot{\gamma}} \dot{\gamma})_{\gamma (t)} - g_\theta
(T_\nabla (X , \dot{\gamma}) , \dot{\gamma})_{\gamma (t)}$. Let us
take the sum over $0 \leq i \leq k_j -1$. The lengths
$L(\gamma^s_{ji})$ ad up to $L(\gamma^s_j )$. Taking into account
that at the points $\gamma (c_j )$ only the lateral limits of
$\dot{\gamma}$ are actually defined, we obtain
\[ \frac{d}{d s} \{ L(\gamma^s_j )\}_{s=0} = \frac{1}{r} \{
g_{\theta , \gamma (c_{j+1})}(X_{\gamma (c_{j+1})} ,
\dot{\gamma}(c_{j+1}^{-})) - \] \[ - g_{\theta , \gamma
(c_{j})}(X_{\gamma (c_{j})} , \dot{\gamma}(c_{j}^{+})) -
\int_{c_j}^{c_{j+1}} F(X, \dot{\gamma}) d t \} \] and taking the
sum over $0 \leq j \leq h$ leads to (\ref{e:V2}) (as $X_{\gamma
(c_0 )} = 0$ and $X_{\gamma (c_{h+1})} = 0$). Q.e.d.
\par
{\em Proof of Corollary} \ref{c:V1}. Let $\gamma (t) \in M$ be a
lengthy curve such that $\gamma \in \Gamma$. If $\gamma$ is a
geodesic of $\nabla$ then $\nabla_{\dot{\gamma}} \dot{\gamma} = 0$
implies (by Theorem \ref{t:V1}) \[ (d_\gamma L )X = \frac{1}{r}
\int_a^b g_\theta (T_\nabla (X , \dot{\gamma}),
\dot{\gamma})_{\gamma (t)} d t \] for any $X \in T_\gamma (\Gamma
)$ and then
\[ T_\nabla (X, \dot{\gamma}) = - 2 \Omega (X_H , \dot{\gamma}) T
+ \theta (X) \tau (\dot{\gamma}), \;\;\; g_\theta (T ,
\dot{\gamma}) = 0, \] yield (\ref{e:V3}). Viceversa, let $\gamma
\in \Gamma$ be a lengthy curve such that (\ref{e:V3}) holds. There
is a partition $a = c_0 < c_1 < \cdots < c_{h+1} = b$ such that
$\gamma$ is differentiable in $[c_j , c_{j+1}]$, $0 \leq j \leq
h$. Let $f$ be a continuous function defined along $\gamma$ such
that $f(\gamma (c_j )) = 0$ for $1 \leq j \leq h$ and $f(\gamma
(t)) > 0$ elsewhere. We may apply (\ref{e:V2}) in Theorem
\ref{t:V2} to the vector field $X = f \; \nabla_{\dot{\gamma}}
\dot{\gamma}$ so that to get
\begin{equation}
(d_\gamma L) X = - \frac{1}{r} \int_a^b f \, \{
|\nabla_{\dot{\gamma}} \dot{\gamma}|^2 - g_\theta (T_\nabla
(\nabla_{\dot{\gamma}} \dot{\gamma} , \dot{\gamma}) ,
\dot{\gamma}) \} d t. \label{e:broken}
\end{equation}
As $\gamma$ is lengthy and $H(M)$ is parallel with respect to
$\nabla$ one has $\nabla_{\dot{\gamma}} \dot{\gamma} \in H(M)$
hence (by (\ref{e:V3})) $(d_\gamma L) (f \nabla_{\dot{\gamma}}
\dot{\gamma} ) = 0$ and
\[ g_\theta (T_\nabla (\nabla_{\dot{\gamma}} \dot{\gamma} ,
\dot{\gamma}), \dot{\gamma}) = - 2 \Omega (\nabla_{\dot{\gamma}}
\dot{\gamma} , \dot{\gamma}) g_\theta (T , \dot{\gamma}) = 0 \] so
that (by (\ref{e:broken})) it must be $\nabla_{\dot{\gamma}}
\dot{\gamma} = 0$ whenever $\nabla_{\dot{\gamma}} \dot{\gamma}$
makes sense, i.e. $\gamma$ is a broken geodesic of $\nabla$. It
remains that we prove differentiability of $\gamma$ at the points
$c_j$, $1 \leq j \leq h$. Let $j \in \{ 1 , \cdots , h \}$ be a
fixed index and let us consider a vector field $X_j \in T_\gamma
(\Gamma )$ such that $X_{j , \gamma (c_j )} = \dot{\gamma}(c_j^- )
- \dot{\gamma}(c_j^+ )$ and $X_{j, \gamma (c_k )} = 0$ for any $k
\in \{ 1 , \cdots h \} \setminus \{ j \}$. Then (by
(\ref{e:V2})-(\ref{e:V3})) one has $|\dot{\gamma}(c_j^- ) -
\dot{\gamma}(c_j^+ )|^2 = 0$. Q.e.d. \vskip 0.1in {\bf Remark}.
The following alternative proof of Theorem \ref{t:V1} is also
available. Since $(M , g_\theta )$ is a Riemannian manifold and
$L(\gamma^s )$ is the Rieman\-nian length of $\gamma^s$ we have
(cf. Theorem 5.1 in \cite{kn:KoNo}, Vol. II, p. 80)
\begin{equation}
\frac{d}{ds} \{ L(\gamma^s )\}_{s=0} = \frac{1}{r} \{ g_\theta (X
, \dot{\gamma})_{\gamma (b)} - g_\theta (X , \dot{\gamma})_{\gamma
(a)} \} - \label{e:fvf}
\end{equation}
\[ - \frac{1}{r} \; \int_a^b g_\theta (X , D_{\dot{\gamma}}
\dot{\gamma})_{\gamma (t)} \; d t  \] where $D$ is the Levi-Civita
connection of $(M , g_\theta )$. On the other hand (cf. e.g.
\cite{kn:DrTo}, section 1.3) $D$ is related to the Tanaka-Webster
connection of $(M , \theta )$ by $D = \nabla + (\Omega - A)
\otimes T + \tau \otimes \theta + 2 \theta \odot J$ hence
\[ g_\theta (X , D_{\dot{\gamma}} \dot{\gamma}) = g_\theta (X ,
\nabla_{\dot{\gamma}} \dot{\gamma}) - \theta (X) A(\dot{\gamma} ,
\dot{\gamma}) + \] \[ + \theta (\dot{\gamma}) A(X , \dot{\gamma})
+ 2 \theta (\dot{\gamma}) \Omega (X , \dot{\gamma}) = g_\theta (X
, \nabla_{\dot{\gamma}} \dot{\gamma}) - g_\theta (T_\nabla (X ,
\dot{\gamma}), \dot{\gamma}) \] so that (\ref{e:fvf}) yields
(\ref{e:V1}). Q.e.d.

\section{The second variation of the length integral}
We introduce the Hessian $I$ of $L$ at a geodesic $\gamma \in
\Gamma$ as follows. Given $X \in T_\gamma (\Gamma )$ let us
consider a $1$-parameter family of curves $\{ \gamma^s \}_{|s| <
\epsilon}$ as in the definition of $(d_\gamma L) X$. Let $I(X,X)$
be given by
\[ I(X,X) = \frac{d^2}{d s^2} \left\{ L(\gamma^s ) \right\}_{s=0}
\] and define $I(X,Y)$ by polarization. By analogy to Riemannian geometry
(cf. e.g. \cite{kn:KoNo}, Vol. II, p. 81) $I(X,Y)$ is referred to
as the {\em index form}. The scope of this section is to establish
\begin{theorem} Let $(M , \theta )$ be a Sasakian manifold.
If $\gamma \in \Gamma$ is a lengthy geodesic of the Tanaka-Webster
connection $\nabla$ of $(M , \theta )$ and $X,Y \in T_\gamma
(\Gamma )$ then \begin{equation}\label{e:index} I(X,Y) =
\frac{1}{r} \int_a^b \{ g_\theta (\nabla_{\dot{\gamma}} X^\bot ,
\nabla_{\dot{\gamma}} Y^\bot ) - g_\theta (R(X^\bot ,
\dot{\gamma}) \dot{\gamma} , Y^\bot ) -
\end{equation}
\[ - 2 \Omega (X^\bot , \dot{\gamma}) \theta (\nabla_{\dot{\gamma}}
Y^\bot ) - 2 [\theta (\nabla_{\dot{\gamma}} X^\bot ) - 2
\Omega(X^\bot , \dot{\gamma})] \Omega (Y^\bot , \dot{\gamma})\} d
t
\]
where $X^\bot = X - (1/r^2 )\, g_\theta (X , \dot{\gamma}) \,
\dot{\gamma}$. \label{t:9}
\end{theorem} We shall need the following reformulation of Theorem
\ref{t:9}
\begin{theorem} Let $(M , \theta )$, $\gamma$ and $X,Y$
be as in Theorem {\rm \ref{t:9}}. Then
\begin{equation}
I(X,Y) = - \frac{1}{r} \int_a^b \{ g_\theta ({\mathcal J}_\gamma
X^\bot , Y^\bot ) + 2 [\theta (\nabla_{\dot{\gamma}} X^\bot ) -
\label{e:index2}
\end{equation}
\[ - 2 \Omega (X^\bot , \dot{\gamma})] \Omega (Y^\bot ,
\dot{\gamma})\} d t + \]
\[ + \frac{1}{r} \sum_{j=1}^h g_{\theta , \gamma (t_j )}
((\nabla_{\dot{\gamma}} X^\bot )^-_{\gamma (t_j )} -
(\nabla_{\dot{\gamma}} X^\bot )^+_{\gamma (t_j )} \, , \,
Y^\bot_{\gamma (t_j )}) \] where ${\mathcal J}_\gamma X \equiv
\nabla^2_{\dot{\gamma}} X - 2 \Omega (X^\prime , \dot{\gamma}) T +
R(X, \dot{\gamma}) \dot{\gamma}$ is the Jacobi operator and $a =
t_0 < t_1 < \cdots < t_h < t_{h+1} = b$ is a partition of $[a,b]$
such that $X$ is differentiable in each interval $[t_j ,
t_{j+1}]$, $0 \leq j \leq h$, and $(\nabla_{\dot{\gamma}} X^\bot
)_{\gamma (t_j )}^{\pm} = \lim_{t \to t_j^{\pm}}
(\nabla_{\dot{\gamma}} X^\bot )_{\gamma (t)}$. \label{t:10}
\end{theorem}
This will be seen to imply
\begin{corollary} Let $(M , \theta )$ be a Sasakian manifold, $\gamma \in \Gamma$
a lengthy geodesic of the Tanaka-Webster connection of $(M ,
\theta )$, and $X \in T_\gamma (\Gamma )$. Then $X^\bot$ is a
Jacobi field if and only if there is $\alpha (X) \in {\mathbb R}$
such that
\begin{equation}\label{e:lJ2} \frac{d}{d t} \{ \theta (X^\bot )
\circ \gamma \} (t) - 2 \Omega (X^\bot , \dot{\gamma})_{\gamma
(t)} = \alpha (X)
\end{equation} for any $a \leq t \leq b$, and
\begin{equation}
\label{e:indlJ2}
I(X,Y) = - (2/r) \alpha (X) \int_a^b \Omega (Y^\bot ,
\dot{\gamma})_{\gamma (t)} \, d t \, ,
\end{equation}
for any $Y \in T_\gamma (\Gamma )$.
\label{c:4}
\end{corollary}

{\em Proof of Theorem} \ref{t:9}. We adopt the notations and
conventions in the proof of Theorem \ref{t:V1}. As a byproduct of
the proof of (\ref{e:F}) we have the identity
\begin{equation}
\frac{1}{2} \; {\mathbb S}(F^2 ) = \langle {\mathbb T}(\mu^*
({\mathbb S})), \mu^* ({\mathbb T}) \rangle +
\label{e:K}\end{equation}
\[ + \langle \omega^* ({\mathbb T}) \cdot \mu^* ({\mathbb S}) ,
\mu^* ({\mathbb T}) \rangle + 2 \langle \Theta^* ({\mathbb S},
{\mathbb T}), \mu^* ({\mathbb T})\rangle . \] Applying ${\mathbb
S}$ we get
\[ \frac{1}{2} \; {\mathbb S}^2 (F^2 ) = \langle {\mathbb S} {\mathbb
T}(\mu^* ({\mathbb S})), \mu^* ({\mathbb T}) \rangle + \langle
{\mathbb T}(\mu^* ({\mathbb S})), {\mathbb S}(\mu^* ({\mathbb T}))
\rangle + \]
\[ + \langle {\mathbb S}(\omega^* ({\mathbb T})) \cdot \mu^*
({\mathbb S}) , \mu^* ({\mathbb T})\rangle + \langle \omega^*
({\mathbb T}) \cdot {\mathbb S} (\mu^* ({\mathbb S})), \mu^*
({\mathbb T})\rangle + \] \[ +  \langle \omega^* ({\mathbb T})
\cdot \mu^* ({\mathbb S}) , {\mathbb S}(\mu^* ({\mathbb T}))
\rangle + \]
\[ + 2 \langle {\mathbb S}(\Theta^* ({\mathbb S},
{\mathbb T})) , \mu^* ({\mathbb T})\rangle + 2 \langle \Theta^*
({\mathbb S}, {\mathbb T}), {\mathbb S}(\mu^* ({\mathbb T}))
\rangle . \] When calculated at points of the form $(t,0) \in Q$
the $4^{\rm th}$ and $5^{\rm th}$ terms vanish (by (\ref{e:B})).
We proceed by calculating the remaining terms (at $(t,0)$). By
(\ref{e:A})
\[ 1^{\rm st} \; {\rm term} = \langle {\mathbb S} {\mathbb
T}(\mu^* ({\mathbb S})), \mu^* ({\mathbb T}) \rangle = \langle
{\mathbb T} {\mathbb S}(\mu^* ({\mathbb S})), \mu^* ({\mathbb T})
\rangle = \]
\[ = {\mathbb T} \left( \langle {\mathbb S}(\mu^* ({\mathbb S})), \mu^*
({\mathbb T}) \rangle \right) -  \langle {\mathbb S}(\mu^*
({\mathbb S})), {\mathbb T}(\mu^* ({\mathbb T})) \rangle . \] Yet
$\gamma \in \Gamma$ is a geodesic hence (by (\ref{e:J})) ${\mathbb
T}(\mu^* ({\mathbb T}))_{(t,0)} = 0$. Hence
\[ 1^{\rm st} \; {\rm term} =
{\mathbb T} \left( \langle {\mathbb S}(\mu^* ({\mathbb S})), \mu^*
({\mathbb T}) \rangle \right)_{(t,0)} . \] Next (by (\ref{e:C}))
\[ 2^{\rm nd} \; {\rm term} =
\langle {\mathbb T}(\mu^* ({\mathbb S})), {\mathbb S}(\mu^*
({\mathbb T})) \rangle = \langle {\mathbb T}(\mu^* ({\mathbb S})),
{\mathbb T}(\mu^* ({\mathbb S})) \rangle + \]
\[ + \langle {\mathbb T}(\mu^* ({\mathbb S})), \omega^* ({\mathbb T})\cdot
\mu^* ({\mathbb S}) \rangle - \langle {\mathbb T}(\mu^* ({\mathbb
S})), \omega^* ({\mathbb S})\cdot \mu^* ({\mathbb T}) \rangle +
\]
\[ + 2 \langle {\mathbb T}(\mu^* ({\mathbb S})), \Theta^* ({\mathbb
S}, {\mathbb T}) \rangle . \] Again terms are evaluated at $(t,0)$
hence $\omega^* ({\mathbb T}) = 0$ (by (\ref{e:B})). On the other
hand $\omega^* ({\mathbb S})$ is ${\bf o}(2n+1)$-valued hence
\[ 2^{\rm nd} \; {\rm term} = \langle {\mathbb T}(\mu^* ({\mathbb
S})), {\mathbb T}(\mu^* ({\mathbb S})) \rangle + \langle \omega^*
({\mathbb S})\cdot {\mathbb T}(\mu^* ({\mathbb S})), \mu^*
({\mathbb T}) \rangle + \]
\[ + 2 \langle {\mathbb T}(\mu^* ({\mathbb S})), \Theta^* ({\mathbb
S}, {\mathbb T}) \rangle  \] at each $(t, 0) \in Q$. Next (by
(\ref{e:D}))
\[ 3^{\rm rd} \; {\rm term} = \langle {\mathbb S}(\omega^* ({\mathbb T})) \cdot \mu^*
({\mathbb S}) , \mu^* ({\mathbb T})\rangle = \langle {\mathbb
T}(\omega^* ({\mathbb S})) \cdot \mu^* ({\mathbb S}) , \mu^*
({\mathbb T})\rangle + \]
\[ + \langle \omega^* ({\mathbb T}) \omega^* ({\mathbb S}) \cdot \mu^*
({\mathbb S}) , \mu^* ({\mathbb T})\rangle - \langle \omega^*
({\mathbb T}) \omega^* ({\mathbb S}) \cdot \mu^* ({\mathbb S}) ,
\mu^* ({\mathbb T})\rangle + \] \[ + 2 \,  \langle \Omega^*
({\mathbb S} , {\mathbb T}) \cdot \mu^* ({\mathbb S}) , \mu^*
({\mathbb T}) \rangle \] or (by (\ref{e:B}))
\[ 3^{\rm rd} \; {\rm term} = \langle {\mathbb T}(\omega^* ({\mathbb S})) \cdot \mu^*
({\mathbb S}) , \mu^* ({\mathbb T})\rangle + 2 \, \langle \Omega^*
({\mathbb S} , {\mathbb T}) \cdot \mu^* ({\mathbb S}) , \mu^*
({\mathbb T}) \rangle
\]
at each $(t,0) \in Q$. Finally (by (\ref{e:C}))
\[ 7^{\rm th} \; {\rm term} = 2 \langle \Theta^*
({\mathbb S}, {\mathbb T}), {\mathbb S}(\mu^* ({\mathbb T}))
\rangle = 2 \langle \Theta^* ({\mathbb S}, {\mathbb T}), {\mathbb
T}(\mu^* ({\mathbb S})) \rangle + \]
\[ + 2 \langle \Theta^*
({\mathbb S}, {\mathbb T}), \omega^* ({\mathbb T}) \cdot \mu^*
({\mathbb S}) \rangle - 2 \langle \Theta^* ({\mathbb S}, {\mathbb
T}), \omega^* ({\mathbb S}) \cdot \mu^* ({\mathbb T}) \rangle + 4
|\Theta^* ({\mathbb S} , {\mathbb T})|^2 \] or (by (\ref{e:B}) and
the fact that $\omega^* ({\mathbb S})$ is skew)
\[ 7^{\rm th} \; {\rm term} = 2 \langle \Theta^* ({\mathbb S}, {\mathbb T}), {\mathbb
T}(\mu^* ({\mathbb S})) \rangle + \]
\[ +  2 \langle \omega^*
({\mathbb S}) \cdot \Theta^* ({\mathbb S}, {\mathbb T}),  \mu^*
({\mathbb T}) \rangle + 4 |\Theta^* ({\mathbb S} , {\mathbb T})|^2
\]
at each $(t,0) \in Q$. Summing up the various expressions and
noting that (again by (\ref{e:J}))
\[ {\mathbb T}\left( \langle {\mathbb S}(\mu^* ({\mathbb S})) ,
\mu^* ({\mathbb T})\rangle \right) + \langle \omega^* ({\mathbb
S}) \cdot {\mathbb T}(\mu^* ({\mathbb S})) , \mu^* ({\mathbb T})
\rangle + \] \[ + \langle {\mathbb T}(\omega^* ({\mathbb S}))
\cdot \mu^* ({\mathbb S}) , \mu^* ({\mathbb T}) \rangle = \]
\[ = {\mathbb T} \left( \langle {\mathbb S}(\mu^* ({\mathbb
S})) + \omega^* ({\mathbb S}) \cdot \mu^* ({\mathbb S}) , \mu^*
({\mathbb T}) \rangle \right) \] we obtain
\begin{equation}
\label{e:M} \frac{1}{2} \; {\mathbb S}^2 (F^2 ) = |{\mathbb
T}(\mu^* ({\mathbb S}))|^2 + 2 \, \langle \Omega^* ({\mathbb S} ,
{\mathbb T}) \cdot \mu^* ({\mathbb S}) , \mu^* ({\mathbb T})
\rangle +
\end{equation}
\[ + {\mathbb T} \left( \langle {\mathbb S}(\mu^* ({\mathbb
S})) + \omega^* ({\mathbb S}) \cdot \mu^* ({\mathbb S}) , \mu^*
({\mathbb T}) \rangle \right) + \] \[ + 4 \langle \Theta^*
({\mathbb S}, {\mathbb T}), {\mathbb T}(\mu^* ({\mathbb S}))
\rangle + 2 \langle {\mathbb S}(\Theta^* ({\mathbb S}, {\mathbb
T})) , \mu^* ({\mathbb T})\rangle + \]
\[ +  2 \langle \omega^*
({\mathbb S}) \cdot \Theta^* ({\mathbb S}, {\mathbb T}),  \mu^*
({\mathbb T}) \rangle + 4 |\Theta^* ({\mathbb S} , {\mathbb T})|^2
\]
at each $(t,0) \in Q$. Since $F {\mathbb S}^2 (F) = \frac{1}{2}
{\mathbb S}^2 (F^2 ) - {\mathbb S}(F)^2$ we get (by (\ref{e:F})
and (\ref{e:M}))
\begin{equation}
\label{e:N} r {\mathbb S}^2 (F) = |{\mathbb T}(\mu^* ({\mathbb
S}))|^2 + 2 \, \langle \Omega^* ({\mathbb S} , {\mathbb T}) \cdot
\mu^* ({\mathbb S}) , \mu^* ({\mathbb T}) \rangle +
\end{equation}
\[ + {\mathbb T} \left( \langle {\mathbb S}(\mu^* ({\mathbb
S})) + \omega^* ({\mathbb S}) \cdot \mu^* ({\mathbb S}) , \mu^*
({\mathbb T}) \rangle \right) + \] \[ + 4 \langle \Theta^*
({\mathbb S}, {\mathbb T}), {\mathbb T}(\mu^* ({\mathbb S}))
\rangle + 2 \langle {\mathbb S}(\Theta^* ({\mathbb S}, {\mathbb
T})) , \mu^* ({\mathbb T})\rangle + \]
\[ +  2 \langle \omega^*
({\mathbb S}) \cdot \Theta^* ({\mathbb S}, {\mathbb T}),  \mu^*
({\mathbb T}) \rangle + 4 |\Theta^* ({\mathbb S} , {\mathbb T})|^2
- \]
\[ - \frac{1}{r^2} \{ \langle {\mathbb T}(\mu^* ({\mathbb S}))
, \mu^* ({\mathbb T})\rangle^2 + 4 \langle \Theta^* ({\mathbb S},
{\mathbb T}), \mu^* ({\mathbb T})\rangle^2 + \]
\[ + 4 \langle {\mathbb T}(\mu^* ({\mathbb S})) , \mu^*
({\mathbb T}) \rangle \; \langle \Theta^* ({\mathbb S}, {\mathbb
T}) , \mu^* ({\mathbb T}) \rangle \} \] at any $(t,0) \in Q$.
Moreover (by (\ref{e:H})) \[ {\mathbb T}(\mu^* ({\mathbb
S}))_{(t,0)} = \frac{d}{d t} \left\{ \mu^* ({\mathbb S}) \circ
\alpha^0 \right\} (t) =
\] \[ =  \lim_{h \to 0} \frac{1}{h} \{ \mu^* ({\mathbb S})_{(t+h, 0)}
- \mu^* ({\mathbb S})_{(t,0)} \} = \]
\[ = \lim_{h \to 0} \frac{1}{h} \{ f^0 (t+h)^{-1} X_{\gamma (t+h)}
- f^0 (t)^{-1} X_{\gamma (t)} \} = \] (as $f^0 : [a,b] \to O(M ,
g_\theta )$ is a horizontal curve)
\[ = \lim_{h \to 0} \frac{1}{h} \{ f^0 (t)^{-1} \tau^{t+h}_t
X_{\gamma (t+h)} - f^0 (t)^{-1} X_{\gamma (t)} \} = \]
\[ = f^0 (t)^{-1} \left( \lim_{h \to 0} \frac{1}{h} \{
\tau^{t+h}_t X_{\gamma (t+h)} - X_{\gamma (t)} \} \right) \] that
is
\begin{equation}
{\mathbb T}(\mu^* ({\mathbb S}))_{(t,0)} = f^0 (t)^{-1}
(\nabla_{\dot{\gamma}} X)_{\gamma (t)} \, . \label{e:O}
\end{equation}
Consequently (by (\ref{e:I}) and (\ref{e:O}))
\begin{equation}
\label{e:TS} |{\mathbb T}(\mu^* ({\mathbb S}))|^2 + 4 \langle
\Theta^* ({\mathbb S}, {\mathbb T}), {\mathbb T}(\mu^* ({\mathbb
S})) \rangle -
\end{equation}
\[ - \frac{1}{r^2} \{ \langle {\mathbb T}(\mu^* ({\mathbb S}))
, \mu^* ({\mathbb T})\rangle^2 + 4 \langle \Theta^* ({\mathbb S},
{\mathbb T}), \mu^* ({\mathbb T})\rangle^2 + \]
\[ + 4 \langle {\mathbb T}(\mu^* ({\mathbb S})) , \mu^*
({\mathbb T}) \rangle \; \langle \Theta^* ({\mathbb S}, {\mathbb
T}) , \mu^* ({\mathbb T}) \rangle \} = \]
\[ = |\nabla_{\dot{\gamma}} X|^2 + 2 g_\theta (T_\nabla (X,
\dot{\gamma}), \nabla_{\dot{\gamma}} X) - \frac{1}{r^2} \{
g_\theta (\nabla_{\dot{\gamma}} X , \dot{\gamma})^2 + \] \[ +
g_\theta (T_\nabla (X , \dot{\gamma}) , \dot{\gamma})^2 + 2
g_\theta (\nabla_{\dot{\gamma}} X , \dot{\gamma}) g_\theta
(T_\nabla (X , \dot{\gamma}), \dot{\gamma}) \} = \]
\[ = |\nabla_{\dot{\gamma}} X^\bot |^2 + 2 g_\theta (T_\nabla (X^\bot ,
\dot{\gamma}), \nabla_{\dot{\gamma}} X) - \] \[ - \frac{1}{r^2} \{
g_\theta (T_\nabla (X^\bot , \dot{\gamma}) , \dot{\gamma})^2 + 2
g_\theta (\nabla_{\dot{\gamma}} X , \dot{\gamma}) g_\theta
(T_\nabla (X^\bot , \dot{\gamma}), \dot{\gamma}) \} \] and (by
(\ref{e:I}) and (\ref{e:H}))
\begin{equation}
{\rm the \; curvature \; term} = 2 \, \langle \Omega^* ({\mathbb
S} , {\mathbb T}) \cdot \mu^* ({\mathbb S}) , \mu^* ({\mathbb T})
\rangle = \label{e:curv}
\end{equation}
\[ = g_\theta (R(X , \dot{\gamma}) X , \dot{\gamma})_{\gamma (t)}
= - g_\theta (R(X^\bot , \dot{\gamma}) \dot{\gamma} , X^\bot ). \]
On the other hand $\pi (f(a,s)) = y$ and $\pi (f(b,s)) = z$ imply
that $(d_{(a,s)} f) {\mathbb S}_{(a,s)}$ and $(d_{(b,s)} f)
{\mathbb S}_{(b,s)}$ are vertical hence
\begin{equation}
\mu^* ({\mathbb S})_{(a,s)} = 0, \;\;\; \mu^* ({\mathbb
S})_{(b,s)} = 0. \label{e:ab}
\end{equation}
Next, we wish to compute ${\mathbb S}(\mu^* ({\mathbb
S}))_{(t,0)}$. To do so we need to further specialize the choice
of $f(t,s)$. Precisely, let $v \in \pi^{-1} (\gamma (a))$ be a
fixed orthonormal frame and let
\begin{equation}
\label{e:lift} f(t,s) = \sigma_t^\uparrow (s), \;\;\; a \leq t
\leq b, \;\; |s| < \epsilon ,
\end{equation}
where $\sigma_t^\uparrow : (-\epsilon , \epsilon ) \to O(M ,
g_\theta )$ is the unique horizontal lift of $\sigma_t :
(-\epsilon , \epsilon ) \to M$ issuing at $\sigma_t (0) =
\gamma^\uparrow (t)$. Also $\gamma^\uparrow : [a,b] \to O(M ,
g_\theta )$ is the horizontal lift of $\gamma : [a,b] \to M$
determined by $\gamma^\uparrow (a) = v$. Therefore $f^0 =
\gamma^\uparrow$ is a horizontal curve, as required by the
previous part of the proof. In addition (\ref{e:lift}) possesses
the property that for each $t$ the curve $s \mapsto f(t,s)$ is
horizontal, as well. Then
\[ {\mathbb S}(\mu^* ({\mathbb S}))_{(t,0)} = \frac{d}{d s}
\left\{ \mu^* ({\mathbb S}) \circ \beta_t \right\} (0) = \] \[ =
\lim_{s \to 0} \frac{1}{s} \{ f(t,s)^{-1} \dot{\sigma}_t (s) -
f(t,0)^{-1} \dot{\sigma}_t (0) \} = \] (as $f_t : (-\epsilon ,
\epsilon ) \to O(M , g_\theta )$, $f_t (s) = f(t,s)$, $|s| <
\epsilon$, is horizontal)
\[ = \lim_{s \to 0} \frac{1}{s} \{ f(t,0)^{-1} \tau^s
\dot{\sigma}_t (s) - f(t,0)^{-1} \dot{\sigma}_t (0) \} = \] \[ =
f(t,0)^{-1} \left( \lim_{s \to 0} \frac{1}{s} \{ \tau^s
\dot{\sigma}_t (s) - \dot{\sigma}_t (0) \} \right) \] where
$\tau^s : T_{\sigma_t (s)} (M) \to T_{\sigma_t (0)}(M)$ is the
parallel displacement along $\sigma_t$ from $\sigma_t (s)$ to
$\sigma_t (0)$, i.e.
\begin{equation}
{\mathbb S}(\mu^* ({\mathbb S}))_{(t,0)} = f(t,0)^{-1}
(\nabla_{\dot{\sigma}_t} \dot{\sigma}_t )_{\gamma (t)} \, .
\label{e:SS}
\end{equation}
By (\ref{e:SS}), $\dot{\sigma}_a (s) = 0$ and $\dot{\sigma}_b (s)
= 0$ (as $\sigma_a (s) = \gamma^s (a) = y =$ const. and $\sigma_b
(s) = \gamma^s (b) = z =$ const.) it follows that
\begin{equation}
\label{e:ab1} {\mathbb S}(\mu^* ({\mathbb S}))_{(a,0)} = 0, \;\;\;
{\mathbb S}(\mu^* ({\mathbb S}))_{(b,0)} = 0.
\end{equation}
Using (\ref{e:ab}) and (\ref{e:ab1}) we may conclude that
\begin{equation}
\int_a^b {\mathbb T} \left( \langle {\mathbb S}(\mu^* ({\mathbb
S})) + \omega^* ({\mathbb S}) \cdot \mu^* ({\mathbb S}) , \mu^*
({\mathbb T}) \rangle \right)_{(t,0)} d t = \label{e:int}
\end{equation}
\[ = \langle {\mathbb S}(\mu^* ({\mathbb
S})) + \omega^* ({\mathbb S}) \cdot \mu^* ({\mathbb S}) , \mu^*
({\mathbb T}) \rangle_{(b,0)} - \] \[ - \langle {\mathbb S}(\mu^*
({\mathbb S})) + \omega^* ({\mathbb S}) \cdot \mu^* ({\mathbb S})
, \mu^* ({\mathbb T}) \rangle_{(a,0)} = 0. \] Similarly
\[ 2\; {\mathbb S}(\Theta^* ({\mathbb S}, {\mathbb T}))_{(t,0)} =  2 \;
\lim_{s \to 0} \frac{1}{s} \{ \Theta^* ({\mathbb S}, {\mathbb
T})_{(t,s)} - \Theta^* ({\mathbb S}, {\mathbb T})_{(t,0)} \} =
\]
\[ =  \lim_{s\to 0} \frac{1}{s} \{ f(t,s)^{-1} T_\nabla
(\dot{\sigma}_t (s) , \dot{\gamma}^s (t) ) - f(t,0)^{-1} T_\nabla
(\dot{\sigma}_t (0), \dot{\gamma}^0 (t)) \} =
\]
\[ = f(t,0)^{-1} \left( \lim_{s \to 0} \frac{1}{s} \{ \tau^s \,
V_{\sigma_t (s)} - V_{\sigma_t (0)} \} \right) \] where $V$ is the
vector field defined at each $a(t,s) = \sigma_t (s)$ by
\[ V_{a(t,s)} = T_{\nabla , a(t,s)} (\dot{\sigma}_t
(s) , \dot{\gamma}^s (t)), \;\;\; a \leq t \leq b, \;\; |s| <
\epsilon .
\] Let us assume from now on that $\tau = 0$, i.e. $(M , \theta )$ is Sasakian.
Then
\[ V_{a(t,s)} = - 2 \Omega_{a(t,s)} (\dot{\sigma}_t (s) ,
\dot{\gamma}^s (t)) T_{a(t,s)} \] and $\nabla T = 0$ yields
\[ \tau^s V_{a(t,s)} = - 2 \Omega_{a(t,s)} (\dot{\sigma}_t (s) ,
\dot{\gamma}^s (t)) T_{\gamma (t)} \, . \] Finally
\begin{equation} 2 \; \langle {\mathbb S}(\Theta^* ({\mathbb S}, {\mathbb
T})), \mu^* ({\mathbb T})\rangle_{(t,0)} =  \label{e:stor}
\end{equation}
\[ = \lim_{s \to 0}
\frac{1}{s} \{ g_{\theta , a(t,s)}(\tau^s V_{a(t,s)} ,
\dot{\gamma}(t)) - g_{\theta , a(t,0)} (V_{a(t,0)} ,
\dot{\gamma}(t)) \} = 0, \] as $\dot{\gamma}(t) \in H(M)_{\gamma
(t)}$. It remains that we compute the term $2 \, \langle \omega^*
({\mathbb S}) \cdot \Theta^* ({\mathbb S}, {\mathbb T}), \mu^*
({\mathbb T}) \rangle$. As $f(t,s)$ is a linear frame at $a(t,s)$
\[ f(t,s) = ( a(t,s), \{ X_{i, a(t,s)} : 1 \leq i \leq 2n+1 \} ),
\]
where $X_i \in T_{a(t,s)}(M)$. Let $(U , x^i )$ be a local
coordinate system on $M$ and let us set $X_i = X^j_i \; \partial
/\partial x^j$. Let $(\Pi^{-1} (U), \tilde{x}^i , g^i_j )$ be the
naturally induced local coordinates on $L(M)$, where $\Pi : L(M)
\to M$ is the projection. Then $g^i_j (f(t,s)) = X^i_j (a(t,s))$.
As $\omega$ is the connection $1$-form of a linear connection
\[ \omega = \omega^j_i \otimes E^i_j \]
where $\omega^i_j$ are scalar $1$-forms on $L(M)$ and $\{ E^i_j :
1 \leq i , j \leq 2n+1 \}$ is the basis of the Lie algebra ${\bf
gl}(2n+1)$ given by $E^i_j = [\delta^i_\ell \; \delta_j^k ]_{1
\leq k,\ell \leq 2n+1}$. Let $\{ e_1 , \cdots , e_{2n+1} \}$ be
the canonical linear basis of ${\mathbb R}^{2n+1}$. Then
\[ \left. \mu^* ({\mathbb T})_{(t,0)} = f(t,0)^{-1} \dot{\gamma}(t) = \frac{d x^i}{d t}
f(t,0)^{-1} \frac{\partial}{\partial x^i} \right|_{\gamma (t)} =
\frac{d x^i}{d t} Y^j_i e_j \] where $[Y^i_j ] = [X^i_j ]^{-1}$.
Therefore
\[ \omega^* ({\mathbb S})_{(t,0)} \cdot \mu^* ({\mathbb
T})_{(t,0)} = \frac{d x^k}{d t} \, Y^i_k \, (f^* \omega^j_i
)({\mathbb S})_{(t,0)}  \, e_j \] (because of $E^i_j \, e_k =
\delta^i_k \, e_j$). On the other hand (by Prop. 1.1 in
\cite{kn:KoNo}, Vol. I, p. 64) $\omega^* ({\mathbb S})_{(t,0)} =
A$ where the left invariant vector field $A \in {\bf gl}(2n+1)$ is
given by
\begin{equation} \label{e:definv}
A^*_{f(t,0)} = (d_{(t,0)} f) {\mathbb S}_{(t,0)} - \ell_{f(t,0)}
\dot{\sigma}_t (0)
\end{equation}
and $\ell_u : T_{\Pi (u)}(M) \to H_u$ is the inverse of $d_u \Pi :
H_u \to T_{\Pi (u)} (M)$, $u \in L(M)$ (the horizontal lift
operator with respect to $H$). Here $A^*$ is the fundamental
vector field associated to $A$, i.e.
\[ A^*_{f(t,0)} = (d_e L_{f(t,0)} ) A_e \]
where $L_u : {\rm GL}(2n+1) \to L(M)$, $u \in L(M)$, is given by
$L_u (g) = u g$ for any $g \in {\rm GL}(2n+1)$, and $e \in {\rm
GL}(2n+1)$ is the unit matrix. If $A = A^j_i E^i_j$ then $A^i_j =
(f^* \omega^i_j )({\mathbb S})_{(t,0)}$. Let $(g^i_j )$ be the
natural coordinates on ${\rm GL}(2n+1)$ so that $L_{f(t,0)}$ is
locally given by
\[ L^i (g) = \tilde{x}^i \, , \;\;\; L^i_j (g) = X^i_k g^k_j \, , \] and then
$(d_e L_{f(t,0)}) (\partial /\partial g^i_j )_e =  X_i^k (\partial
/\partial g^k_j )_{f(t,0)}$. Next (cf. \cite{kn:KoNo}, Vol. I, p.
143)
\[ \ell \; \frac{\partial}{\partial x^j} = \partial_j -
(\Gamma^i_{jk} \circ \Pi ) g^k_\ell \frac{\partial}{\partial
g^i_\ell}
\] (where $\partial_i = \partial /\partial \tilde{x}^i$) and
(\ref{e:definv}) lead to
\[ \left. A^k_\ell X^i_k \frac{\partial}{\partial g^i_\ell} \right|_{f(t,0)}
= (d_{(t,0)} f) {\mathbb S}_{(t,0)} - X^j (\gamma (t)) \{
\partial_j - (\Gamma^i_{jk} \circ \Pi ) g^k_\ell  \frac{\partial}{\partial
g^i_\ell} \}_{f(t,0)} \] or (by applying this identity to the
coordinate functions $g^i_\ell$)
\begin{equation}
\label{e:64} A^k_\ell X^i_k = {\mathbb S}_{(t,0)} (g^i_\ell \circ
f) + X^j (\gamma (t)) \Gamma^i_{jk} (\gamma (t)) X^k_\ell \, .
\end{equation}
If $f^i_j = g^i_j \circ f$ then
\[ {\mathbb S}_{(t,0)} (f^i_j ) = \frac{d}{d s} \{ f^i_j \circ
\beta_t \} (0) = \frac{\partial f^i_j}{\partial s}(t,0) \]
Therefore (by (\ref{e:64}))
\begin{equation} \label{e:65} A^k_\ell = Y^k_i \{ \frac{\partial
f^i_\ell}{\partial s} (t,0) + X^j (\gamma (t))
\Gamma^i_{jm}(\gamma (t)) X^m_\ell \} \, .
\end{equation}
So far we got (by (\ref{e:65}))
\[ \omega^* ({\mathbb S})_{(t,0)} \cdot \mu^* ({\mathbb
T})_{(t,0)} = Y^\ell_k \; \frac{d x^k}{dt} \; \{ \frac{\partial
f^i_\ell}{\partial s}(t,0) + \] \[ + X^j (\gamma (t))
\Gamma^i_{jm}(\gamma (t)) X^m_\ell \} f(t,0)^{-1} \left.
\frac{\partial}{\partial x^i}\right|_{\gamma (t)} \, . \] Let us
observe that
\[ \frac{\partial f^i_j}{\partial s}(t,0) = \frac{\partial
X^i_j}{\partial x^k}(\gamma (t)) \frac{\partial a^k}{\partial
s}(t,0) = \frac{\partial X^i_j}{\partial x^k} (\gamma (t)) X^k
(\gamma (t)) \] hence
\[ \frac{\partial f^i_j}{\partial s}(t,0) +  X^k (\gamma (t))
\Gamma^i_{k\ell}(\gamma (t)) X^\ell_j = (\nabla_X X_j )^i_{\gamma
(t)} \] and we may conclude that
\begin{equation}
\omega^* ({\mathbb S})_{(t,0)} \cdot \mu^* ({\mathbb T})_{(t,0)} =
Y^j_k \; \frac{d x^k}{dt} \; f(t,0)^{-1} (\nabla_X X_j )_{\gamma
(t)} = 0 . \label{e:66}
\end{equation} Indeed
\[ (\nabla_X X_i )_{\gamma (t)} = (\nabla_{\dot{\sigma}_t} X_i
)_{\sigma_t (0)} = \lim_{s \to 0} \frac{1}{s} \{ \tau^s X_{i ,
\sigma_t (s)} - X_{i , \sigma_t (0)} \} = \]
\[ =   \lim_{s \to 0} \frac{1}{s} \{ \tau^s f(t,s) e_i - f(t,0)
e_i \} = 0 \] because $f_t$ is horizontal (yielding $\tau^s f(t,s)
= f(t,0)$). By (\ref{e:TS})-(\ref{e:curv}),
(\ref{e:int})-(\ref{e:stor}) and (\ref{e:66}) the identity
(\ref{e:N}) may be written
\[
\frac{d^2}{d s^2} \{ L(\gamma^s )\}_{s=0} = \frac{1}{r} \int_a^b
\{ |\nabla_{\dot{\gamma}} X^\bot |^2 - g_\theta (R(X^\bot ,
\dot{\gamma}) \dot{\gamma} , X^\bot ) + \]
\[ + 2 g_\theta (T_\nabla (X^\bot , \dot{\gamma}) ,
\nabla_{\dot{\gamma}} X ) + |T_\nabla (X^\bot , \dot{\gamma})|^2 -
\]
\[ - \frac{1}{r^2} [ g_\theta (T_\nabla (X^\bot , \dot{\gamma}) , \dot{\gamma})^2 +
2 \, g_\theta (\nabla_{\dot{\gamma}} X , \dot{\gamma}) g_\theta
(T_\nabla (X^\bot , \dot{\gamma}), \dot{\gamma})] \} d t  \] or
(by $T_\nabla (X^\bot  , \dot{\gamma}) = - 2 \Omega (X^\bot ,
\dot{\gamma}) T$ and $\theta (\dot{\gamma}) = 0$)
\begin{equation} \label{e:68} I(X,X) = \frac{1}{r} \int_a^b
\{ |\nabla_{\dot{\gamma}} X^\bot |^2 - g_\theta (R(X^\bot ,
\dot{\gamma}) \dot{\gamma} , X^\bot ) + \end{equation} \[ + 4
\Omega (X^\bot , \dot{\gamma})^2 - 4 \Omega (X^\bot ,
\dot{\gamma}) \theta (\nabla_{\dot{\gamma}} X) \} \; d t . \]
Finally, by polarization $I(X,Y) = \frac{1}{2} \{ I(X+Y,X+Y) -
I(X,X) - I(Y,Y)\}$ the identity (\ref{e:68}) leads to
(\ref{e:index}).
\par
{\em Proof of Theorem} \ref{t:10}. As $\nabla g_\theta = 0$
\[ \int_{t_j}^{t_{j+1}} \{ g_\theta (\nabla_{\dot{\gamma}} X^\bot
, \nabla_{\dot{\gamma}} Y^\bot ) - 2 \Omega (X^\bot ,
\dot{\gamma}) \theta (\nabla_{\dot{\gamma}} Y^\bot ) \} d t = \]
\[ = \int_{t_j}^{t_{j+1}} \{ \frac{d}{d t}[g_\theta
(\nabla_{\dot{\gamma}} X^\bot , Y^\bot ) - 2 \Omega (X^\bot ,
\dot{\gamma}) \theta (Y^\bot )] - \] \[ - g_\theta
(\nabla^2_{\dot{\gamma}} X^\bot - 2 \Omega (\nabla_{\dot{\gamma}}
X^\bot , \dot{\gamma}) T , Y^\bot ) \} d t =  \]
\[ = g_{\theta , \gamma (t_{j+1})} ((\nabla_{\dot{\gamma}} X^\bot
)^-_{\gamma (t_{j+1})} , Y^\bot_{\gamma (t_{j+1})}) - g_{\theta ,
\gamma (t_j )} ((\nabla_{\dot{\gamma}} X^\bot )^+_{\gamma (t_j )}
, Y^\bot_{\gamma (t_j )}) - \]
\[ - 2 \Omega (X^\bot , \dot{\gamma})_{\gamma (t_{j+1})} \theta
(Y^\bot )_{\gamma (t_{j+1})} + 2 \Omega (X^\bot ,
\dot{\gamma})_{\gamma (t_j )} \theta (Y^\bot )_{\gamma (t_j )} -
\]
\[ - \int_{t_j}^{t_{j+1}} \{ g_\theta
(\nabla^2_{\dot{\gamma}} X^\bot - 2 \Omega (\nabla_{\dot{\gamma}}
X^\bot , \dot{\gamma}) T , Y^\bot ) \} d t \] and (\ref{e:index})
implies (\ref{e:index2}). Q.e.d.
\par
{\em Proof of Corollary} \ref{c:4}. If $X^\bot \in J_\gamma$ then
$X^\bot$ is differentiable in $[a,b]$ hence the last term in
(\ref{e:index2}) vanishes. Also ${\mathcal J}_\gamma X^\bot = 0$
and (\ref{e:index2}) yield
\[ I(X,Y) = - \frac{2}{r} \int_a^b \{ \theta (\nabla_{\dot{\gamma}}
X^\bot ) - 2 \Omega (X^\bot , \dot{\gamma})\} \Omega (Y^\bot ,
\dot{\gamma}) \, d t \] which implies (by Lemma \ref{l:J2}) both
(\ref{e:lJ2})-(\ref{e:indlJ2}). Viceversa, let us assume that
(\ref{e:lJ2}) holds for some $\alpha (X) \in {\mathbb R}$. Let $f$
be a smooth function on $M$ such that $f(\gamma (t_j )) = 0$ for
any $0 \leq j \leq h$ and $f(\gamma (t)) > 0$ for any $t \in [a,b]
\setminus \{ t_0 , t_1 , \cdots , t_h \}$ and let us consider the
vector field $Y = f \, {\mathcal J}_\gamma X^\bot$. As ${\mathcal
J}_\gamma X^\bot$ is orthogonal to $\dot{\gamma}$ the identity
(\ref{e:indlJ2}) implies
\[ \int_a^b f(\gamma (t)) |{\mathcal J}_\gamma X^\bot |^2 \, d t =
0 \] hence ${\mathcal J}_\gamma X^\bot = 0$ in each interval $[t_j
, t_{j+1}]$. To prove that $X^\bot \in J_\gamma$ it suffices (by
Prop. 1.1 in \cite{kn:KoNo}, Vol. II, p. 63) to check that
$X^\bot$ is of class $C^1$ at each $t_j$. To this end, for each
fixed $j$ we consider a vector field $Y$ along $\gamma$ such that
\[ Y_{j , \gamma (t)} = \begin{cases} (\nabla_{\dot{\gamma}}
X^\bot )^-_{\gamma (t_j )} - (\nabla_{\dot{\gamma}} X^\bot
)^+_{\gamma (t_j )} \, , & {\rm for} \;\; t = t_j \cr 0, & {\rm
for} \;\;\; t = t_k \, , \;\; k \neq j. \cr \end{cases} \] Then
(by (\ref{e:indlJ2})) $|(\nabla_{\dot{\gamma}} X^\bot )^-_{\gamma
(t_j )} - (\nabla_{\dot{\gamma}} X^\bot )^+_{\gamma (t_j )}|^2 =
0$. Q.e.d. \vskip 0.1in {\bf Remark}. Let $\gamma (t) \in M$ be a
lengthy $C^1$ curve. Then $D_{\dot{\gamma}} \dot{\gamma} =
\nabla_{\dot{\gamma}} \dot{\gamma} - A(\dot{\gamma}, \dot{\gamma})
T$ hence on a Sasakian manifold $\gamma$ is a geodesic of $\nabla$
if and only if $\gamma$ is a geodesic of the Riemannian manifold
$(M , g_\theta )$. This observation leads to the following
alternative proof of Theorem \ref{t:9}. Let $\gamma \in \Gamma$
and $X,Y \in T_\gamma (\Gamma )$ be as in Theorem \ref{t:9}. By
Theorem 5.4 in \cite{kn:KoNo}, Vol. II, p. 81, we have
\begin{equation} I(X,Y) = \frac{1}{r} \int_a^b \{ g_\theta (D_{\dot{\gamma}}
X^\bot , D_{\dot{\gamma}} Y^\bot ) - g_\theta (R^D (X^\bot ,
\dot{\gamma}) \dot{\gamma} , Y^\bot )\} dt .
\label{e:rem0}
\end{equation}
Now on one hand
\begin{equation}
D_{\dot{\gamma}} X^\bot = \nabla_{\dot{\gamma}} X^\bot +
\Omega (\dot{\gamma} , X^\bot ) T + \theta (X^\bot ) J
\dot{\gamma}
\label{e:rem1}
\end{equation}
and on the other the identity
\[ R^D (X,Y) Z = R(X,Y) Z + (J X \wedge J Y) Z - 2 \Omega (X,Y) J
Z + \] \[ + 2 g_\theta ((\theta \wedge I)(X,Y) , Z) T - 2 \theta
(Z) (\theta \wedge I)(X,Y), \;\;\; X,Y,Z \in {\mathcal X}(M), \]
yields
\begin{equation}
R^D (X^\bot , \dot{\gamma})\dot{\gamma} = R(X^\bot ,
\dot{\gamma})\dot{\gamma} - 3 \Omega (X^\bot , \dot{\gamma}) J
\dot{\gamma} + r^2 \theta(X^\bot ) T. \label{e:rem2}
\end{equation}
Let us substitute from (\ref{e:rem1})-(\ref{e:rem2}) into
(\ref{e:rem0}) and use the identity
\[ \theta (X^\bot ) \Omega (\nabla_{\dot{\gamma}} Y^\bot ,
\dot{\gamma}) + \theta (Y^\bot ) \Omega (\nabla_{\dot{\gamma}}
X^\bot , \dot{\gamma}) + \]\[ + \Omega (\dot{\gamma} , X^\bot )
\theta (\nabla_{\dot{\gamma}} Y^\bot ) + \Omega (\dot{\gamma} ,
Y^\bot ) \theta (\nabla_{\dot{\gamma}} X^\bot ) = \]
\[ = \frac{d}{d t} \{ \theta (X^\bot ) \Omega (Y^\bot ,
\dot{\gamma}) + \theta (Y^\bot ) \Omega (X^\bot , \dot{\gamma}) \}
- \]
\[ - 2 \{ \Omega (X^\bot , \dot{\gamma}) \theta
(\nabla_{\dot{\gamma}} Y^\bot ) + \Omega (Y^\bot , \dot{\gamma})
\theta (\nabla_{\dot{\gamma}} X^\bot ) \} \] (together with
$X_{\gamma (a)} = X_{\gamma (b)} = 0$) so that to derive
(\ref{e:index}). Q.e.d.

\vskip 0.1in As an application of Theorems \ref{t:conj3} and
\ref{t:9} we shall establish
\begin{theorem} Let $(M , \theta )$ be a Sasakian manifold of CR dimension $n$
and $\nabla$ its Tanaka-Webster connection. Let $\gamma : [a,b]
\to M$ be a lengthy geodesic of $\nabla$, parametrized by arc
length. If there is $c \in (a,b)$ such that the point $\gamma (c)$
is horizontally conjugate to $\gamma (a)$ and for any $\delta > 0$
with $[c-\delta , c + \delta ] \subset (a,b)$ the space ${\mathcal
H}_{\gamma_\delta}$ has maximal dimension $4n$ {\rm (}where
$\gamma_\delta$ is the geodesic $\gamma : [c-\delta , c + \delta ]
\to M${\rm )} then $\gamma$ is not a minimizing geodesic joining
$\gamma (a)$ and $\gamma (b)$, that is the length of $\gamma$ is
greater than the Riemannian distance {\rm (}associated to $(M ,
g_\theta )${\rm )} between $\gamma (a)$ and $\gamma (b)$.
\label{t:14}
\end{theorem}
{\em Proof}. Let $\gamma : [a,b] \to M$ be a geodesic of the
Tanaka-Webster connection of the Sasakian manifold $(M , \theta
)$, obeying to the assumptions in Theorem \ref{t:14}. Then (by
Theorem \ref{t:conj3}) there is a piecewise differentiable vector
field $X$ along $\gamma$ such that 1) $X$ is orthogonal to
$\dot{\gamma}$ and $J \dot{\gamma}$, 2) $X_{\gamma (a)} =
X_{\gamma (b)} = 0$, and 3) $I_a^b (X) < 0$. Let $\{\gamma^s
\}_{|s| < \epsilon}$ be a $1$-parameter family of curves as in the
definition of $(d_\gamma L) X$ and $I(X,X)$. By Corollary
\ref{c:V1} (as $\gamma$ is a geodesic of $\nabla$) one has
\[ \frac{d}{d s} \left\{ L(\gamma^s )\right\}_{s=0} = 0. \]
On the other hand (by Theorem \ref{t:9} and $X^\bot = X$)
\[  I(X,X) = I_a^b (X) + 4 \int_a^b \Omega (X , \dot{\gamma}) \{ \Omega (X ,
\dot{\gamma}) - \theta (X^\prime )\} d t \] hence (as $X$ is
orthogonal to $J \dot{\gamma}$)
\[ \frac{d^2}{d s^2} \left\{ L(\gamma^s ) \right\}_{s=0} = I_a^b
(X) < 0 \] so that there is $0 < \delta <\epsilon$ such that
$L(\gamma^s ) < L(\gamma )$ for any $|s| < \delta$. \vskip 0.1in
{\bf Remark}. If there is a $1$-parameter variation of $\gamma$
(inducing $X$) by {\em lengthy} curves then $L(\gamma )$ is
greater than the Carnot-Carath\'edory distance between $\gamma
(a)$ and $\gamma (b)$.

\section{Final comments and open problems}
Manifest in R. Strichartz's paper (cf. \cite{kn:Str}) is the
absence of covariant derivatives and curvature. Motivated by our
Theorem \ref{t:2} we started developing a theory of geodesics of
the Tanaka-Webster connection $\nabla$ on a Sasakian manifold $M$,
with the hope that although lengthy geodesics of $\nabla$ form
(according to Corollary \ref{c:relation}) a smaller family than
that of sub-Riemannian geodesics, the former may suffice for
establishing an analog to Theorem 7.1 in \cite{kn:Str}, under the
assumption that $\nabla$ is complete (as a linear connection on
$M$). The advantage of working within the theory of linear
connections is already quite obvious (e.g. any $C^1$ geodesic of
$\nabla$ is automatically of class $C^\infty$, as an integral
curve of some $C^\infty$ basic vector field, while sub-Riemannian
geodesics are assumed to be of class $C^2$, cf. \cite{kn:Str}, p.
233, and no further regularity is to be expected {\em a priori})
and doesn't contradict R. Strichartz's observation that
sub-Riemannian manifolds, and in particular strictly pseudoconvex
CR manifolds endowed with a contact form $\theta$, exhibit no
approximate Euclidean behavior (cf. \cite{kn:Str}, p. 223).
Indeed, while Riemannian curvature measures the higher order
deviation of the given Riemannian manifold from the Euclidean
model, the curvature of the Tanaka-Webster connection describes
the pseudoconvexity properties of the given CR manifold, as
understood in several complex variables analysis. The role as a
possible model space played by the {\em tangent cone} of the
metric space $(M , \rho )$ at a point $x \in M$ (such as produced
by J. Mitchell's Theorem 1 in \cite{kn:Mic}, p. 36) is unclear.
\par
Another advantage of our approach stems from the fact that the
exponential map on $M$ thought of as a sub-Riemannian manifold is
never a diffeomorphism at the origin (because all sub-Riemannian
geodesics issuing at $x \in M$ must have tangent vectors in
$H(M)_x$) in contrast with the ordinary exponential map associated
to the Tanaka-Webster connection $\nabla$. In particular {\em cut
points} (as introduced in \cite{kn:Str}, p. 260) do not possess
the properties enjoyed by conjugate points in Riemannian geometry
because (by Theorem 11.3 in \cite{kn:Str}, p. 260) given $x \in M$
cut points occur arbitrary close to $x$. On the contrary (by
Theorem 1.4 in \cite{kn:KoNo}, Vol. II, p. 67) given $x \in M$ one
may speak about the {\em first} point conjugate to $x$ along a
geodesic of $\nabla$ emanating from $x$, therefore the concept of
conjugate locus $C(x)$ may be defined in the usual way (cf. e.g.
\cite{kn:Man}, p. 117). The systematic study of the properties of
$C(x)$ on a strictly pseudoconvex CR manifold is an open problem.
\par
Yet another concept of exponential map was introduced by D.
Jerison \& J. M. Lee, \cite{kn:JeLe} (associated to {\em parabolic
geodesics} i.e. solutions $\gamma (t)$ to $\left(
\nabla_{\dot{\gamma}} \dot{\gamma} \right)_{\gamma (t)} = 2c
T_{\gamma (t)}$ for some $c \in {\mathbb R}$). A comparison
between the three exponential formalisms (in \cite{kn:Str},
\cite{kn:JeLe}, and the present paper) hasn't been done as yet. We
conjecture that given a $2$-plane $\sigma \subset T_x (M)$ its
pseudohermitian sectional curvature $k_\theta (\sigma )$ measures
the difference between the length of a circle in $\sigma$ (with
respect to $g_{\theta , x}$) and the length of its image by
$\exp_x$ (the exponential mapping at $x$ associated to $\nabla$).
Also a useful relationship among $\exp_x$ and the exponential
mapping associated to the Fefferman metric $F_\theta$ on $C(M)$
should exist (and then an understanding of the singular points of
the latter, cf. e.g. M. A. Javaloyes \& P. Piccione,
\cite{kn:JaPi}, should shed light on the properties of singular
points of the former).
\par
Finally, the analogy between Theorem 7.3 in \cite{kn:Str}, p. 245
(producing ``approximations to unity'' on Carnot-Carath\'eodory
complete sub-Riemannian manifolds) and Lemma 2.2 in
\cite{kn:Str2}, p. 50 (itself a corrected version of a result by
S.-T. Yau, \cite{kn:Yau}) indicates that Theorem 7.3 is the proper
ingredient for proving that the sublaplacian $\Delta_b$ is
essentially self-adjoint on $C^\infty_0 (M)$ and the corresponding
heat operator is given by a positive $C^\infty$ kernel. These
matters are relegated to a further paper.

\begin{appendix}
\section{Contact forms of constant pseudohermitian sectional curvature}
The scope of this section is to give a proof of Theorem
\ref{t:J3}. Let $(M , \theta )$ be a nondegenerate CR manifold and
$\theta$ a contact form on $M$. Let $\nabla$ be the Tanaka-Webster
connection of $(M , \theta )$. We recall the first Bianchi
identity
\begin{equation} \sum_{XYZ} R(X,Y)Z = \sum_{XYZ} \{ T_\nabla
(T_\nabla (X,Y), Z) + (\nabla_X T_\nabla )(Y,Z) \} \label{e:A1}
\end{equation}
for any $X,Y,Z \in T(M)$, where $\sum_{XYZ}$ denotes the cyclic
sum over $X,Y,Z$. Let $X,Y, Z \in H(M)$ and note that
\[ T_\nabla (T_\nabla (X,Y) , Z) = - 2 \Omega (X,Y) \tau (Z), \]
\[ (\nabla_X T_\nabla )(Y,Z) = - 2 (\nabla_X \Omega )(Y,Z) T = 0.
\]
Indeed $\nabla g_\theta = 0$ and $\nabla J = 0$ yield $\nabla
\Omega = 0$. Thus (\ref{e:A1}) leads to
\begin{equation}
\sum_{XYZ} R(X,Y)Z = - 2 \sum_{XYZ} \Omega (X,Y) \tau (Z),
\label{e:A2}
\end{equation}
for any $X,Y,Z \in H(M)$. Let us define a $(1,2)$-tensor field $S$
by setting $S(X,Y) = (\nabla_X \tau )Y - (\nabla_Y \tau )X$. Next,
we set $X,Y \in H(M)$ and $Z = T$ in (\ref{e:A1}) and observe that
\[ T_\nabla (T_\nabla (X,Y), T) + T_\nabla (T_\nabla (Y,T), X) +
T_\nabla (T_\nabla (T,X), Y) = \]
\[ = - T_\nabla (\tau (Y), X) + T_\nabla (\tau (X), Y) = \;\;\;\;\;\; {\rm
(as \; \tau \; is \; H(M)-valued)} \]
\[ = 2 \{ \Omega (\tau (Y) , X) - \Omega (\tau (X) , Y) \} T =
2 g_\theta ((\tau J + J \tau )X , Y) T = 0, \] (by the purity
axiom) and
\[ (\nabla_X T_\nabla )(Y, T) + (\nabla_Y T_\nabla )(T , X) +
(\nabla_T T_\nabla )(X,Y) = \]
\[ = - (\nabla_X \tau )Y + (\nabla_Y \tau )X - 2 (\nabla_T \Omega
)(X,Y) T = - S(X,Y). \] Finally (\ref{e:A1}) becomes
\begin{equation}
R(X,T) Y + R(T, Y) X = S(X,Y), \label{e:A3}
\end{equation}
for any $X,Y \in H(M)$. The $4$-tensor $R$ enjoys the properties
\begin{equation}
R(X,Y,Z,W) = - R(Y,X, Z, W), \label{e:A4}
\end{equation}
\begin{equation}
R(X,Y, Z,W) = - R(X,Y, W, Z), \label{e:A5}
\end{equation}
for any $X,Y,Z,W \in T(M)$. Indeed (\ref{e:A4}) follows from
$\nabla g_\theta = 0$ while (\ref{e:A5}) is obvious. We may use
the reformulation (\ref{e:A2})-(\ref{e:A3}) of the first Bianchi
identity to compute $\sum_{YZW} R(X,Y,Z,W)$ for arbitrary vector
fields. For any $X \in T(M)$ we set $X_H = X - \theta (X) T$ (so
that $X_H \in H(M)$). Then
\[ \sum_{YZW} R(X,Y,Z,W) = \sum_{YZW} g_\theta (R(Z,W) Y_H , X) =
\]
\[ = \sum_{YZW} g_\theta (R(Z_H , W_H ) Y_H  +  \theta (Y)
[R(W_H , T) Z_H + R(T, Z_H ) W_H ] , X)   \] hence
\begin{equation}
\sum_{YZW} R(X,Y,Z,W) = \label{e:A6}
\end{equation}
\[ = - \sum_{YZW} \{ 2 \Omega (Y,Z) A(W, X) +
\theta (Y) g_\theta (X , S(Z_H , W_H )) \} \] for any $X,Y,Z,W \in
T(M)$. Next, we set \[ K(X,Y,Z,W) = \sum_{YZW} R(X,Y,Z,W) \] and
compute (by (\ref{e:A4})-(\ref{e:A5}))
\[ K(X,Y,Z,W) - K(Y,Z,W,X) - K(Z,W,X,Y) + K(W,X,Y,Z) = \]
\[ = 2 R(X,Y,Z,W) - 2 R(Z,W,X,Y) \]
hence (by (\ref{e:A6}))
\[ 2 R(X,Y,Z,W) - 2 R(Z,W,X,Y) = \] \[ = - \sum_{YZW} \{ 2 \Omega (Y,Z)
A(X,W) + \theta (Y) g_\theta (X , S(Z_H , W_H )) \} + \]
\[ + \sum_{ZWX} \{ 2 \Omega (Z,W) A(Y,X) + \theta (Z) g_\theta (Y,
S(W_H , X_H )) \} + \]
\[ + \sum_{WXY} \{ 2 \Omega (W,X) A(Y,Z) + \theta (W) g_\theta (Z
, S(X_H , Y_H )) \} - \]
\[ - \sum_{XYZ} \{ 2 \Omega (X,Y) A(Z,W) + \theta (X) g_\theta (W,
S(Y_H , Z_H )) \} \] or
\begin{equation} 2 R(X,Y,Z,W) - 2 R(Z,W,X,Y) = \label{e:A7}
\end{equation}
\[ = - 4 \Omega (Y,Z) A(X,W) + 4 \Omega (Y,W) A(X,Z) - \] \[ - 4 \Omega
(X,W) A(Y,Z) + 4 \Omega (X,Z) A(Y,W) + \]
\[ + \theta (X) [ g_\theta (Y, S(Z_H , W_H )) + g_\theta (Z, S(Y_H
, W_H )) - g_\theta (W, S(Y_H , Z_H )) ] + \]
\[ + \theta (Y) [ g_\theta (Z, S(W_H , X_H )) - g_\theta (W, S(Z_H
, X_H )) - g_\theta (X , S(Z_H , W_H )) ] + \]
\[ + \theta (Z) [g_\theta (Y, S(W_H , X_H )) - g_\theta (X, S(W_H
, Y_H )) - g_\theta (W, S(X_H , Y_H ))] + \]
\[ + \theta (W)[ g_\theta (Y, S(X_H , Z_H )) - g_\theta (X, S(Y_H
, Z_H )) + g_\theta (Z, S(X_H , Y_H ))]. \] As $\nabla_X \tau$ is
symmetric one has
\[ g_\theta (Y, S(X,Z)) - g_\theta (X, S(Y,Z)) = g_\theta
(S(X,Y),Z) \] for any $X,Y,Z \in H(M)$, so that (\ref{e:A7}) may
be written
\begin{equation} R(X,Y,Z,W) = R(Z,W,X,Y) -
\label{e:A8}
\end{equation}
\[ - 2 \Omega (Y,Z) A(X,W) + 2 \Omega (Y,W) A(X,Z) - \]
\[ - 2 \Omega (X,W) A(Y,Z) + 2 \Omega (X,Z) A(Y,W) + \]
\[ + \theta (X) g_\theta (S(Z_H , W_H ), Y) + \theta (Y) g_\theta
(S(W_H , Z_H ) , X) + \]
\[ + \theta (Z) g_\theta (S(Y_H , X_H ), W) + \theta (W) g_\theta
(S(X_H , Y_H ), Z) , \] for any $X,Y,Z,W \in T(M)$.
\par
The properties (\ref{e:A4})-(\ref{e:A6}) and (\ref{e:A8}) may be
used to compute the full curvature of a manifold of constant
pseudohermitian sectional curvature (the arguments are similar to
those in the proof of Prop. 1.2 in \cite{kn:KoNo}, Vol. I, p.
198). Assume from now on that $M$ is strictly pseudoconvex and
$G_\theta$ positive definite. Let us set
\[ R_1 (X,Y,Z,W) = g_\theta (X,Z) g_\theta (Y,W) - g_\theta (Y,Z)
g_\theta (W,X) \] so that
\begin{equation}
R_1 (X,Y,Z,W) = - R_1 (Y,X,Z,W), \label{e:A9}
\end{equation}
\begin{equation}
R_1 (X,Y,Z,W) = - R_1 (X,Y,W,Z), \label{e:A10}
\end{equation}
\begin{equation}
\sum_{YZW} R_1 (X,Y,Z,W) = 0. \label{e:A11}
\end{equation}
Assume from now on that $k_\theta = c =$ const. Let us set $L = R
- 4 c R_1$ and observe that
\begin{equation}
L(X,Y,X,Y) = 0 \label{e:A12}
\end{equation}
for any $X,Y \in T(M)$. Indeed, if $X, Y$ are linearly dependent
then (\ref{e:A12}) follows from the skew symmetry of $L$ in the
pairs $(X,Y)$ and $(Z,W)$, respectively. If $X , Y$ are
independent then let $\sigma \subset T_x (M)$ be the $2$-plane
spanned by $\{ X_x , Y_x \}$, $x \in M$. Then
\[ L(X,Y,X,Y)_x  = R(X,Y,X,Y)_x - 4 c R_1 (X,Y,X,Y)_x = \]
\[ = 4 k_\theta (\sigma ) [|X|^2 |Y|^2 - g_\theta (X,Y)^2 ]_x - 4 c R_1
(X,Y,X,Y)_x = 0. \] Next (by (\ref{e:A12}))
\[ 0 = L(X,Y+W,X,Y+W) = L(X,Y,X,W) + L(X,W,X,Y) \]
i.e.
\begin{equation}
L(X,Y,X,W) = - L(X,W,X,Y) \label{e:A13}
\end{equation}
for any $X,Y,W \in T(M)$. As well known (cf. e.g. Prop. 1.1 in
\cite{kn:KoNo}, Vol. I, p. 198) the properties
(\ref{e:A9})-(\ref{e:A11}) imply as well the symmetry property
\begin{equation}
R_1 (X,Y,Z,W) = R_1 (Z,W, X,Y). \label{e:A14}
\end{equation}
Therefore $L(X,Y,Z,W) - L(Z,W,X,Y) = R(X,Y,Z,W) - R(Z,W,X,Y)$
hence (by (\ref{e:A8}))
\begin{equation}
L(X,Y,Z,W) = L(Z,W,X,Y) + \label{e:A15}
\end{equation}
\[ + 2 \Omega (Y,W) A(X,Z) - 2 \Omega (Y,Z) A(X,W) + \]
\[ + 2 \Omega (X,Z) A(Y,W) - 2 \Omega (X,W) A(Y,Z) + \]
\[ + \theta (X) g_\theta (S(Z_H , W_H ), Y) + \theta (Y) g_\theta
(S(W_H , Z_H ) , X) + \]
\[ + \theta (Z) g_\theta (S(Y_H , X_H ), W) + \theta (W) g_\theta
(S(X_H , Y_H ), Z) . \] Applying (\ref{e:A15}) (to interchange the
pairs $(X,W)$ and $(X,Y)$) we get
\[ L(X,W,X,Y) = L(X,Y,X,W) + \]
\[ + 2 \Omega (W,Y) A(X,X) - 2 \Omega (W,X) A(X,Y) - 2 \Omega
(X,Y) A(W,X) + \]
\[ + \theta (X) g_\theta (S(W_H , Y_H ) , X) + \theta (Y) g_\theta
(S(X_H , W_H ), X) + \] \[ + \theta (W) g_\theta (S(Y_H , X_H ) ,
X)
\] hence (\ref{e:A13}) may be written
\begin{equation} L(X,Y,X,W) = \Omega (W,X) A(X,Y) +  \label{e:A16}
\end{equation}
\[ + \Omega (X,Y)
A(W,X) - \Omega (W,Y) A(X,X) - \]
\[ - \frac{1}{2} \{ \theta (X) g_\theta (S(W_H , Y_H ) , X) +
\theta (Y) g_\theta (S(X_H , W_H ) , X) + \]
\[ + \theta (W) g_\theta (S(Y_H , X_H ) , X) \} . \]
Consequently
\[ L(X+Z,Y,X+Z,W) = \Omega (W, X+Z) A(X+Z, Y)  + \]
\[ + \Omega (X+Z, Y) A(W, X+Z) - \Omega (W,Y) A(X+Z, X+Z) - \]
\[ - \frac{1}{2} \; g_\theta (X+Z \; , \; \theta (X+Z) S(W_H , Y_H ) + \]
\[ + \theta (Y) S(X_H + Z_H , W_H ) + \theta (W) S(Y_H , X_H + Z_H ))
\]
or (using (\ref{e:A16}) to calculate $L(X,Y,X,W)$ and
$L(Z,Y,Z,W)$)
\begin{equation}
L(X,Y,Z,W) + L(Z,Y,X,W) = \Omega (X,Y) A(W , Z) + \label{e:A17}
\end{equation}
\[ + \Omega (W,X) A(Z,Y) + \Omega (W,Z) A(X,Y) + \] \[ + \Omega
(Z,Y) A(W,X) - 2 \Omega (W,Y) A(X,Z) - \]
\[ - \frac{1}{2} \; g_\theta (X , \theta (Z) S(W_H , Y_H ) +
\theta (Y)S(Z_H , W_H ) + \theta (W) S(Y_H , Z_H )) - \]
\[ - \frac{1}{2} \; g_\theta (Z , \theta (X) S(W_H , Y_H ) +
\theta (Y) S(X_H , W_H ) + \theta (W) S(Y_H , X_H )). \] On the
other hand, by the skew symmetry of $L$ in the first pair of
arguments and by (\ref{e:A15}) (used to interchange the pairs
$(Y,Z)$ and $(X,W)$)
\[ L(Z,Y,X,W) = - L(Y,Z, X,W) = - L(X,W,Y,Z) + \]
\[ + 2 \Omega (Z,X) A(Y,W) - 2 \Omega (Z,W) A(Y,X) + \]
\[ + 2 \Omega (Y,W) A(Z,X) - 2 \Omega(Y,X)A(Z,W) - \]
\[ - \theta (Y) g_\theta (S(X_H , W_H ) , Z) - \theta (Z) g_\theta
(S(W_H , X_H ) , Y) - \]
\[ - \theta (X) g_\theta (S(Z_H , Y_H ) , W) - \theta (W) g_\theta
(S(Y_H , Z_H ) , X) \] so that (\ref{e:A17}) becomes
\begin{equation}
L(X,Y,Z,W) = L(X,W,Y,Z) + 2 \Omega (X,Z) A(Y,W) - \label{e:A18}
\end{equation}
\[  - \Omega (W,Z) A(X,Y) - \Omega (X,Y) A(Z,W) + \] \[ + \Omega (W,X) A(Z,Y) +
\Omega (Z,Y) A(W,X) + \]
\[ + \frac{1}{2}\; \theta (X) \{ g_\theta (S(Z_H , Y_H ), W) +
g_\theta (S(Z_H , W_H ) , Y) \} - \]
\[ - \frac{1}{2} \; \theta (Y) \{ g_\theta (S(Z_H , W_H ) , X) +
g_\theta (S(W_H , X_H ) , Z) \} + \]
\[ + \frac{1}{2} \; \theta (Z) \{ g_\theta (S(W_H , X_H ) , Y) +
g_\theta (S(Y_H , X_H ) , W) \} - \]
\[ - \frac{1}{2} \; \theta (W) \{ g_\theta (S(Y_H , X_H ) , Z) +
g_\theta (S(Z_H , Y_H ) , X) \} . \] By cyclic permutation of the
variables $Y,Z,W$ in (\ref{e:A18}) we obtain another identity of
the sort
\[ L(X,Y,Z,W) = L(X,Z,W,Y) - 2 \Omega (X,W) A(Z,Y) + \]
\[ + \Omega (Y,W) A(X,Z) + \Omega (X,Z) A(W,Y) - \]
\[ - \Omega (Y,X) A(W,Z) - \Omega (W,Z) A(Y,X) - \]
\[ - \frac{1}{2} \; \theta (X) \{ g_\theta (S(W_H , Z_H ), Y) +
g_\theta (S(W_H , Y_H ) , Z) \} + \]
\[ + \frac{1}{2} \; \theta (Y) \{ g_\theta (S(Z_H , X_H ) , W) +
g_\theta ( S(W_H , Z_H ) , X)\} - \]
\[ + \frac{1}{2} \; \theta (Z) \{ g_\theta (S(W_H , Y_H ) , X) +
g_\theta (S(Y_H , X_H ) , W) \} + \]
\[ - \frac{1}{2} \; \theta (W) \{ g_\theta (S(Y_H, X_H ) , Z) +
g_\theta (S(Z_H , X_H ) , Y) \} \] which together with
\eqref{e:A18} leads to
\[ 3 L(X,Y,Z,W) = \sum_{YZW} L(X,Y,Z,W) - 2 \Omega (W,Z) A(X,Y) + \]
\[ + 3 \Omega (X,Z) A(Y,W) - 3 \Omega (X,W) A(Y,Z) +  \]
\[ + \Omega (Z,Y) A(W,X) + \Omega (Y, W) A(X,Z) +  \]
\[ + \frac{3}{2} \; \theta (X) g_\theta (S(Z_H , W_H ) , Y) - \frac{1}{2} \;
\theta (Y) g_\theta (S(Z_H , W_H ) , X) + \]
\[ + \frac{1}{2} \; \theta (Z) \{ 2 g_\theta (S(Y_H , X_H ) , W) +
\]
\[ + g_\theta (S(W_H , X_H ) , Y) + g_\theta (S( W_H , Y_H ) , X)
\} - \]
\[ - \frac{1}{2} \; \theta (W) \{ 2 g_\theta (S(Y_H , X_H ) , Z) + \]
\[ + g_\theta (S(Z_H , Y_H ) , X) + g_\theta (S(Z_H , X_H ) , Y ) \}  \]
or
\[ L(X,Y,Z,W) = \Omega (Y,W) A(X,Z) - \Omega (Y,Z) A(X,W) + \]
\[ + \Omega (X,Z) A(Y,W) - \Omega (X,W) A(Y,Z) + \]
\[ + \frac{1}{2} \; \{ \theta (X) g_\theta (S(Z_H , W_H ) , Y) - \theta (Y)
g_\theta (S(Z_H , W_H ) , X) +  \]
\[ +  \theta (Z) g_\theta (S(Y_H , X_H ) , W) -
\theta (W) g_\theta (S(Y_H , X_H , Z) \}  \] or
\begin{equation} R(X,Y,Z,W) = 4c \{ g_\theta (X,Z) g_\theta (Y,W) -
g_\theta (Y,Z) g_\theta (X,W) \} + \label{e:A19}
\end{equation}
\[ + \Omega (Y,W) A(X,Z) - \Omega (Y,Z) A(X,W) + \]
\[ + \Omega (X,Z) A(Y,W) - \Omega (X,W) A(Y,Z) + \]
\[ + g_\theta (S(Z_H , W_H ) , (\theta \wedge I)(X,Y)) -  g_\theta (S(X_H ,
Y_H ), (\theta \wedge I)(Z,W)) \] for any $X,Y,Z, W \in T(M)$,
where $I$ is the identical transformation and $(\theta \wedge
I)(X,Y) = \frac{1}{2} \{ \theta (X) Y - \theta (Y) X \}$. Using
(\ref{e:A19}) one may prove Theorem \ref{t:J3} as follows. Let $Y
= T$ in (\ref{e:A19}). As $R(Z,W) T = 0$ and $S$ is $H(M)$-valued
we get
\begin{equation} 0 = 4c \{ g_\theta (X,Z) \theta (W) -
g_\theta (X,W) \theta (Z) \} - \frac{1}{2} \; g_\theta (S(Z_H ,
W_H ) , X), \label{e:A20}
\end{equation}
for any $X,Z,W \in T(M)$. In particular for $Z,W \in H(M)$
\[ S(Z,W) = 0. \]
Hence $S(Z_H , W_H ) = 0$ and (\ref{e:A20}) becomes
\[ c \{ g_\theta (X,Z) \theta (W) - g_\theta (X,W) \theta (Z) \}
= 0, \] for any $X,Z,W \in T(M)$. In particular for $Z = X \in
H(M)$ and $W = T$ one has $c |X|^2 = 0$ hence $c = 0$ and
(\ref{e:A19}) leads to (\ref{e:J7}). Then $\tau = 0$ yields $R =
0$. To prove the last statement in Theorem \ref{t:J3} let us
assume that $M$ has CR dimension $n \geq 2$ (so that the Levi
distribution has rank $> 3$). Assume that $R = 0$ i.e.
\[ \Omega (X,Z) \tau (Y) - \Omega (Y,Z) \tau (X) = A(X,Z) J Y -
A(Y,Z) J X \] (by (\ref{e:J7})). In particular for $Z = Y$
\begin{equation}
\Omega (X,Y) \tau (Y) = A(X,Y) J Y - A(Y,Y) J X. \label{e:A21}
\end{equation}
Let $X \in H(M)$ such that $|X| = 1$, $g_\theta (X,Y) = 0$ and
$g_\theta (X, J Y) = 0$. Taking the inner product of (\ref{e:A21})
with $J X$ gives $A(Y,Y) = 0$, hence $A = 0$ (as $A$ is
symmetric). Q.e.d.
\end{appendix}

\end{document}